\documentclass[11pt]{amsart}
\usepackage{mathrsfs,amssymb,amsmath,url,a4wide}

\theoremstyle{plain}

    \newtheorem{thm}{Theorem}

    \newtheorem{lem}[thm]{Lemma}

    \newtheorem{prop}[thm]{Proposition}

    \newtheorem{cor}[thm]{Corollary}

\theoremstyle{definition}

    \newtheorem{defn}[thm]{Definition}
    \newtheorem{definition}[thm]{Definition}

\theoremstyle{remark}

\newcommand{\D}{{\mathscr D}}

\newcommand{\I}{{\mathscr I}}
\newcommand{\A}{{\mathscr A}}

\renewcommand{\S}{{\mathscr S}}

\renewcommand{\L}{{\mathfrak L}}

\newcommand{\ignore}[1]{}

\newcommand{\nin}{\notin}

\newcommand{\cl}[1]{\langle #1 \rangle}

\newcommand{\rev}{\leftrightarrow}

\newcommand{\To}{\rightarrow}
\newcommand{\mult}{\times}

\newcommand{\mix}[2]{\binom{#1}{#2}}
\newcommand{\turn}{\circlearrowleft}

\DeclareMathOperator{\id}{id} \DeclareMathOperator{\Aut}{Aut}

 \DeclareMathOperator{\tp}{tp}
 \DeclareMathOperator{\cyc}{Cycl}
\DeclareMathOperator{\sw}{sw}\DeclareMathOperator{\sep}{Sep}\DeclareMathOperator{\betw}{Betw}

\usepackage{graphicx}

\newcommand{\Fresse}{Fra\"{i}ss\'{e}}
\renewcommand{\D}{D}
\DeclareMathOperator{\tsw}{tsw}
\DeclareMathOperator{\Betw}{Betw}
\DeclareMathOperator{\Cycl}{Cycl}

\newcommand{\adel}{\Aut(\Delta)}
\newcommand{\adelcn}{\Aut(\Delta,c_1,\ldots,c_n)}
\newcommand{\delcn}{(\Delta,c_1,\ldots,c_n)}

\DeclareMathOperator{\Sep}{Sep}
\DeclareMathOperator{\JI}{JI}
\DeclareMathOperator{\JIo}{\JI^{o}}
\DeclareMathOperator{\OP}{OP}

\DeclareMathOperator{\dcl}{dcl}

\newcommand{\apl}{\mix l\id }
\newcommand{\apu}{\mix u\id}

\title{The 42 reducts of the random ordered graph}
\author
{Manuel Bodirsky}
    \address{Institut f\"{u}r Algebra\\TU Dresden\\01062 Dresden\\Germany}
    \email{Manuel.Bodirsky@tu-dresden.de}
   \urladdr{http://www.math.tu-dresden.de/~bodirsky/}
\thanks{The first and third author have received funding from the European Research Council under the European Community's Seventh Framework Programme (FP7/2007-2013 Grant Agreement no. 257039).}
\author{Michael Pinsker}\thanks{The second author has been funded through project I836-N23 of the  Austrian Science Fund (FWF) as well as through an APART fellowship of the Austrian Academy of Science.}
    \address{Institut f\"{u}r Computersprachen\\Theory and Logic Group\\Technische Universit\"{a}t Wien\\Favoritenstrasse 9/E1852\\A-1040 Wien\\
Austria}
    \email{marula@gmx.at}
    \urladdr{http://dmg.tuwien.ac.at/pinsker/}
\author{Andr\'{a}s Pongr\'{a}cz}
    \address{Laboratoire d'Informatique  (LIX), CNRS UMR 7161\\
    \'{E}cole Polytechnique \\91128 Palaiseau\\
    France}
\email{pongeee@cs.elte.hu}
\urladdr{http://www.lix.polytechnique.fr/~andras.pongracz/}

\date{\today}


\begin{document}
\maketitle

\begin{abstract}
The \emph{random ordered graph} is the up to isomorphism unique countable homogeneous linearly ordered graph that embeds all finite linearly ordered graphs. We determine the reducts of the random ordered graph up to first-order interdefinability.

\end{abstract}

\section{Introduction}
A famous result on infinite permutation groups, due to Cameron~\cite{Cameron5}, says
that there are %up to isomorphism 
exactly five closed permutation groups on a countably infinite set $D$ that are \emph{highly set-transitive}, that is, for all finite $A,B \subseteq D$ of the same cardinality 
the groups contain a permutation sending $A$ to $B$. This result is also famous in model theory, albeit in a different form: a relational structure $\Gamma$
is a \emph{reduct} of a structure $\Delta$ iff they
have the same domain, and every relation of $\Gamma$ has a first-order definition in $\Delta$ without parameters. 
For two structures $\Gamma_1,\Gamma_2$, we
set $\Gamma_1 \preceq \Gamma_2$ iff $\Gamma_1$ is a reduct of $\Gamma_2$; this defines a quasi-order
on the set of all structures. We call
two structures $\Gamma_1,\Gamma_2$ \emph{(first-order) equivalent} iff $\Gamma_1 \preceq \Gamma_2$ and $\Gamma_2 \preceq \Gamma_1$.
It can be shown (see Section 3.4 in~\cite{Oligo}) 
that Cameron's result on highly set-transitive permutation groups is equivalent to the fact 
 that the order of the rationals $({\mathbb Q};<)$ has precisely five reducts up to equivalence. Interestingly, many other prominent 
 homogeneous structures (see e.g.~\cite{Hodges}) have only finitely many inequivalent reducts: for instance the random graph~\cite{RandomReducts}, %(there are again exactly five reducts), 
 the Henson graphs~\cite{RandomReducts}, 
the random tournament~\cite{Bennett-thesis}, the expansion $({\mathbb Q};<,0)$ of $({\mathbb Q};<)$ by a constant~\cite{JunkerZiegler}, 
and the random partial order~\cite{Poset-Reducts}.

Thomas~\cite{RandomReducts} conjectured that \emph{every}
countable homogeneous structure with a finite relational
language has only finitely many inequivalent reducts. 
We know very little about this conjecture, beyond the fact that it is true for some fundamental 
homogeneous structures. When we factor the quasi-order $\preceq$ of the reducts of such a homogeneous structure by equivalence, then the resulting order forms a lattice; we do not even know
how to show that this lattice has only finitely many atoms, or no infinite ascending chains. 
In order to learn more, it seems to be
unavoidable to verify the conjecture 
for more of the classical
structures from model theory, independently from whether or not we believe the conjecture.
Let us mention that apart from the intrinsic 
interest of this problem at the intersection of model-theory and permutation group theory, 
there are other contexts in which
the conjecture plays an interesting role:
complexity classification in infinite-domain constraint satisfaction (see e.g.~\cite{BP-reductsRamsey}), 
and minimal flows
and universal minimal flows
in topological dynamics (see Section~4.2 in~\cite{TodorHDR}).

There is an important construction to systematically produce
new examples of homogeneous structures from old ones:
suppose that $\Gamma_1$ and $\Gamma_2$ are
two homogeneous structures 
with disjoint relational signatures $\tau_1$ and $\tau_2$
such that for each $i \in \{1,2\}$
the class of 
all finite structures that embed into $\Gamma_i$ has the
strong amalgamation property (see e.g.~\cite{Oligo}). 
%\footnote{The class of all finite
%substructures that embed into a homogeneous structure $\Gamma$ has the strong amalgamation property if and only if
%the algebraic closure of $\Gamma$ is trivial; these concepts will not be important in the further course of this article, so for further background we refer to~\cite{Oligo}.} (such as, for example, $({\mathbb Q};<)$ and the random graph);
%\footnote{A homogeneous structure has trivial algebraic closure if and only if its age has the strong amalgamation property; these concepts will not be important in the further course of this article, so for further background we refer to~\cite{Oligo}.}), 
Then the class of all finite structures whose $\tau_1$-reduct\footnote{The $\tau$-reduct of a $(\tau \cup \sigma)$-structure $\Gamma$ is the reduct of $\Gamma$ with signature $\tau$ obtained from $\Gamma$ by dropping the relations for the symbols from $\sigma$ and keeping the relations for the symbols from $\tau$.} embeds into $\Gamma_1$ and whose $\tau_2$-reduct embeds into $\Gamma_2$ is a strong amalgamation class. It can be verified by a straightforward back-and-forth argument that 
the $\tau_i$-reduct of the
\Fresse-limit  $\Delta$ of this class is isomorphic to $\Gamma_i$ (for an introduction to these basic concepts from model theory, see e.g.~\cite{Hodges}). We therefore call $\Delta$ the \emph{(free) superposition} of $\Gamma_1$ and $\Gamma_2$.
%Many important structures arise in this way from some base structures. 
The example we are concerned with in this paper is the
\emph{random ordered graph}, which is the free superposition of $({\mathbb Q};<)$ and the random graph, and denoted by $(\D;<,E)$. 
By the above, the reduct $(D;<)$ of the random ordered graph is isomorphic to the order of the rationals, 
and the reduct $(D;E)$ to the random graph. 
We can
thus interpret the construction of $(D;<,E)$ as one where we freely add
the order of the rationals to the random graph $(D;E)$, and consequently
could have called $(D;<,E)$ the \emph{ordered random graph}.
%As we have already noted, obvious reducts of the random ordered graph are 
%the reducts of the random graph and the reducts of $({\mathbb Q};<)$. 
%So we could as well have called $(\D;<,E)$ the \emph{ordered random graph}.

In this article we present a complete classification of the
reducts of the random ordered graph up to equivalence. Without counting the obvious reducts $(D;<,E)$ and $(D;=)$, there are precisely 42 such reducts (also see~\cite{HG2G}). This is the first time that the reducts of a free superposition of two homogeneous structures are determined.

There are several reasons why the random ordered 
graph is an interesting next candidate 
for verifying Thomas' conjecture. 
First of all, when $\Delta$ is the superposition
of two homogeneous structures $\Gamma_1$ and $\Gamma_2$, and when we know the reducts
of $\Gamma_1$ and $\Gamma_2$, we would
like to be able to compute the reducts of $\Delta$
from this information. The superposition of
$({\mathbb Q};<)$ and the random graph is one of the simplest examples of this type. While there are some reducts of $(\D;<,E)$
that are indeed built from the reducts of $({\mathbb Q};<)$ and 
the reducts of the random graph in a systematic way, additional 
reducts appear in our classification. 
%So our result shows that we are still far from a general reduct classification result for superpositions.
Moreover, there is a surprising asymmetry between the roles of $({\mathbb Q};<)$ and
$(V;E)$ in the reduct lattice of the random ordered graph. It is
probably fair to say that this shows that there is no general reduct classification result for superpositions. We remark that while we could have investigated the superposition of $({\mathbb Q};<)$ with itself or the random graph with itself, this would not have allowed us to discover the kind of asymmetry which we encounter in the present result. The reducts of the superposition of two copies of $({\mathbb Q};<)$ have been determined in subsequent work~\cite{LP-permutations-on-the}.

Another reason for studying the reducts of the random ordered graph is the technical dimension of this endeavor:
when it comes to reduct classification,
the complexity of this structure
%is
%in a different order of magnitude than the complexity of the previous structures for which the reducts have been classified with our method, and 
challenges the current methods. 
And indeed, not only the classifications of the reducts 
of $({\mathbb Q};<)$
and the random graph, but also the classification
of the reducts of the random 
tournament~\cite{Bennett-thesis} 
appears as a subcase of our
result, since the random tournament itself
turns out to be a reduct of the random ordered graph.

In our classification proof, we apply a technique 
that is based on \emph{Ramsey theory} and \emph{canonical functions}; an introduction to this technique can be
found in~\cite{BP-reductsRamsey}. 
In the setting of the random ordered graph,
the Ramsey-type result that we are going to
employ is the well-known fact that
the class of all finite linearly ordered graphs forms a Ramsey class (due to~\cite{AbramsonHarrington,NesetrilRoedlOrderedStructures}; also see~\cite{NesetrilRoedlPartite}). 
%While the technique applied by Junker and Ziegler to show that $({\mathbb Q};<,0)$ has 116 reducts 
%is quite specific to this structure, 
The strategy we apply here is generally applicable 
when the homogeneous structure under consideration has a expansion by finitely many relations which is homogeneous and whose class of finite induced substructures is a Ramsey class. 
We do not know of a single structure which is homogeneous in a finite relational signature and which does \emph{not} have such an expansion. Therefore, 
arbitrary homogeneous structures in a finite relational language might well be within
the scope of this method. 

When applying the method in practice for the random ordered graph one faces a severe
combinatorial explosion compared to previous classifications. Even though in principle
parts of the classification could have been automated, the sheer number of cases is prohibitive for such a computer-aided approach. An important idea to keep 
the description of the reducts manageable is to identify and work 
with the \emph{join-irreducible elements} of the reduct lattice.
In previous situations where the same method has been applied, 
the size of the reduct lattice was so
small that neither the advantage nor the possibility of reducing the
work to the join irreducibles became visible.

\vspace{.4cm}
\noindent {\bf Outline of the article.}
Section 2 is a catalog %presents the catalog 
of the reducts of the
random ordered graph.
Section 3 gives an overview of the classification proof. 
Section 4 shows how to reduce the 
classification to those reducts that contain the relation $<$. 
Section 5 then classifies those. 
In Section 6 we verify that the reduct lattice 
is indeed as it has been drawn in Section 2. 

%Coupling
%Reducts of the random tournament
%Special reducts -- describe only $S_4$ reduct
%Show that there are 42
%hile the technique applied here is generally applicable, the size of the reduct classification proof
%has a different order of magnitude than previous
%classifications. 

%do not include:
%open: having finitely many reducts preserved by interdefinability
%assumption on homogeneous in finite signature not really robust
%how about assumption: omega-cat in finite signature?
\section{The Reduct Catalog}
In this section we describe the reducts of the random ordered graph $(\D;<,E)$. 
%As we have mentioned in the introduction, we identify reducts $\Gamma_1$ and $\Gamma_2$ when
%$\Gamma_1 \preceq \Gamma_2$ and $\Gamma_2 \preceq \Gamma_1$, i.e., if they are reducts of each other.
%We write 
%$\Gamma_1 \prec \Gamma_2$, and say that $\Gamma_1$ is a \emph{proper} reduct of $\Gamma_2$, if $\Gamma_1 \preceq \Gamma_2$ but not $\Gamma_2 \preceq \Gamma_1$. 
We write $\Aut(\Gamma)$ for the set of automorphisms of a structure $\Gamma$.
For reducts $\Gamma_1$, $\Gamma_2$
of the random ordered graph we have
$\Gamma_1 \preceq \Gamma_2$ 
if and only if $\Aut(\Gamma_1) \supseteq \Aut(\Gamma_2)$.
This holds more generally for $\omega$-categorical structures, and is a consequence of the theorem of Engeler, Svenonius, and Ryll-Nardzewski (see, e.g.,~\cite{Hodges}). 

%Automorphism groups carry
The permutations on a countable set $X$ carry 
 a natural topology, the 
\emph{topology of pointwise convergence}: 
a set $F$ of permutations of $X$
is \emph{closed} iff it contains every permutation $g$ 
of $X$ such that
for every finite $S \subseteq X$ there exists
an element $f \in F$ satisfying $f(x)=g(x)$ for all $x \in S$.
Note that automorphism groups of structures are closed. 
%We write $\bar F$ for the smallest closed set of permutations that contains $F$, and
%The \emph{closure} $\cl{F}$ is the smallest
%closed permutation group that contains $F$.
% change this
%We have 
%note that $\overline{\Aut(\Gamma)} = \Aut(\Gamma)$ for every relational structure $\Gamma$.
Conversely, every closed group is the automorphism
group of a relational structure.
%for every closed group $G$ that contains
%the automorphism group of an $\omega$-categorical structure
%$\Delta$ there exists a reduct $\Gamma$ of $\Delta$ 
%such that $G = \Aut(\Gamma)$. 
When $F$ is a set of permutations, we write $\cl{F}$
for the smallest closed permutation group that contains $F$. 
The closed groups on a countable set
%that contain $\Aut(\Delta)$ 
form a lattice with
respect to inclusion: for two such groups $G_1$ and $G_2$,
the \emph{meet} $G_1 \wedge G_2$ of $G_1$ and $G_2$ is their intersection $G_1 \cap G_2$,
whereas their \emph{join} $G_1 \vee G_2$ is $\cl{G_1 \cup G_2}$.
From the above it follows that the
poset $\preceq$ of the reducts of an $\omega$-categorical structure $\Delta$ up to
equivalence is antiisomorphic to the lattice of closed groups
containing $\Aut(\Delta)$. 
In particular, it is itself a lattice. This lattice has a largest element, which is the equivalence class of $\Delta$; all other elements will be called the \emph{proper} reducts of $\Delta$.
It also has a smallest element, the equivalence class of $(D;=)$, which
will be called the \emph{trivial} reduct.

%It is now easy to see the claim from the introduction
%that the poset obtained from factoring the order 
%$\preceq$ of the reducts of an $\omega$-categorical structure $\Delta$ by equivalence 
%is in fact a lattice, dually isomorphic to the 
%lattice of closed groups containing $\Aut(\Delta)$. 
%The \emph{join} $\Gamma_1 \vee \Gamma_2$ of two reducts $\Gamma_1$ and $\Gamma_2$ is simply the reduct that contains the relations of both 
%$\Gamma_1$ and $\Gamma_2$. The \emph{meet} $\Gamma_1 \wedge \Gamma_2$ is the reduct that contains all relations with a first-order definition in $\Delta$ that are preserved by 
%both $\Aut(\Gamma_1)$ and $\Aut(\Gamma_2)$.

We write $\mathfrak L$ for the closed groups that
contain $\Aut(D;<,E)$.
Since we are mostly taking the group-theoretic perspective,
$\mathfrak L$ is the lattice that we are going to work with
in the following.  
The key to efficiently describe $\mathfrak L$ 
is to first identify the \emph{join irreducible elements} of $\mathfrak L$, that is, the closed groups $G_0$ that properly contain $\Aut(\D;<,E)$ and that have
the property that whenever $G_0 = G_1 \vee G_2$ for $G_1,G_2 \in \mathfrak L$, then $G_0 = G_1$ or $G_0 = G_2$. We will prove that $\mathfrak L$ has eleven join irreducibles, and that every element of $\mathfrak L$ is a join of join irreducibles. From this it already follows that $\mathfrak L$ is finite; by a finer study of the inclusions
between the groups we prove that there are 42 proper non-trivial reducts of the random ordered graph. 

The remainder of this section is organized as follows: we first 
recall the classical descriptions of the reducts of $({\mathbb Q};<)$ and the random graph in Subsection~\ref{sect:straightforward}, since those descriptions
appear as subclassifications of our result. 
It turns out that the random tournament is also a reduct 
of the random ordered graph. 
We recall the descriptions of the reducts of the random tournament,
and then describe the groups that arise as intersections 
of automorphism groups of previously described reducts, in Subsection~\ref{sect:coupled}. 
All but three reducts of the random ordered graph have an automorphism group that arises in this way; those three `sporadic' reducts are described in Subsection~\ref{sect:sporadic}. Finally, we show a picture of $\mathfrak L$ in Subsection~\ref{sect:picture}.

%The domain of the ordered random graph will be denoted by $V$, and its relations $E$ and $<$. We denote reducts by tuples of relations that define them. The relations that appear are those of the random graph ($E$, $\R3$, $\R4$, $\R5$,=) and those of the rationals ($<$, $\betw$, $\cyc$, $\sep$,=).
%In the sequel, we will always consider the order of the lattice of closed permutation groups, i.e., we say that a reduct $\Gamma$ is larger than a reduct $\Gamma'$ iff $\Aut(\Gamma)$ contains $\Aut(\Gamma')$. We will denote permutation groups by tuples of functions by which they are generated over $\Aut(E,<)$; when the tuple has just one element, we simply specify the function. We moreover identify reducts with their permutation groups, i.e., we may say things like ``$(f,g)$ is minimal over $(U,V)$'', for functions $f,g$ and relations $U,V$.

\subsection{Straightforward reducts}
\label{sect:straightforward}
For $i \in \{1,2\}$, let $\Gamma_i$ be the \Fresse-limit
of two strong amalgamation classes, with signature $\tau_i$ and with $n_i$ reducts up to equivalence. 
Then the superposition 
$\Delta$ of $\Gamma_1$ and $\Gamma_2$ has
at least $n_1n_2$ reducts up to equivalence, obtained in the following
straightforward way. Recall that the $\tau_i$-reduct of 
$\Delta$ is isomorphic to $\Gamma_i$, so we identify
$\Gamma_i$ with this reduct. Pick a reduct $\Gamma_1'$ of $\Gamma_1$ and a reduct $\Gamma_2'$ of $\Gamma_2$, and consider the structure $\Gamma$ with the same domain as $\Delta$ obtained by adding the relations 
from both $\Gamma_1'$ and $\Gamma_2'$. 
Observe that the automorphism group of $\Gamma$
equals $\Aut(\Gamma_1') \cap \Aut(\Gamma_2')$. 
It is clear that when in the above construction we replace $\Gamma_1'$ or $\Gamma_2'$ by an inequivalent reduct of $\Gamma_1$ or $\Gamma_2$, respectively, then we obtain an inequivalent reduct of $\Delta$. 

Applying this to the random ordered graph, we
obtain 25 inequivalent reducts.
Among them we find
in particular the reducts of $({\mathbb Q};<)$,
and the reducts of the random graph, which we
discuss next.

%In this subsection we describe those reducts of $(\D;<,E)$ that have already been discussed in the literature;
%these are the reducts of $({\mathbb Q};<)$, of the random graph, and of the random tournament. 

\subsubsection{The reducts of $({\mathbb Q};<)$}
Consider the structure obtained from $({\mathbb Q};<)$
by picking an irrational number $\pi$, and flipping the
order between the intervals $(-\infty,\pi)$ and $(\pi,\infty)$.
The resulting structure is isomorphic to $({\mathbb Q};<)$;
write $\turn$ for such an isomorphism (a permutation of ${\mathbb Q}$).
Write $\leftrightarrow$ for the operation $x \mapsto -x$.  
The automorphism groups of the proper non-trivial reducts of $({\mathbb Q};<)$ can now be described
as follows: %besides the trivial $\Aut({\mathbb Q};<)$ and $\Aut({\mathbb Q};=)$, 
there are $\cl{\Aut({\mathbb Q};<) \cup \{\leftrightarrow\}}$, $\cl{\Aut({\mathbb Q};<) \cup \{\turn\}}$, 
and $\cl{\Aut({\mathbb Q};<) \cup \{\leftrightarrow,\turn\}}$. This follows from the result due to Cameron mentioned in the beginning of Section~1.
%~\cite{Cameron5}.
%; also see~\cite{BP-reductsRamsey}.
%All of those groups except for $\cl{\Aut({\mathbb Q};<) \cup \{\leftrightarrow,\turn\}}$ are join irreducible in the lattice of closed supergroups of $\Aut({\mathbb Q};<)$ (but only one is join irreducible in $\mathfrak L$; more on that later). 

It will be useful later to also know relational descriptions
of the reducts of $({\mathbb Q};<)$. 
Let $\Betw$, $\Cycl$, and $\Sep$ be the relations
with the following definitions over $({\mathbb Q};<)$.
\begin{itemize}
\item $\Betw(x,y,z) \Leftrightarrow  \big((x< y \wedge y<z) \vee (z<y \wedge y<x)\big)$, 
\item $\Cycl(x,y,z) \Leftrightarrow \big((x< y \wedge y<z) \vee (y< z \wedge z<x) \vee (z<x \wedge x<y)\big)$, and
\item $\Sep(x,y,u,v) \Leftrightarrow \big((\Cycl(x,y,u) \wedge \Cycl(x,v,y))  \vee (\Cycl(x,u,y) \wedge \Cycl(x,y,v))\big )$. 
\end{itemize}

The following is well-known (see e.g.~\cite{JunkerZiegler}): 
\begin{itemize}
\item $\cl{\Aut({\mathbb Q};<) \cup \{\leftrightarrow\}} = \Aut({\mathbb Q};\Betw)$, 
\item $\cl{\Aut({\mathbb Q};<) \cup \{\turn\}} = \Aut({\mathbb Q};\Cycl)$, and 
\item $\cl{\Aut({\mathbb Q};<) \cup \{\leftrightarrow,\turn\}} = \Aut({\mathbb Q};\Sep)$.
\end{itemize}

\subsubsection{The reducts of the random graph}
\label{sect:random-graph}
The \emph{random graph} is the up to isomorphism unique countable homogeneous graph $(V;E)$ that is \emph{universal} 
in the sense that it contains all countable graphs as an induced subgraph. Equivalently, the random graph
is the up to isomorphism unique countable graph
with the following graph extension property (see~\cite{Hodges,Oligo}): 
for all disjoint finite subsets $U_1,U_2$ of $V$, there exists a vertex $v \in V \setminus (U_1 \cup U_2)$ that is connected to all vertices in $U_1$ and to none in $U_2$. 

To describe the reducts of the random graph $(V;E)$, it will again be more instructive to describe their automorphism groups. 
The complement graph of $(V;E)$ 
%graph $(V;{{V}\choose{2}} \setminus E)$ 
satisfies the graph extension property as well, and hence there exists an isomorphism between these two structures, a permutation of $V$ we denote by $-$. 
Let $v \in V$. Consider the graph $(V;E')$ obtained from $(V;E)$ by `switching' between 
edges and non-edges on all pairs $(u,v)$ with $u \in V \setminus \{v\}$. That is, if $(u,v) \in E$ then $(u,v) \notin E'$,
and if $(u,v) \notin E$ then $(u,v) \in E'$. For all $u,u' \in V \setminus \{v\}$ we have $(u,u') \in E$ if and only if 
$(u,u') \in E'$. Then $(V;E')$ also satisfies the graph extension property, and hence there exists an isomorphism
between $(V;E)$ and $(V;E')$, which we denote by $\sw$. 
Now the automorphism groups of the  non-trivial  reducts
of $(V;E)$ are $\cl{\Aut(V;E) \cup \{-\}}$, $\cl{\Aut(V;E) \cup \{\sw\}}$, 
and $\cl{\Aut(V;E) \cup \{-,\sw\}}$. This result is due to Thomas~\cite{RandomReducts}. 
%; also see~\cite{Thomas96,BP-reductsRamsey}.

Again, we present relational descriptions of those groups. 
For $k \geq 2$, let $R^{(k)}$ be the $k$-ary relation that
contains all $k$-tuples of pairwise distinct elements $x_1,\dots,x_k$ in $V$
such that the number of (undirected) edges 
between those elements is odd.
The following is well-known (see~\cite{Thomas96}).
% (see e.g.~\cite{RandomMinOps}): 
\begin{itemize}
\item $\cl{\Aut(V;E) \cup \{-\}} = \Aut(V;R^{(4)})$, 
\item $\cl{\Aut(V;E) \cup \{\sw\}} = \Aut(V;R^{(3)})$, and 
\item $\cl{\Aut(V;E) \cup \{-,\sw\}} = \Aut(V;R^{(5)})$.
\end{itemize}

%\subsection{The reducts of the random tournament}
\subsection{Coupled reducts}
\label{sect:coupled}
Consider the reduct of the random ordered graph that contains the binary relation $T$ defined over $(D;<,E)$ by the formula $$x \neq y \wedge (x<y \Leftrightarrow E(x,y)) \; .$$
The reduct $(D;T)$ is isomorphic to the \emph{random tournament}, defined below, and
it is straightforward to show that the automorphism group of $(D;T)$
is \emph{not} an intersection of automorphism groups of reducts of $(D;<)$ or reducts of $(D;E)$.
The reducts of the random tournament have
been classified by Bennett~\cite{Bennett-thesis},
and in Subsection~\ref{sect:tournament} we 
give a brief description of those reducts. 
%recall threducts of the random tournament.
In Subsection~\ref{sect:coupled-ops} we give an operational description of those groups that arise as intersections of automorphism groups of
reducts of $(D;<)$, $(D;E)$, or $(D;T)$. 

\subsubsection{The reducts of the random tournament}
\label{sect:tournament}
A \emph{tournament} is a directed graph such
that for all distinct vertices $u,v$ exactly one of $(u,v)$ and $(v,u)$ is an edge of the graph. 
In the following, we will use basic concepts from model theory, as they are used e.g. in~\cite{Hodges}. 
The class of all finite tournaments forms an amalgamation class, and the corresponding \Fresse-limit 
will be called the \emph{random tournament}. 
The random tournament is the up to isomorphism unique countable tournament with the following tournament
extension property: for all finite subsets $U_1,U_2$ of vertices, there exists a vertex $v \nin U_1 \cup U_2$ 
such that $(u,v)$ is an edge for all $u \in U_1$,
and $(v,u)$ is an edge for all $u \in U_2$. 
%Now consider the relation $T$ defined by $x \neq y \wedge (x<y \Leftrightarrow E(x,y))$ over $(\D;<,E)$.
It is straightforward to verify that the relation $T(x,y)$ defined by $x \neq y \wedge (x<y \Leftrightarrow E(x,y))$ over $(\D;<,E)$ 
satisfies the tournament extension property.
Thus, the random tourament is (isomorphic to) the reduct $(\D;T)$ of $(\D;<,E)$.

The tournament obtained by \emph{flipping} the orientation
of all edges again satisfies the tournament extension property;
we denote the isomorphism by $\rightleftharpoons$. 
The tournament obtained by flipping the orientation 
of all edges that are adjacent to some fixed vertex 
satisfies the tournament extension property; denote
by $\tsw$ the respective permutation of $D$. 
Then the automorphism groups of the non-trivial reducts
of $(D;T)$ are: 
$\cl{\Aut(D;T) \cup \{\rightleftharpoons\}}$, $\cl{\Aut(D;T) \cup \{\tsw\}}$, 
and $\cl{\Aut(D;T) \cup \{\leftrightharpoons,\tsw\}}$. This result is due to Bennet~\cite{Bennett-thesis}.

Let $\Betw_T$ be the ternary relation with the following first-order definition over $T$. 
$$\Betw_T(x,y,z) \Leftrightarrow \big (T(x,y) \wedge T(y,z) \wedge T(z,x)\big) \vee \big(T(z,y) \wedge T(y,x) \wedge T(x,z) \big )$$ 
Define the ternary relation $\Cycl_T$ by
\begin{align*} \Cycl_T(x,y,z) \Leftrightarrow  \big ( & (T(x,y) \wedge T(y,z) \wedge T(z,x)) \\
\vee & (T(x,z) \wedge T(z,y) \wedge T(x,y)) \\
\vee & (T(y,x) \wedge T(x,z) \wedge T(y,z)) \\
\vee & (T(z,y) \wedge T(y,x) \wedge T(z,x)) \big )
\end{align*}

Finally, let $\Sep_T$ be the relation of arity four that
contains all tuples $(x,y,u,v) \in D^4$ such that
$\big|\big\{ T \cap (\{x,y\} \times \{u,v\})\big\}\big|$ is even.  
It is clear that also $\Sep_T$ is first-order definable over $(D;T)$. 

%\begin{align*}
%\Sep_T(x,y,u,v) \Leftrightarrow \big( & 
%(T(x,u) \wedge T(x,v) \wedge T(y,u) \wedge T(y,v)) \\
%\vee & (T(x,u) \wedge T(x,v) \wedge T(y,u) \wedge T(y,v)) \\
%\vee & (T(x,u) \wedge T(x,v) \wedge T(y,u) \wedge T(y,v)) \\
%\vee & (T(x,u) \wedge T(x,v) \wedge T(y,u) \wedge T(y,v)) \\
%\vee & (T(x,u) \wedge T(x,v) \wedge T(y,u) \wedge T(y,v)) \\
%\vee &  ((\neg T(x,u) \vee \neg T(x,v)) \wedge (\neg T(y,u) \vee \neg T(y,v))) \big )
%\end{align*}

The following is a consequence of Bennett's classification. 
\begin{itemize}
\item $\cl{\Aut({\mathbb Q};T) \cup \{\leftrightharpoons\}} = \Aut({\mathbb Q};\Betw_T)$, 
\item $\cl{\Aut({\mathbb Q};T) \cup \{\tsw\}} = \Aut({\mathbb Q};\Cycl_T)$, and 
\item $\cl{\Aut({\mathbb Q};T) \cup \{\leftrightharpoons,\tsw\}} = \Aut({\mathbb Q};\Sep_T)$.
\end{itemize}

\subsubsection{Reducts from coupling operations}
\label{sect:coupled-ops}
In this section we will define some more reducts by constructing permutations of $\D$
that, intuitively, combine the behavior of 
$\turn$, $\leftrightarrow$, and the identity with the behavior of $\sw$, $-$, and the identity in various ways. 
In the following we exclusively work with groups that
contain $\Aut(\D;<,E)$.
For this reason, and in order to be concise, 
we write $\cl{F}$ instead of $\cl{F \cup \Aut(\D;<,E)}$.

Consider the structure obtained from $(\D;<,E)$ by flipping simultaneously edges with non-edges and the order. The resulting structure is isomorphic with $(\D;<,E)$; let $\mix - \rev$ be an isomorphism witnessing this fact.
Similarly, the structure obtained by flipping edges with non-edges whilst keeping the order is isomorphic to $(\D;<,E)$; let $\mix - \id$ be an isomorphism witnessing this fact.
We furthermore define in an analogous fashion permutations $\mix \id \rev$, $\mix \id \turn$, and  $\mix \sw \id$.
Slightly less obvious is ``coupling'' $\sw$ with $\turn$: Pick an irrational $\pi$, and consider the structure obtained by flipping edges and non-edges which cross $\pi$, as well as flipping the order between the intervals $(-\infty,\pi)$ and $(\pi,\infty)$. The resulting structure is isomorphic to $(\D;<,E)$, and we pick an isomorphism $\mix \sw \turn$ witnessing this.

With the relational descriptions of the reducts of $({\mathbb Q};<)$, of the random graph, and of the random tournament it is straightforward to show that these operations do not generate each other. In Figure~\ref{fig:JI1} we introduce names for the groups generated by each of these operations; it will turn out that these groups are join irreducible in $\mathfrak L$. %We will use those letters to also denote the respective reducts; whether we mean reduct or group will always be clear from the context. 
In fact, we will later see that the groups $a,b,c,d,e$ and $f$ are precisely the \emph{atoms} of $\mathfrak L$, that is, every closed group properly containing $\Aut(\D;<,E)$ contains
at least one of $a,b,c,d,e,f$. The remaining join irreducible elements of $\mathfrak L$ will be defined in Section~\ref{sect:sporadic} and listed in Figure~\ref{fig:JI2}.

\begin{figure}
\begin{center}
\begin{tabular}{|l|cccccc}
\hline
Name & a & b & c & d & e & f \\
\hline
& \\[-1ex]
Description & $\cl{\mix\id\rev}$ & $\cl{\mix\id\turn}$ & $\cl{\mix -\id}$ & $\cl{\mix \sw\id}$ & $\cl{\mix -\rev}$ & $\cl{\mix \sw\turn}$\\[1.1ex]
\hline
\end{tabular}
\end{center}
\caption{Join irreducible elements of $\mathfrak L$ generated by coupled operations introduced in Section~\ref{sect:coupled-ops}}
\label{fig:JI1}
\end{figure}

At this point, one is tempted to conjecture that all the reducts
of the random ordered graph are obtained from
reducts of $(\D;<)$, $(\D;E)$, and $(\D;T)$
by combining the relations from those reducts in all possible ways. Indeed, all reducts that we have encountered so far can be obtained in this form. 
%straightforward combinations of the reducts of $({\mathbb Q};<)$ and the random graph.
%an example of such a straightforward combination would be
%$\Aut(D,E,\Betw)$.
However, as we will see in the following, there are more reducts. 

%the groups that we describe in Section~\ref{sect:medium} are \emph{not} such straightforward combinations of reducts of $({\mathbb Q};<)$ and the random graph. 

\subsection{Sporadic Reducts}
\label{sect:sporadic}
Let $\apl$ be a permutation of $\D$ that
preserves $<$ and switches the graph relation 
below some irrational $\pi$, and leaves it unaltered otherwise. Analogously, let $\apu$ be a permutation which preserves
$<$ and switches the graph relation 
above $\pi$. To show that the corresponding reducts are pairwise distinct, and also distinct from all the previous reducts, we make the following definitions.
\begin{align*} 
R^l_3 := & \{(a_1,a_2,a_3) \; | \; a_1 < a_2 < a_3 \text{ and } E(a_1,a_3) \Leftrightarrow E(a_2,a_3)\} \\
R^u_3 := & \{(a_1,a_2,a_3) \; | \; a_1 < a_2 < a_3 \text{ and } E(a_1,a_3) \Leftrightarrow E(a_1,a_2)\} 
\end{align*}
It will be shown in Section~\ref{sect:medium-sized} that $\cl{\apl} = \Aut(D;R^l_3)$,
and $\cl{\apu} = \Aut(D;R^u_3)$. 

In Figure~\ref{fig:JI2}, we assign names to the remaining join irreducibles of $\mathfrak L$, including the groups we just defined.

\begin{figure}
\begin{center}
\begin{tabular}{ccccc|l|}
\hline
 g & h & i & j & k & Name \\
\hline
& & & & & \\[-1ex]
$\cl{\apl}$ & $\cl{\apu}$ & $\Aut(\D;E)$ & $\Aut(\D;<)$ & $\Aut(\D;T)$ & Description \\[1.1ex]
\hline
\end{tabular}
\end{center}
\caption{The remaining join irreducible elements of $\mathfrak L$, including some sporadic ones introduced in Section~\ref{sect:sporadic}}
\label{fig:JI2}
\end{figure}

% can be found in Section~\ref{sect:verification}; those descriptions also show that 
%the corresponding reducts are pairwise distinct, and also
%distinct from all the previously described reducts. 
%The relation $$R^-_3 = \{(a_1,a_2,a_3) \; | \; a_1<a_2<a_3 \text{ and } E(a_1,a_3) \Leftrightarrow E(a_2,a_3)\}$$ is preserved by $\apl$, but not by $\apu$.
%Dually, the relation $$R^+_3 = \{(a_1,a_2,a_3) \; | \; a_1<a_2<a_3 \text{ and } E(a_1,a_3) \Leftrightarrow E(a_1,a_2)\}$$ is preserved by $\apu$, but not by $\apl$.
%We will later show the following:
%\begin{itemize}
%\item $\cl{\apl} = \Aut(\D;R^-_3)$, and
%\item $\cl{\apu} = \Aut(\D;R^+_3)$.
%\end{itemize}

\subsection{The big picture}
\label{sect:picture}
The main result of the paper is the following. 

\begin{thm}
The random ordered graph $(D;<,E)$ has exactly 42 proper non-trivial reducts up to equivalence. 
\end{thm}

As discussed before, this is equivalent to
the statement that there are 44 closed permutation 
groups that contain the automorphism group of the random ordered graph.

%We will later see that the groups $\Aut(\D;<)$, $\Aut(\D;T)$, and $\Aut(\D;E)$ that we have
%seen above are join irreducibles of $\mathfrak L$. 
%We will call those the \emph{large} join irreducibles. 
%We now describe the remaining join irreducibles, which
%fall into two groups, the \emph{small}, and the \emph{medium-sized} join irreducibles. 
%\subsubsection{The small join irreducibles}

\begin{figure}[h]
\begin{center}
\includegraphics[scale=0.6]{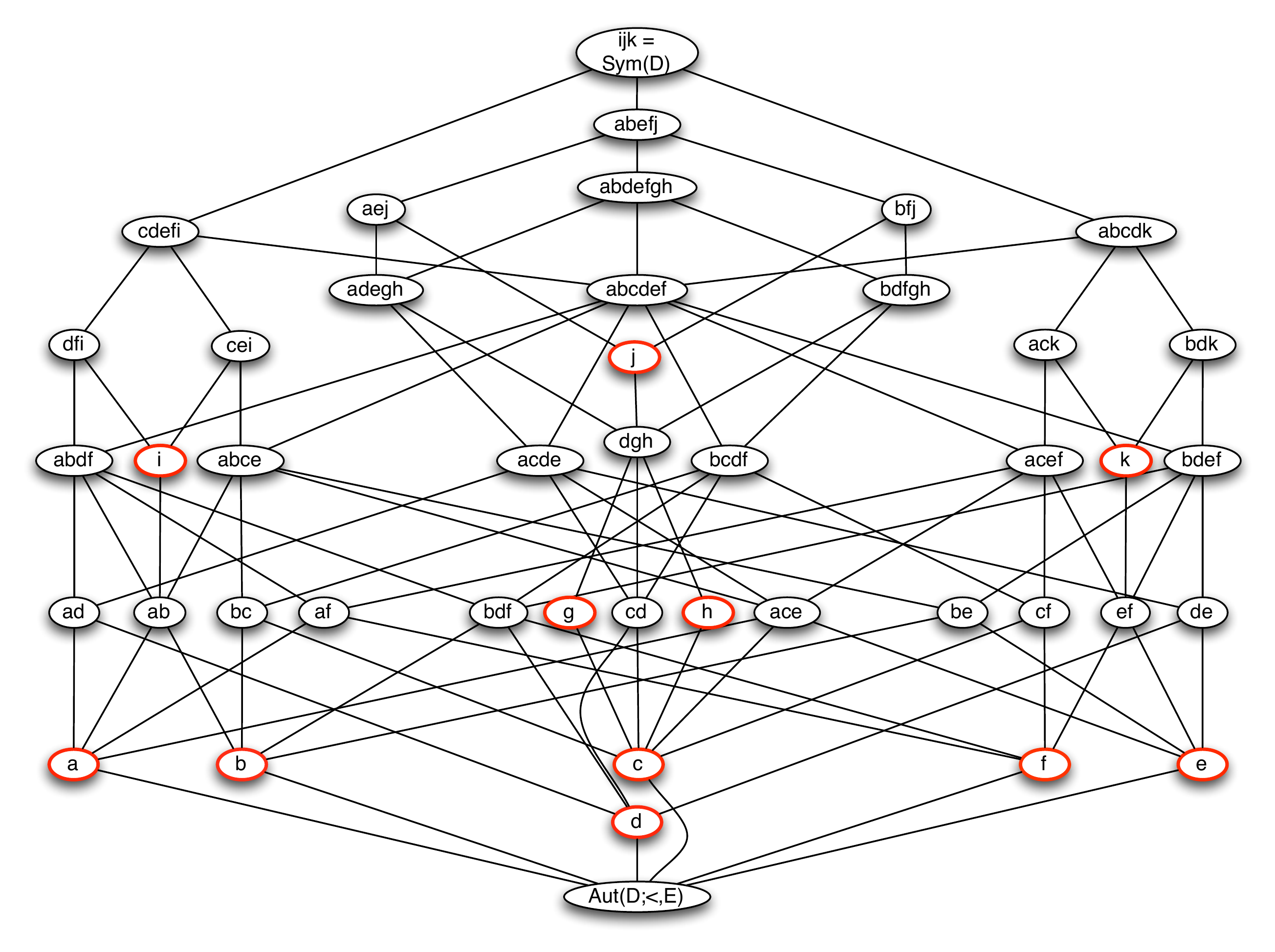}
\end{center}
\caption{The lattice $\mathfrak L$ of closed supergroups of $\Aut(\D;<,E)$.}
\label{fig:reducts}
\end{figure}

The Hasse diagram of the lattice $\mathfrak L$ of those groups is as depicted in Figure~\ref{fig:reducts}. In this figure, the bottom vertex
denotes $\Aut(\D;<,E)$. All other vertices are labeled with the set of maximal join irreducibles below them. For instance, the vertex labeled $dgh$ lies above $c$, $d$, $g$, $h$, with $d$, $g$, $h$ maximal; and it can be variously defined as $d\vee g$, $d\vee h$, or $g\vee h$.

%labeled with a sequence of join irreducibles whose join is the respective group $G$. We list all those join irreducibles contained in $G$ that are \emph{maximal in $G$},
%i.e., that are not contained in another join irreducible that is contained in $G$.
%For instance, the vertex labeled $ab$
%$denotes the group $a \vee b$ in the lattice, that is,
%$\cl{{\mix \id \leftrightarrow},{\mix \id \turn}}$. 
%However, for the group $\Aut(\D;T)$, we only
%write $k$, since the other join irreducibles $f$ and $e$ in $\Aut(\D;T)$ are not maximal in $\Aut(\D;T)$. 

The verification that the lattice indeed has the shape that is presented in~Figure~\ref{fig:reducts} will
be completed in Section~6. 
%We make crucial use of the fact
%that the lattice $\mathfrak L$ of closed supergroups
%of $\Aut(D;<,E)$ is isomorphic to the lattice of ideals of 
%the join irreducibles of $\mathfrak L$ with the order of inclusion.

Note that the groups from Section~\ref{sect:sporadic}, labeled `g', `h', and `dgh' in the picture, show that $\mathfrak L$ is not symmetric with respect to the roles of $(\D;<)$ and $(\D;E)$, that is, there is no
automorphism of $\mathfrak L$ sending $\Aut(\D;<)$
to $\Aut(\D;E)$.

\section{Overview of the Proof}\label{sect:overview}

The main part of our proof identifies the join irreducible elements of $\L$, and shows that every closed group containing $\Aut(D;<,E)$ is the join of such join irreducibles.

\begin{defn}\label{defn:JI}
We set $\Delta:=(D;<,E)$. 
Let $\JI$ be the set of those groups which were given a name in $\{a,\ldots,k\}$ in the previous section (Figures~\ref{fig:JI1} and~\ref{fig:JI2}):
\begin{itemize}
\item the ``small'' groups $\cl{\mix\id\turn},\cl{\mix\sw\turn},\cl{\mix\id\rev},\cl{\mix -\rev},\cl{\mix\sw\id},\cl{\mix -\id}$;
\item the ``medium-sized'' groups $\cl{\apl}$ and $\cl{\apu}$;
\item the ``large'' groups $\Aut(D;<), \Aut(D;E), \Aut(D;T)$.
\end{itemize}
\end{defn}

It will turn out that the groups of $\JI$ are precisely the join irreducibles of $\L$. In Sections~\ref{sect:unordered} and~\ref{sect:ordered} we will obtain the following statement.

\begin{prop}\label{prop:ji}
Let $G,H\supseteq \adel$ be closed groups such that $H\setminus G\neq \emptyset$. Then there exists an element of $\JI$ which is contained in $H$ but not in $G$.
\end{prop}

From this we immediately get

\begin{cor}\label{cor:join}
Every closed group containing $\adel$ is the join of elements of $\JI$.
\end{cor}
\begin{proof}
Let $H$ be a closed group containing $\adel$, and let $G$ be the join of all elements of $\JI$ which are contained in $H$. By Proposition~\ref{prop:ji} we must have $H\setminus G=\emptyset$, so $H=G$.
\end{proof}

It follows that $|\L|\leq 2^{|\JI|}$, and in particular that $\L$ is finite. To obtain a full proof of our classification, it remains to show for all $G\in \JI$ and all subsets $\S$ of $\JI$ that $G\subseteq \bigvee \S$ if and only if it has been drawn that way in the picture. We postpone this task to Section~\ref{sect:verification}, and first concentrate on the mathematically more interesting and challenging proof of Proposition~\ref{prop:ji}. This proof will be obtained by a Ramsey-theoretic analysis of functions on $\Delta$, and follows~\cite{BPT-decidability-of-definability, BodPin-Schaefer,RandomMinOps,BP-reductsRamsey}. Before starting out, we provide the necessary definitions for this analysis.

\begin{definition}\label{defn:type}
    Let $\Lambda$ be a structure. The \emph{type} $\tp(a)$ of an $n$-tuple $a$ of elements in $\Lambda$ is the set of first-order formulas with
     free variables $x_1,\dots,x_n$ that hold for $a$ in $\Lambda$.
\end{definition}

\begin{definition}\label{defn:behaviour}
    Let $\Lambda, \Omega$ be structures. A \emph{type condition} between $\Lambda$ and $\Omega$ is a pair $(t,s)$, where $t$ is a type of an $n$-tuple in $\Lambda$, and $s$ is a type of an $n$-tuple in $\Omega$, for some $n\geq 1$. 
    
      A function $f\colon\Lambda\To \Omega$ \emph{satisfies} a type condition $(t,s)$ between $\Lambda$ and $\Omega$ iff for all $n$-tuples $a=(a_1,\ldots,a_n)$ of elements of $\Lambda$ with $\tp(a)=t$ the $n$-tuple $f(a):=(f(a_1),\ldots,f(a_n))$ has type $s$ in $\Omega$. A \emph{behavior} is a set of type conditions between structures $\Lambda$ and $\Omega$. A function from $\Lambda$ to $\Omega$ \emph{satisfies a  behavior $B$} iff it satisfies all the type conditions of $B$.
\end{definition}

\begin{definition}\label{defn:canonical}
    Let $\Lambda, \Omega$ be structures. A function $f\colon \Lambda \To \Omega$ is \emph{canonical} iff for all types $t$ of $n$-tuples in $\Lambda$ there exists a type $s$ of $n$-tuples in $\Omega$ such that $f$ satisfies the type condition $(t,s)$. In other words, $n$-tuples of equal type in $\Lambda$ are sent to $n$-tuples of equal type in $\Omega$ under $f$, for all $n\geq 1$.
\end{definition}

Note that for a canonical function $f\colon \Lambda \To \Omega$, the set of all type conditions satisfied by $f$ is a function from the types over $\Lambda$ to the types over $\Omega$.

\begin{definition}
For a set $F$ of functions from $D$ to $D$, we denote by $\overline F$ the closure of $F$ in $D^D$, i.e., the set of all functions in $D^D$ which agree with some function in $F$ on every finite subset of $D$.
\end{definition}

\begin{definition}
We say that a set $S\subseteq D^D$ \emph{generates a function $f\in D^D$ over $\Delta$} iff $f$ is contained in the smallest closed monoid containing $S\cup\adel$; equivalently, $f$ is contained in the closure of the set of all functions that can be composed from $S$ and $\adel$. Since we will only be interested in sets of functions containing $\Aut(\Delta)$, we will say ``generates'' rather than ``generates over $\Delta$''.
\end{definition}

 The combinatorial kernel of our method is the following proposition which follows from~\cite{BPT-decidability-of-definability, BP-reductsRamsey} and the fact that the set of finite linearly ordered graphs is a \emph{Ramsey class} (due to~\cite{AbramsonHarrington,NesetrilRoedlOrderedStructures};~\cite{NesetrilRoedlPartite}). It states that if a function on $\Delta$ does something of interest on a finite set, then it produces a canonical function which still does the same thing. For $c_1,\ldots,c_n\in D$ we denote by $(\Delta,c_1,\ldots,c_n)$ the structure obtained from $\Delta$ by adding the constants $c_1,\ldots,c_n$ to the language.

\begin{prop}\label{prop:ramseystrong}
Let $f\colon D \To D$ be a function, and let $c_1,\ldots,c_n\in D$. Then 
$$
\overline{\{\alpha\circ f\circ \beta\;|\; \alpha\in \adel,\; \beta\in \adelcn\}}
$$ 
contains a function $g$ such that
\begin{itemize}
\item $g$ is canonical as a function from $(\Delta,c_1,\ldots,c_n)$ to $\Delta$;
\item $g$ agrees with $f$ on $\{c_1,\ldots, c_n\}$.
\end{itemize}
In particular, $f$ generates a function $g$ with these properties.
\end{prop}

It follows that we can distinguish the elements of $\L$ by canonical functions, as follows (see~\cite{BPT-decidability-of-definability, BP-reductsRamsey}). 

\begin{prop}\label{prop:ramsey}
Let $G,H\supseteq\adel$ be closed groups such that $H\setminus G\neq\emptyset$. Then there exist $c_1,\ldots,c_n\in D$ and a function $f\colon D\To D$ such that
\begin{itemize}
\item $f$ is canonical as a function from $(\Delta,c_1,\ldots,c_n)$ to $\Delta$;
\item $f$ is generated by $H$ but not by $G$.
\end{itemize}
\end{prop}

In Section~5, we will prove that if $G,H\supseteq \adel$ are closed groups such that $H\setminus G\neq \emptyset$, then there exists a group in $\JI$ which is contained in $H$ but not in $G$, or the function $f$ of Proposition~\ref{prop:ramsey} can be assumed to preserve the order $<$. This is achieved in Proposition~\ref{prop:unorderedCanonical}. In Section~6 we show that in the latter situation, we also find an element of $\JI$ which is contained in $H$ but not in $G$ (Proposition~\ref{prop:ordered}). This completes the proof of Proposition~\ref{prop:ji}.

\section{The Unordered Case}\label{sect:unordered}

The goal of this section is the following

\begin{prop}\label{prop:unorderedCanonical}
Let $G,H\supseteq \adel$ be closed groups such that $H\setminus G\neq \emptyset$. Then one of the following holds:
\begin{itemize}
\item There exists an element of $\JI$ which is contained in $H$ but not in $G$;
\item there exist $c_1,\ldots,c_n\in D$ and an order preserving canonical function $f\colon \delcn\To\Delta$ which is generated by $H$ but not by $G$.
\end{itemize}
\end{prop}

We will prove this proposition by analyzing possible behaviors of canonical functions guaranteed by Proposition~\ref{prop:ramsey}. To this end, let us introduce some notation and terminology.

\begin{defn}
We define a binary relation $N(x,y)$ on $D$ by the formula $\neg E(x,y)\wedge x\neq y$.
\end{defn}

\begin{defn}
Let $f:D\To D$ be injective, and let $X,Y\subseteq D$. We say that $f$
\begin{itemize}
\item \emph{preserves a relation $R(x_1,\ldots,x_n)$ on $X$} iff $R(x_1,\ldots,x_n)$ implies $R(f(x_1),\ldots,f(x_n))$ for all $x_1,\ldots,x_n\in X$;
\item \emph{reverses the order on $X$} iff $x<y$ implies $f(x)>f(y)$ for all $x,y\in X$;
\item \emph{keeps the graph relation on $X$} iff $E(x,y)\leftrightarrow E(f(x),f(y))$ for all $x,y\in X$;
\item \emph{flips the graph relation on $X$} iff $E(x,y)\leftrightarrow N(f(x),f(y))$ for all $x,y\in X$;
\item \emph{eradicates edges on $X$} iff $N(f(x),f(y))$ for all $x,y\in X$;
\item \emph{eradicates non-edges on $X$} iff $E(f(x),f(y))$ for all $x,y\in X$;
\item\emph{preserves a relation $R(x,y)$ between $X,Y$} iff $R(x,y)$ implies $R(f(x),f(y))$ for all $x\in X$ and $y\in Y$;
\item \emph{reverses the order between $X,Y$} iff $x<y$ implies $f(x)>f(y)$ for all $x\in X$ and $y\in Y$;
\item \emph{keeps the graph relation between $X,Y$} iff $E(x,y)\leftrightarrow E(f(x),f(y))$ for all $x\in X$ and $y\in Y$;
\item \emph{flips the graph relation between $X,Y$} iff $E(x,y)\leftrightarrow N(f(x),f(y))$ for all $x\in X$ and $y\in Y$;
\item \emph{eradicates edges between $X,Y$} iff $N(f(x),f(y))$ for all $x\in X$ and $y\in Y$;
\item \emph{eradicates non-edges between $X,Y$} iff $E(f(x),f(y))$ for all $x\in X$ and $y\in Y$.
\end{itemize}
\end{defn}

\begin{defn}
Let $c_1,\ldots,c_n, d\in D$ be distinct. The \emph{level} of $d$ in $\delcn$ is the number of $c_i$ below $d$. For $0\leq i \leq n$ we call the set of all $d\in D\setminus\{c_1,\ldots,c_n\}$ which are on level~$i$ the \emph{$i$-th level of  $\delcn$} .
\end{defn}

\begin{defn}
Let $\Omega$ be a structure. Then any set of the form $\{\alpha(x): \alpha\in\Aut(\Omega)\}$, where $x$ is an element of the domain of $\Omega$, is called an \emph{orbit} of $\Omega$.
\end{defn}

Note that if $c_1,\ldots,c_n\in D$ and $X$ is an orbit of $\delcn$, then all elements of $X$ are on the same level, which we will call the level of $X$. Observe also that $X$ induces a copy of $\Delta$.

\begin{defn}
For subsets $X,Y$ of $D$ we write $X<Y$ iff $x<y$ for all $x\in X$ and all $y\in Y$. We write $X>Y$ iff $Y<X$.
\end{defn}

\subsection{Possible behaviors}

\begin{lem}\label{lem:orderonorbit}
Let $c_1,\ldots,c_n\in D$, let $f\colon\delcn\To\Delta$ be canonical, and let $X$ be an infinite orbit of $\delcn$. Then $f$ preserves the order on $X$, or it reverses the order on $X$.
\end{lem}
\begin{proof}
Pick $x,y\in X$ with $x<y$ and $E(x,y)$. Assume that $f(x)<f(y)$. Pick $z\in X$ such that $y<z$, $E(y,z)$, and $N(x,z)$. Then since $f$ is canonical and $f(x)<f(y)$ we must have $f(y)<f(z)$, and hence $f(x)<f(z)$. Now if $u,v\in X$ with $u<v$ are arbitrary, then $f(u)<f(v)$ is implied by $f(x)<f(y)$ in case $E(u,v)$, and by $f(x)<f(z)$ in case $N(u,v)$ holds. Hence, $f$ preserves the order on $X$. Dually, we can derive that $f$ is order reversing on $X$ when $f(x)>f(y)$.
\end{proof}

\begin{lem}\label{lem:orderonsamelevel}
Let $c_1,\ldots,c_n\in D$ and let $f\colon\delcn\To\Delta$ be canonical. Let $X,Y$ be orbits of $\delcn$ on the same level, and assume that neither $f[X]<f[Y]$ nor $f[X]>f[Y]$ holds. Then either $f$ preserves the order on $X\cup Y$, or it reverses the order on $X\cup Y$.
\end{lem}
\begin{proof}
By Lemma~\ref{lem:orderonorbit}, $f$ preserves or reverses the order on $X$ and $Y$. Assuming it preserves the order on $X$, we show that it preserves the order on $X\cup Y$; the order reversing case is dual. Suppose that $x\in X$, $y\in Y$, $x<y$ and $E(x,y)$ imply $f(x)>f(y)$. Then for arbitrary $u\in X$ and $v\in Y$ we can pick $x\in X$ with $x<u$, $x<v$, and $E(x,v)$. By our assumption, we then have $f(x)>f(v)$, and moreover $f(u)>f(x)>f(v)$ since $f$ is order preserving on $X$. Hence, $f[X]>f[Y]$, a contradiction. We therefore conclude that $x\in X$, $y\in Y$, $x<y$ and $E(x,y)$ imply $f(x)<f(y)$; by duality, $x\in X$, $y\in Y$, $x<y$ and $N(x,y)$ imply $f(x)<f(y)$ as well, and so $x\in X$, $y\in Y$ and $x<y$ imply $f(x)<f(y)$.

Suppose that $f$ reverses the order on $Y$. Then the dual argument of the argument above shows that $x\in X$, $y\in Y$ and $x>y$ imply $f(x)<f(y)$. But then putting together our information we have $f[X]<f[Y]$, a contradiction. Hence, $f$ preserves the order on $Y$. Then the argument above shows that $x\in X$, $y\in Y$ and $x>y$ imply $f(x)>f(y)$, and so $f$ preserves the order on $X\cup Y$, proving the lemma.
\end{proof}

\subsection{Compactness and eradicating structure}

We will need some general facts that are often needed in the investigation of reducts of homogeneous structures. We chose to formulate them not only for $\Delta$ to emphasize their generality. The first statement can be proven by a standard compactness argument.

\begin{lem}\label{lem:compactness}
Let $\Omega$ be an $\omega$-categorical structure with domain $D$, and let $T$ be a set of finite partial functions with the following properties:
\begin{itemize}
\item $T$ is closed under restrictions;
\item every finite subset of $D$ is the domain of some function in $T$;
\item whenever $p\in T$ and $\alpha\in\Aut(\Omega)$, then $\alpha\circ p\in T$.
\end{itemize}
Let $F\subseteq D^D$ be so that every $p\in T$ is the restriction of some $f\in F$. Then
$$
\overline{\{\alpha\circ f\;|\; \alpha\in \Aut(\Omega),\; f\in F\}}
$$
contains a function $g\colon D\To D$ such that every restriction of $g$ to a finite set belongs to $T$.
\end{lem}

%\begin{defn}
%A structure is called \emph{trivial} iff it is a reduct of the \emph{empty} structure on its domain, i.e., the structure which has no relations or functions.
%\end{defn}

\begin{cor}\label{cor:eradicateglobally}
Let $G\supseteq \adel$ be a closed group, and assume that $G$ generates for every finite $A\subseteq D$ a function $f\colon D\To D$ which eradicates edges or non-edges on $A$.  Then $G$ generates a canonical function $g\colon D\To D$ which eradicates edges or non-edges.
\end{cor}
\begin{proof}
By Lemma~\ref{lem:compactness}, $G$ generates a function $h\colon D\To D$ which eradicates edges or non-edges. By Proposition~\ref{prop:ramseystrong},
$$
\overline{\{\alpha\circ h\circ\beta\;|\; \alpha, \beta\in \Aut(\Delta)\}}
$$
contains a canonical function; clearly, this function still eradicates edges or non-edges.
\end{proof}

\subsection{Moving above $\Aut(D;<)$}

\begin{lem}\label{lem:eradicatesonfinite}
Let $G\supseteq \adel$ be a closed group, and assume it generates for every finite $A\subseteq D$ a function $f\colon D\To D$ which eradicates edges or non-edges on $A$.  Then $G\supseteq \Aut(D;<)$.
\end{lem}
\begin{proof}
By Lemma~\ref{lem:orderonorbit} and Corollary~\ref{cor:eradicateglobally}, $G$ generates an order preserving or order reversing function $e\colon D\To D$ which eradicates edges or non-edges on $D$. Replacing $e$ by $e^2$, we then have that $e$ preserves the order. 
Let $\Gamma$ be a reduct of $\Delta$ whose automorphism group equals $G$, and let $\Gamma'$ be the structure induced by $e[D]$ in $\Gamma$.  Because $\Delta$ has quantifier elimination, the restriction of any formula over $\Delta$ to $e[D]$ is equivalent to a formula over $(D;<)$. Hence, $\Gamma'$ is definable in $(e[D];<)$. Since $e^{-1}$ is an isomorphism from $\Gamma'$ onto $\Gamma$ as well as from $(e[D];<)$ to $(D;<)$, $\Gamma$ is definable in $(D;<)$, and so its automorphism group $G$ contains $\Aut(D;<)$.
\end{proof}

\begin{lem}\label{lem:eradicatesonorbit}
Let $G\supseteq \adel$ be a closed group, and assume it generates a canonical function $f\colon\delcn\To\Delta$, where $c_1,\ldots,c_n\in D$, which eradicates edges or non-edges on an infinite orbit of $\delcn$. Then $G\supseteq \Aut(D;<)$.
\end{lem}
\begin{proof}
Let $X$ be an infinite orbit on which $f$ eradicates edges or non-edges. The structure induced by $X$ is isomorphic to $\Delta$. Hence, any tuple of elements of $D$ can be sent into $X$ by an automorphism of $\Delta$, and the claim follows from Lemma~\ref{lem:eradicatesonfinite}.
\end{proof}

\begin{lem}\label{lem:eradicatebetweenorbits}
Let $G\supseteq \adel$ be a closed group, let $c_1,\ldots,c_n\in D$, and let $f\colon\delcn\To\Delta$ be a canonical function generated by $G$. If there exist infinite orbits $X,Y$ of $\delcn$ such that $f$ eradicates edges or non-edges between $X$ and $Y$, then $G\supseteq \Aut(D;<)$.
\end{lem}
\begin{proof}
If $f$ eradicates edges or non-edges on $X$ or $Y$ then we are done by Lemma~\ref{lem:eradicatesonorbit}, so we may assume that it keeps or flips the graph relation on $X$ and $Y$. 

\textit{Case~1.} $f$ keeps the graph relation on both $X$ and $Y$. Let $A\subseteq D$ be finite. If there exist $a,a'\in A$ with $E(a,a')$, then let $\alpha\in\adel$ be so that $\alpha[A]\subseteq X\cup Y$ and such that $\alpha(a)\in X$ iff $\alpha(a')\in Y$. Then $f\circ\alpha[A]$ has less edges than $A$, and so the iteration of this process allows us to send $A$ to an independent set. The statement now follows from Lemma~\ref{lem:eradicatesonfinite}.

\textit{Case~2.} $f$ flips the graph relation on both $X$ and $Y$. We may assume that $f[D]\subseteq Y$. Then replacing $f$ by $f^2$ brings us back to Case~1.

\textit{Case~3.} $f$ flips the graph relation on precisely one of the sets $X$ and $Y$, say without loss of generality on $X$. 
%We may further assume that $X<Y$ or that $X,Y$ are on the same level. 
If $f$ preserves $<$ between $X$ and $Y$ and between $Y$ and $X$, then we may assume $f[X]\subseteq X$ and $f[Y]\subseteq Y$; replacing $f$ by $f^2$ then brings us back to Case~1. Otherwise we may assume $f[Y]\subseteq X$ and $f[X]\subseteq Y$; replacing $f$ by $f^2$ then brings us back to Case~2. 
\end{proof}

\subsection{Moving above $\Aut(D;<)$ or turns}

\begin{lem}\label{lem:orderbetweenlevels}
Let $G\supseteq \adel$ be a closed group, and let $c_1,\ldots,c_n\in D$. Assume $G$ generates a canonical function $f\colon\delcn\To\Delta$ such that for two infinite orbits $X,Y$ of $\delcn$ with $X<Y$ we have $f[X]>f[Y]$, and such that $f$ is order preserving on either $X$ or $Y$. Then $G$ contains $\mix\id\turn$, $\mix\sw\turn$, or $\Aut(D;<)$.
\end{lem}
\begin{proof}
If $f$ eradicates edges or non-edges on $X$, on $Y$, or between $X$ and $Y$, then we are done by Lemmas~\ref{lem:eradicatesonorbit} and~\ref{lem:eradicatebetweenorbits}; therefore, we may assume that $f$ keeps or flips the graph relation on $X$, on $Y$, and between $X$ and $Y$. Say without loss of generality that $f$ is order preserving on $X$. If it flips the graph relation on $X$, then it generates $\mix -\id$, and replacing $f$ by $\mix -\id\circ f$ gives us a function which keeps the graph relation on $X$. Now let  $a_1,\ldots,a_k\in D$ be arbitrary, and assume without loss of generality $a_1<\cdots < a_k$. Let $\alpha\in\adel$ be so that it sends $a_1,\ldots,a_{k-1}$ into $X$ and $a_{k}$ into $Y$. Set $g:=f\circ \alpha$. Then $g(a_k)<g(a_1)<\cdots<g(a_{k-1})$. If $f$ keeps the graph relation between $X$ and $Y$, this shows that we can change the order between the $a_i$ cyclically without changing the graph relation by repeated applications of $f$ and automorphisms of $\Delta$, and so $f$ generates $\mix \id\turn$. Otherwise, it flips the graph relation between $X$ and $Y$, and hence application of $g$ flips the graph relation between $a_k$ and the other $a_i$. This shows that we can change the order between the $a_i$ cyclically by repeated applications of $f$ and automorphisms of $\Delta$ in such a way that the graph relation changes between $a_i$ and $a_j$ whenever the order changes, and so $f$ generates $\mix \sw\turn$.
\end{proof}

\subsection{Moving above $\Aut(D;<)$, $\Aut(D;E)$, or $\Aut(D;T)$}

\begin{lem}\label{lem:movesctotop}
Let $G\supseteq \Aut(\Delta)$ be a closed group. Assume that it generates a canonical function $f\colon(\Delta,c)\To\Delta$ such that $f$ is order preserving on and between the infinite orbits of $(\Delta,c)$, and such that $f(c)>f(v)$ for all $v\in D\setminus\{c\}$. Then $G$ contains $\Aut(D;<)$, $\Aut(D;E)$, or $\Aut(D;T)$.
\end{lem}
\begin{proof}
By Lemmas~\ref{lem:eradicatesonorbit} and~\ref{lem:eradicatebetweenorbits}, we may assume that $f$ keeps or flips the graph relation on every infinite orbit and between any two infinite orbits. 
Since $f$ is order preserving on and between the infinite orbits, it is easy to see that there exists a self-embedding $e$ of $\Delta$ such that $(e\circ f)[X]\subseteq X$ for every infinite orbit $X$ of $(\Delta,c)$. Replacing $f$ by $e\circ f$, we henceforth assume that $f$ itself has this property. Setting $g:=f^2$, we then have: 
\begin{itemize}
\item $g$ is canonical as a function from $(\Delta,c)$ to $\Delta$;
\item $g(c)\geq g(v)$ for all $v\in D$, and $g$ is order preserving on and between the infinite orbits of $(\Delta,c)$;
\item $(g(x),g(y))\in E$ iff $(x,y)\in E$ for all $x,y\in D\setminus\{c\}$.
\end{itemize}

Suppose that $g$ eradicates edges between $\{c\}$ and level~1. Let $a_1,\ldots,a_n \in D$, and say without loss of generality $a_1<\ldots<a_n$. Let $\alpha\in\adel$ be so that $\alpha(a_1)=c$, and set $h:=g\circ\alpha$. Then $N(h(a_1),h(a_i))$ for all $2\leq i\leq n$, and otherwise $E$ and $N$ are preserved between the $a_i$ under $h$. Moreover, $h(a_2)<\ldots<h(a_n)<h(a_1)$. Iterating this process, we can send the $a_i$ to an independent set, and so Lemma~\ref{lem:eradicatesonfinite} implies that $G$ contains $\Aut(D;<)$. The same argument works when $g$ eradicates edges between $\{c\}$ and level~0, and the dual argument works when $g$ eradicates non-edges between $\{c\}$ and an entire level.

\textit{Case 1.} $g$ preserves $E$ and $N$. Then for any $a_1,\ldots,a_n\in D$ there is $h$ generated by $g$ such that $h(a_1)<\cdots<h(a_n)$. Since $h$ preserves $E$ and $N$ we have $G\supseteq \Aut(D;E)$.

\textit{Case 2.} $g$ flips $E$ and $N$ between $c$ and $D\setminus \{c\}$. By composing $g$ with a self-embedding of $\Delta$ we can achieve $g(c)=c$ while keeping the properties listed above. Then $g^2$ preserves $E$ and $N$ and we are back in the preceding case.

\textit{Case 3.} $g$ flips $E$ and $N$ between $c$ and level~1, and preserves $E$ and $N$ between $c$ and level~0. We will show that $G$ contains $\Aut(D;T)$. 
Let $a_1,\ldots,a_n,b_1,\ldots,b_n\in D$ be so that the mapping $\xi$ which sends every $a_i$ to $b_i$ is an isomorphism with respect to the structures which $\{a_1,\ldots,a_n\}$ and $\{b_1,\ldots,b_n\}$ induce in $(D;T)$. Say without loss of generality $a_1<\cdots<a_n$, and write $b_{j_1}<\cdots<b_{j_n}$. Let $\alpha_1\in\adel$ send $a_{j_1}$ to $c$, and set $h_1:=g\circ\alpha_1$. Let $\alpha_2\in\adel$ send $h_1(a_{j_2})$ to $c$, and set $h_2:=g\circ\alpha_2\circ h_1$. Continue like this, arriving at $h:=h_n$. We then have that $h(a_{j_1})<\cdots<h(a_{j_n})$. Moreover, under $h$ edges and non-edges are flipped between elements $a_i,a_j$ if and only if the order is flipped between these elements. Hence, $h$ restricted to $\{a_1,\ldots,a_n\}$ is an isomorphism between induced substructures of $(D;T)$. It follows that the mapping which sends every $h(a_{j_i})$ to $b_{j_i}$ is a partial isomorphism on $\Delta$. Let $\beta\in\adel$ agree with this partial isomorphism. Then $\beta\circ h$ sends $a_i$ to $b_i$, and so $G\supseteq \Aut(D;T)$.

\textit{Case 4.} $g$ flips $E$ and $N$ between $c$ and level~0, and preserves $E$ and $N$ between $c$ and level~1. Assume without loss of generality that $g(c)=c$. Then considering $g^2$ instead of $g$ brings us back to Case~3.
\end{proof}

\begin{lem}\label{lem:orderonlevel}
Let $G\supseteq \adel$ be a closed group, and assume it generates a canonical function $f\colon\delcn\To\Delta$, where $c_1,\ldots,c_n\in D$. Then one of the following holds:
\begin{itemize}
\item $f$ preserves the order on all levels of $\delcn$;
\item $f$ reverses the order on all levels of $\delcn$;
\item $G$ contains either $\Aut(D;<)$, $\Aut(D;E)$, or $\Aut(D;T)$.
\end{itemize}
\end{lem}
\begin{proof}

%ignore
\ignore{
Let $X,Y$ be infinite orbits on the same level, and suppose that $g[X]<g[Y]$. Then pick any $c\in Y$. The structure $(X\cup\{c\};<,E,c)$ is isomorphic to $(D;<,E,c)$. Picking any isomorphism $i\colon(D;<,E,c)\To (X\cup\{c\};<,E,c)$, we have that $i$ is generated by $\adel$, and so $f\circ i$ is generated by $f$. But $f\circ i$ satisfies the hypothesis of  Lemma~\ref{lem:movesctotop}, and so we are done. Hence, Lemma~\ref{lem:orderonsamelevel} implies that $f$ preserves or reverses the order on this level.
}
%end ignore

Suppose that the first two cases of the lemma do not apply. Then there are infinite orbits $X$, $Y$ on different levels such that $f$ reverses the order on one of them, and keeps the order on the other one. Assume without loss of generality that $X<Y$. If $f$ does not preserve the order between $X$ and $Y$, then let $O\in\{X,Y\}$ be the orbit on which the order is reversed, and assume that the range of $f$ is contained in $O$ by composing it with a self-embedding of $\Delta$ if necessary. Then $f^2$ still reverses the order on precisely one of the orbits $X$ and $Y$, but preserves the order between $X$ and $Y$. Replacing $f$ by $f^2$ we may henceforth assume this situation. Say without loss of generality that $f$ reverses the order on $X$ and preserves the order on $Y$. 

If $f$ eradicates edges or non-edges on $X$ or $Y$ then we are done by Lemma~\ref{lem:eradicatesonorbit}, so assume that $f$ either keeps or flips the graph relation on $X$ and on $Y$. If $f$ flips the graph relation on $Y$, then assuming $f[D]\subseteq Y$ by virtue of the existence of an appropriate self-embedding of $\Delta$, and replacing $f$ by $f^2$, we obtain that $f$ keeps the graph relation on $Y$, which we will henceforth assume. Referring to Lemma~\ref{lem:eradicatebetweenorbits} we may also assume that $f$ either keeps or flips the graph relation between $X$ and $Y$, leaving us with four cases.

\textit{Case~1.} $f$ keeps the graph relation on $X$ and between $X$ and $Y$. Then by the above assumptions, $f$ keeps the graph relation on $X\cup Y$. We claim that in this case, $G$ contains $\Aut(D;E)$. To see this, let $a_1,\ldots,a_k\in D$ be arbitrary, and assume without loss of generality $a_1<\cdots < a_k$. Let $1\leq i\leq k$, and let $\alpha\in\adel$ be so that it sends $a_1,\ldots,a_{i-1}$ into $X$ and $a_{i},\ldots,a_k$ into $Y$. Then  $g:=f\circ\alpha$ does not change the graph relation between the $a_j$, and $g(a_{i-1})<\cdots<g(a_1)<g(a_{i})<\ldots < g(a_k)$. Now let $\beta\in\adel$ be so that it sends $g(a_{i-1}),\ldots,g(a_1),g(a_{i})$ into $X$ and $g(a_{i+1}),\ldots,g(a_k)$ into $Y$. Then $h:=f\circ\beta\circ g=f\circ\beta\circ f\circ \alpha$ does not change the graph relation on $\{a_1,\ldots,a_k\}$, and $h(a_i)<h(a_1)<\cdots <h(a_{i-1})<h(a_{i+1})<\cdots<h(a_k)$. By repeated application of functions of this form, we can change the order of the $a_j$ ad libitum without changing the graph relation. By the homogeneity of $\Delta$, this implies that any function from $\{a_1,\ldots,a_k\}$ to $D$ which keeps the graph relation can be extended to a function in $G$. Since $G$ is closed we conclude that it contains $\Aut(D;E)$.

\textit{Case~2.} $f$ keeps the graph relation on $X$ and flips it between $X$ and $Y$.  Let again $a_1,\ldots,a_k\in D$, and define $h$ as before. Then as before, $h$ moves $a_i$ below the other $a_j$; however, this time it flips the graph relation between $a_i$ and the other $a_j$. Let $\gamma\in\adel$ be so that it moves $h(a_i)$ into $X$, and the rest of the $h(a_j)$ into $Y$. Then $f\circ \gamma\circ h$ has the properties of $h$ in Case~1, and we again conclude that $G$ contains $\Aut(D;E)$.

\textit{Case~3.} $f$ flips the graph relation on $X$ and keeps it between $X$ and $Y$. Proceeding as in Case~1, we then see that by applications of functions in $G$ we can change the order among any $a_1,\ldots,a_k\in D$ arbitrarily; however, this time with every application of $f$ the order relation between distinct $a_i,a_j$ changes if and only if the graph relation changes. This shows that $G$ contains $\Aut(D;T)$.

\textit{Case~4.} $f$ flips the graph relation on $X$ and flips it between $X$ and $Y$. Again, we proceed as in Case~1. So let $a_1,\ldots,a_k\in D$ and $1\leq i\leq k$ be as in that case, and define the function $h$ as before.
This time, $h$ flips the graph relation between $a_i$ and $a_j$ if and only if $j>i$, for all $1\leq j\leq k$. Let $\gamma\in\adel$ be so that it moves $h(a_i)$ into $X$, and the rest of the $h(a_j)$ into $Y$. Then $f\circ \gamma\circ h$ differs from $h$ in that it flips the graph relation between $a_i$ and $a_j$ if and only if $j<i$, for all $1\leq j\leq k$. In other words, it changes the graph relation between two elements in $\{a_1,\ldots,a_k\}$ if and only if it changes the order relation between them. As in Case~3 we conclude that $G$ contains $\Aut(D;T)$.
\end{proof}

\begin{proof}[Proof of Proposition~\ref{prop:unorderedCanonical}]
Assume first that $H$ contains $\Aut(D;<)$, $\Aut(D;E)$, or $\Aut(D;T)$; we claim that then the first case of the statement applies. To see this, note first that if $G$ does not contain any of these groups then we are done, so we may assume the contrary. Now observe that we know all the closed groups containing  one of the groups $\Aut(D;<)$, $\Aut(D;E)$, or $\Aut(D;T)$ from the reduct classifications for $(D;<), (D;E)$, and $(D;T)$; by our assumption, $G$ and $H$ are among them. It is therefore enough to check that any two such groups can be distinguished by elements in $\JI$. This is a mere checking of containment relations in Figure~\ref{fig:reducts} using the relational descriptions of the groups above $\Aut(D;<)$, $\Aut(D;E)$, and $\Aut(D;T)$ and left to the reader.

Now let $c_1,\ldots,c_n\in D$ with $c_1<\cdots <c_n$ and $f\colon\delcn\To \Delta$ be a canonical function which is generated by $H$ but not by $G$. By Lemma~\ref{lem:orderonlevel} we may assume that $f$ either preserves the order on all levels, or it reverses the order on all levels; for otherwise, $H$ contains ${\Aut(D;<)}$, $\Aut(D;E)$, or $\Aut(D;T)$ and we are done. In the latter case, $H$ contains either $\mix \id\rev$ or $\mix - \rev$; denote this function by $g$. If $g\nin G$ then we are done, so we may assume $g\in G$. Then $g\circ f$ is still not generated by $G$, and preserves the order on all levels; replacing $f$ by $g\circ f$ we may thus assume that $f$ preserves the order on all levels. In particular, orbits on the same level stay on the same level when applying $f$. 

Suppose there are $0\leq i<j\leq n$ such that $f$ sends level~$i$ above level~$j$ with respect to the order. Then let $(i,j)$ be the smallest such pair with respect to the lexicographic order. By Lemma~\ref{lem:orderbetweenlevels}, $H$ contains $\Aut(D;<)$, $\mix\id\turn$, or $\mix\sw\turn$. In the first case we are done, so assume the second or third case and denote the corresponding function by $g'$. As above, if $g'\nin G$ then we are done, so we may assume $g'\in G$. We may assume that the irrational number $\pi$ around which $g'$ turns lies between the images of the $i$-th and $j$-th level under $f$, and that levels are either sent entirely above or below $\pi$ under $f$. Then $g'\circ f$ is still canonical, is still not generated by $G$ because $g'\in G$, and if there are still $0\leq i'<j'\leq m$ as above for this new function, then $(i',j')$ is larger than $(i,j)$ in the lexicographic order. Hence repeating this process we may assume that $f$ preserves the order between all levels. Now suppose that there is $1\leq i\leq n$ such that $c_i$ is either sent by $f$ below level $i-1$ or above level $i$. Assuming without loss of generality the latter, consider the structure induced by the union of level~$i-1$, level~$i$, and $\{c_i\}$. By Lemma~\ref{lem:movesctotop} we are done. Thus, the only remaining possibility is that $f$ is order preserving, and so the second statement of the proposition holds.
\end{proof}

\section{The Ordered Case}
\label{sect:ordered}

\begin{defn}
We will denote the set of those groups in $\JI$ which preserve the order relation $<$ by $\JIo$. We moreover denote the set containing all joins of groups in $\JIo$ as well as $\adel$ (i.e., the set of groups below the label $j$ corresponding to the group $\Aut(D;<)$ in Figure~\ref{fig:reducts}) by $\OP$.
\end{defn}

The main goal of this section is to prove

\begin{prop}\label{prop:ordered}
Let $G,H\supseteq \adel$ be closed groups such that there exist $c_1,\ldots,c_n\in D$ and an order preserving canonical function $f\colon \delcn\To\Delta$ which is generated by $H$ but not by $G$. Then there exists an element of $\JIo$ which is contained in $H$ but not in $G$.
\end{prop}
Together with Proposition~\ref{prop:unorderedCanonical} from the preceding section, this will complete the proof of Proposition~\ref{prop:ji}; confer also the overview in Section~\ref{sect:overview}. The proof of Proposition~\ref{prop:ordered} will be given just before Section~\ref{sect:fitting}, assuming the truth of a certain statement (Lemma~\ref{lem:complete}). That statement will then be shown in Section~\ref{sect:fitting}.

\subsection{Homogeneity of the order preserving reducts}

In this subsection we will prove the following.

\begin{prop}\label{prop:homogeneous}
All groups in $\OP$ are automorphism groups of homogeneous structures in a 4-ary language.
\end{prop}

The homogeneity in a $4$-ary language will later be used in Proposition~\ref{lem:usehomogeneous} when we compare behaviors of canonical functions with groups in $\OP$.

$\JIo$ contains precisely five groups:
the small groups $\cl{\mix -\id}$, $\cl{\mix\sw\id}$, 
the medium-sized groups $\cl{\apl}$, $\cl{\apu}$, and
the large group $\Aut(D;<)$. 
For the two small groups and their join, the statement follows from
a general fact about reducts of superposed homogeneous structures (Section~\ref{sect:superposed-homogen}). 
The statement is most interesting 
for the group $\cl{\apl}$ (Proposition~\ref{prop:apl}), the group $\cl{\apu}$ (Proposition~\ref{prop:apu}), and the group $\cl{\apl,\apu}$ (Proposition~\ref{prop:apl-apu}). We finally show in Section~\ref{sect:op} that $\OP$ consists of precisely eight groups, 
and prove that also the remaining groups are homogeneous in a language of maximal arity four. 

%Since $(D;E)$ is isomorphic to $(V;E)$,
%we will use $\sw$ and $-$ as if they were defined
%on $D$ instead of $V$; this should not cause confusion. 

\subsubsection{Homogeneous reducts of superposed structures}
\label{sect:superposed-homogen}
Note that $\cl{\mix\sw\id}=\Aut(D;<,R^{(3)})$:
this follows from the description of $\Aut(D;R^{(3)})$
 in terms of $\sw$ from Section~\ref{sect:random-graph}
 and the fact that $E$ and $<$, and hence also $R^{3}$ and $<$, are freely superposed. 
To show homogeneity of $(D;<,R^{(3)})$, we will use the following.  

\begin{lem}\label{lem:superposed-homogen}
Let $\Delta$ be the free superposition of two homogeneous structures $\Gamma_1$ and $\Gamma_2$. We may assume that both $\Gamma_1$ and $\Gamma_2$ are reducts of $\Delta$.
Also suppose that $\Gamma_1$ has
a homogeneous reduct $\Gamma_1'$ and $\Gamma_2$ has a homogeneous reduct $\Gamma_2'$. Let $\Delta'$ be the reduct of $\Delta$ that contains precisely the relations of $\Gamma_1'$ and $\Gamma_2'$. Then $\Delta'$ is the free superposition of $\Gamma_1'$ and $\Gamma_2'$, and in particular homogeneous.
\end{lem}
\begin{proof}
This follows from the definition of superpositions, and a straightforward back-and-forth argument. 
\end{proof}

It is well-known that $(D;R^{(3)})$ is homogeneous; see e.g.~\cite{MacphersonSurvey}. Hence, Lemma~\ref{lem:superposed-homogen} implies the homogeneity of $(D;<,R^{(3)})$. 
We would like to use Lemma~\ref{lem:superposed-homogen} to also prove Proposition~\ref{prop:homogeneous} for the group $\cl{\mix-\id}$ and the group $\cl{\mix\sw\id,\mix-\id}$. 
Note that $\cl{\mix-\id} = \Aut(D;R^{(4)},<)$ and that $\cl{\mix\sw\id,\mix-\id} = \Aut(D;R^{(5)},<)$, by 
the results from Section~\ref{sect:random-graph}. 
Let $P$ be the relation $$\{(x,y,z) \in D^3 \; | \; (E(x,y) \wedge N(x,z)\wedge N(y,z)) \vee (N(x,y) \wedge E(x,z) \wedge E(y,z)) \} \; .$$

\begin{prop}\label{prop:arity-3}
We have $\Aut(D;R^{(4)}) = \Aut(D;P)$, and $(D;P)$ is homogeneous.
\end{prop} 
\begin{proof}
It is clear that $P$ is not preserved by $\sw$, but preserved by $-$, and hence, by the result of Thomas described in Section~\ref{sect:random-graph}, 
it follows that $\Aut(D;P)=\Aut(D;R^{(4)})$. 
To prove homogeneity of $(D;P)$,
let $(a_1,\dots,a_n),(b_1,\dots,b_n) \in D^n$, and let $\alpha$ be a partial isomorphism of $(D;P)$ with domain $A := \{a_1,\dots,a_n\}$ 
such that $\alpha(a_i) = b_i$ for all $i \leq n$.
We have to show that $\alpha$ can be extended to an automorphism
of $(D;P)$. For $n=1$ this is trivial. For $n \geq 2$, we first treat the case that $E(a_1,a_2)$ and $E(b_1,b_2)$. 
If for some $i \leq n$ the element $a_i$ is neither adjacent to $a_1$ nor to $a_2$, then $b_i$ is adjacent to neither $b_1$ nor $b_2$, since
$\alpha$ preserves the relation $P$. If $a_i$ is adjacent to $a_1$, but not to $a_2$,
then $b_i$ must be adjacent to $b_1$, but not to $b_2$, again because of preservation of $P$. Similarly, when $a_i$ is adjacent to $a_2$ but not to $a_1$
then $b_i$ must be adjacent to $b_2$ but not to $b_1$. The only
remaining situation is that $a_i$ is adjacent to both $a_1$ and $a_2$.
In this case also $b_i$ is adjacent to both $b_1$ and $b_2$, by the assumption that also $\neg P$ is preserved and reasoning as above. 
Now suppose that $i,j \leq n$. Since $E(a_1,a_i)$ iff $E(b_1,b_i)$
and $E(a_2,a_i)$ iff $E(b_2,b_i)$, and by preservation of $P$ and $\neg P$ we have that $E(a_i,a_j)$ iff $E(b_i,b_j)$. Hence, 
$\alpha$ must be a partial isomorphism of $(D;E)$ and by homogeneity of $(D;E)$ it
can be extended to an automorphism of $(D;E)$, and therefore also
to an automorphism of $(D;P)$. 

Next, consider the situation that $E(a_1,a_2)$ and $N(b_1,b_2)$. 
In this case, the map $x \mapsto -\alpha(x)$ defined on $A$ is a partial isomorphism of $(D;P)$ that satisfies the assumption above, and hence it can be extended to an automorphism $\beta$ of $(D;P)$. Then $-^{-1} \circ \beta$ is an automorphism of $(D;P)$ that extends $\alpha$. 
Similarly, when $N(a_1,a_2)$ we consider $x \mapsto \alpha(-x)$ 
% or $x \mapsto -\alpha(-x)$ 
and thereby reduce the argument to the situation above. 
%We will show that there are $\beta,\beta' \in \cl{\Aut(V;E) \cup \{-\}$
%such that $\beta \circ \alpha \circ \beta'$ is an isomorphism  
\end{proof}

As explained above, the following is a consequence of Proposition~\ref{prop:arity-3} and Lemma~\ref{lem:superposed-homogen}.
 
\begin{cor}
We have $\cl{\mix-\id} = \Aut(D;P,<)$, and $(D;P,<)$
is homogeneous.
\end{cor}
%\begin{proof}
%\end{proof}

%\begin{enumerate}
%\item %Let $\Gamma(3)$ be the reduct of $(V;E)$ that 
%contains all 4-ary relations that are preserved by $-$.
%Then 
%\item 
%Let $\Gamma(4)$ be the reduct of $(V;E)$ that
%contains all 4-ary relations that are preserved by $\sw$ and $-$.
Let $(a,b,c,d) \in D^4$ be such that $(a,d)$ is the only edge induced by $a,b,c,d$.
Let $Q$ be the smallest relation that contains $(a,b,c,d)$ and is preserved by $\Aut(D;E) \cup \{\sw,-\}$. 

\begin{prop}\label{prop:arity-4}
$\Aut(D;R^{(5)}) = \Aut(D;Q)$, and $(D;Q)$ is homogeneous. 
\end{prop}
\begin{proof}
The proof follows the same strategy as the proof of Proposition~\ref{prop:arity-3}. Again it follows from the results mentioned in Section~\ref{sect:random-graph} that $\Aut(D;R^{(5)}) = \Aut(D;Q)$. To prove the homogeneity of
$(D;Q)$, let $(a_1,\dots,a_n),(b_1,\dots,b_n) \in D^n$, and let $\alpha$ be a partial isomorphism of $(D;Q)$ with domain $\{a_1,\dots,a_n\}$ which maps $a_i$ to $b_i$ for all $i \leq n$. 
We have to show that $\alpha$ can be extended to an automorphism of $(D;Q)$. This is easy to see for $n \leq 3$. % since there is only one orbit of triples in $\Aut(D;R^{(5)})$. 
For $n > 3$, we first treat the case that 
$a_1,a_2,a_3$ and $b_1,b_2,b_3$ are independent sets in $(D;E)$, and that for each $i \leq n$ the element $a_i$ is adjacent to at most one of $a_1,a_2,a_3$ and $b_i$ is adjacent to at most one of $b_1,b_2,b_3$. If for some $i \leq n$ the element $a_i$ is adjacent to none of $a_1,a_2,a_3$, then $\neg Q(a_1,a_2,a_3,a_i)$, and therefore also $\neg Q(b_1,b_2,b_3,b_i)$, which implies that $b_i$ is adjacent to none of $b_1,b_2,b_3$. Moreover, preservation of $Q$ and $\neg Q$ implies that for all $i \leq n$ and $j \in \{1,2,3\}$ the element $a_i$ is adjacent to $a_j$ if and only if $b_i$ is adjacent to $b_j$. This fact and again preservation of $Q$ and $\neg Q$ imply that for all $i,j \leq n$, 
$a_i$ is adjacent to $a_j$ if and only if $b_i$ is adjacent to $b_j$. Hence, $\alpha$ must be a partial isomorphism of $(D;E)$,
and by homogeneity of $(D;E)$ it can be extended to an automorphism of $(D;E)$, and therefore also to an automorphism of $(D;Q)$. 

In the general case, when $a_1,a_2,a_3$ and
$b_1,b_2,b_3$ are not necessarily independent sets,
it is easy to see that we can choose $\beta_1,\beta_2 \in \Aut(D;Q)$ such that
$\gamma := \beta_1 \circ \alpha \circ \beta_2^{-1}$ is a partial isomorphism on $(D;Q)$ where 
\begin{itemize}
\item $\{\beta_2(a_1),\beta_2(a_2),\beta_2(a_3)\}$ forms an independent set;
\item $\{\beta_1(b_1),\beta_1(b_2),\beta_1(b_3)\}$ forms an independent set;
\item for $i \leq n$, $\beta_2(a_i)$
is adjacent to at most one of $\beta_2(a_1),\beta_2(a_2),\beta_2(a_3)$;
\item for $i \leq n$, $\beta_1(b_i)$
is adjacent to at most one of $\beta_1(b_1),\beta_1(b_2),\beta_1(b_3)$.
\end{itemize}
The previous paragraph implies that $\gamma$ can be extended to $\gamma' \in \Aut(D;E)$.
Then $\beta_1^{-1} \gamma' \beta_2$ is an automorphism of $(D;Q)$ which extends $\alpha$.
This shows homogeneity of $(D;Q)$. 
\end{proof}

The following is a consequence of Proposition~\ref{prop:arity-4} in combination with Lemma~\ref{lem:superposed-homogen}.
% that $\cl{\mix\sw\id,\mix-\id} = \Aut(D;Q,<)$ is homogeneous in a language of arity four. 

\begin{cor}
We have $\cl{\mix-\id,\mix\sw\id} = \Aut(D;Q,<)$, and $(D;Q,<)$
is homogeneous.
\end{cor}

\subsubsection{Homogeneous structures for the medium-sized groups}
\label{sect:medium-sized}
We now turn to the medium-sized group $\cl{\apl}$. 
Recall that $\apl$ is a permutation which preserves $<$
and switches the graph relation below some irrational $\pi$.
We now come to a statement that has already been announced in Section~\ref{sect:sporadic}.

\begin{prop}\label{prop:apl}
We have $\cl{\apl} = \Aut(D;R^l_3)$. The structure
$(D;R^l_3)$ is homogeneous. 
\end{prop}

To prove this proposition, we need some definitions and preparatory lemmas. 

\begin{definition}\label{def:lower-layered}
Let $S=\{v_1,\dots,v_n\} \subseteq D$ 
with $v_1 < \cdots < v_n$. We say that an order preserving map
$\delta \colon S \to D$ is \emph{lower layered} iff there exists a Boolean vector $t = (t_1,\dots,t_{n-1}) \in \{0,1\}^{n-1}$ such that 
for $i,j \in \{1,\dots,n\}$, $i<j$, we have that
$E(\delta(v_i),\delta(v_j))$ if and only if one of the following holds:
\begin{enumerate}
\item $E(v_i,v_j)$ and $\sum_{j-1 \leq k < n} t_k$ is even 
\item $N(v_i,v_j)$ and $\sum_{j-1 \leq k < n} t_k$ is odd.
\end{enumerate}
A permutation $\delta$ of $D$ is called 
\emph{lower layered} iff all its finite-range restrictions are lower layered. 
\end{definition}

An example of a permutation that is lower
layered is $\apl$.
It will be convenient to be slightly sloppy with
notation by identifying $0$ with $N$ and $1$ with $E$, so that the condition from Definition~\ref{def:lower-layered} can be rewritten to
$$E(\delta(v_i),\delta(v_j)) = E(v_i,v_j) + \sum_{j-1 \leq k < n} t_k \; ;$$
all arithmetic in this section is modulo 2.

%Analogously, we define \emph{lower layered} operations.  
% an example of a lower layered operation is $\apl$.
%Note that the composition of two lower layered
%operations 
%$\delta_1$, $\delta_2$
%is again lower layered: when $S \subseteq D$
%is of cardinality $n$, and for $i \in \{0,1\}$ the
%vector $t_i \in \{0,1\}^n$ witnesses that $\delta_i$ is lower layered on $S$, then
%$t_1 \oplus t_2$ witnesses that $\delta_1 \circ \delta_2$ is lower layered on $S$ as well.

\begin{lem}\label{lem:lower}
Let $S \subseteq D$ be finite, 
and $\gamma \colon S \rightarrow D$ be a partial isomorphism of $(D;R_3^l)$. Then $\gamma$ is lower layered.
\end{lem}
\begin{proof}
Our proof is by induction on the cardinality $n$
of $S$. The statement is clearly true for $n=1$ and $n=2$. 
Now suppose that the statement is true for $n \geq 2$. We want to show it for $n+1$. 
Write $S = \{v_1,\dots,v_n,v_{n+1}\}$ with
$v_1 < \cdots < v_n < v_{n+1}$. 
By induction hypothesis, there exists
$t \in \{0,1\}^{n-1}$ such that for all $i,j \in \{2,\dots,n+1\}$ with $i<j$
$$E(\gamma(v_i),\gamma(v_j)) = E(v_i,v_j) + \sum_{j-2 \leq k < n} t_k \; .$$
Define $t' \in \{0,1\}^{n}$ by $t'_i := t_{i-1}$ for $i \in \{2,\dots,n\}$, and 
$$t'_1 := E(\gamma(v_1),\gamma(v_2)) + E(v_1,v_2) + \sum_{1 \leq k < n} t_k \, .$$
%$, and $t' := (1,t_1,\dots,t_{n-1})$ otherwise. 
We claim that for $i,j \in \{1,\dots,n+1\}$, $i<j$, we have 
\begin{align}
E(\gamma(v_i),\gamma(v_j)) = E(v_i,v_j) + \sum_{j-1 \leq k < n+1} t'_k \; .
\label{eq:upper}
\end{align}
This is true for $i,j \in \{2,\dots,n+1\}$ by induction
assumption. Otherwise, $i = 1$. 
First consider the case that $j=2$.
By definition of $t_1'$ we have
$$E(\gamma(v_1),\gamma(v_2)) = E(v_1,v_2) +  t_1' + \sum_{1 \leq k < n} t_k \, ,$$ 
which equals $E(v_1,v_2) + \sum_{1 \leq k < n+1} t'_k$ by definition of $t_k'$ for $k > 1$.
Hence, Equation~(\ref{eq:upper}) is true in this case. 

Now consider the case that $j>2$.
Since $\gamma$ preserves $R_3^l$ and $\neg
R_3^l$ we have that
$$E(v_1,v_j) + E(v_2,v_j) = E(\gamma(v_1),\gamma(v_j)) + E(\gamma(v_2),\gamma(v_j)) \; .$$
Therefore, and by induction hypothesis we have 
\begin{align*}
E(\gamma(v_1),\gamma(v_j)) 
& = E(v_1,v_j) + E(v_2,v_j) + E(\gamma(v_2),\gamma(v_j)) \\
& = E(v_1,v_j) + E(v_2,v_j)+ E(v_2,v_j) + \sum_{j-2 \leq k < n} t_k \\
& = E(v_1,v_j) + \sum_{j-1 \leq k < n+1} t'_k \; ,
\end{align*} 
which is what we had to show. 
\end{proof}

\begin{lem}\label{lem:lower-2}
Let $S \subseteq D$ be finite,
and let $\gamma \colon S \to D$ be lower layered. 
Then $\gamma$ can be extended to an element
of $\cl{\apl}$. 
\end{lem}
\begin{proof}
Write $S = \{v_1,\dots,v_n\}$ with $v_1< \cdots < v_n$,
and let $t \in \{0,1\}^{n-1}$ be the witness that $\gamma$ is a lower layered permutation on $S$. 
We will show that there exists a sequence $\beta_1,\dots,\beta_n \in \Aut(D;E,<)$ and
a sequence $\alpha_1,\dots,\alpha_{n-1} \in \{\apl,\id\}$
such that $\gamma(x) = \beta_n \alpha_{n-1} \beta_{n-1}  \cdots \alpha_1 \beta_1(x)$ for all $x \in S$.

Choose a $\beta_1 \in \Aut(D;E,<)$ that maps $(v_1,v_2)$ below $\pi$ and $(v_3,\dots,v_n)$ above $\pi$. 
Choose $\alpha_1$ to be $\apl$ if $t_1=0$,
and the identity otherwise. 
To define $\beta_i$ for $2 \leq i \leq n-1$, suppose that $\beta_j$ has been defined for $j < i$, and
write $\gamma'$ for $\apl \beta_{i-1} \apl \cdots \apl \beta_1$. Then 
choose a $\beta_i \in \Aut(D;E,<)$ that maps $\gamma(v_1,\dots,v_{i+1})$  
below $\pi$ and $\gamma(v_{i+2},\dots,v_n)$ above $\pi$. Choose $\alpha_i$ to be $\apl$ if $t_i=0$,
and the identity otherwise. Then it is easy to verify that 
the restriction of $\alpha_{n-1} \beta_{n-1}  \cdots \alpha_1 \beta_1$ to $S$ is a partial isomorphism of $(D;E,<)$; by homogeneity of $(D;E,<)$ there exists a
$\beta_n \in \Aut(D;E,<)$ such that
$\gamma(x) = \beta_n \alpha_{n-1} \beta_{n-1}  \cdots \alpha_1 \beta_1(x)$ for all $x \in S$. 
\end{proof}

We finally give our operational reduct characterization. 

\begin{prop}\label{prop:lower}
Let $\gamma$ be a permutation of $D$. 
Then the following are equivalent. 
\begin{enumerate}
\item $\gamma$ is generated by $\apl$.
\item $\gamma \in \Aut(D;R^l_3)$.
\item $\gamma$ is lower layered. 
\end{enumerate}
\end{prop}
\begin{proof}
For the implication from $(1)$ to $(2)$ it suffices to verify that $\apl$ preserves $R^l_3$, which is straightforward. 
To show that $(2)$ implies $(3)$, 
let $\gamma \in \Aut(D;R^l_3)$ be arbitrary.
Lemma~\ref{lem:lower} shows that 
the restriction of $\gamma$ to finite subsets $S$ of $D$ is layered, and hence $\gamma$ is lower layered. 
The implication from $(3)$ to $(1)$ follows
from Lemma~\ref{lem:lower-2} and local closure. 
\end{proof}

\begin{proof}[Proof of Proposition~\ref{prop:apl}]
The equality $\Aut(D;R^l_3) = \cl{\apl}$ follows immediately from the equivalence of (1) and (2) in 
Proposition~\ref{prop:lower}. To show homogeneity of $(D;R^l_3)$,
let $\gamma$ be a partial isomorphism of $(D;R^l_3)$.
By Lemma~\ref{lem:lower},
$\gamma$ is lower layered. 
By Lemma~\ref{lem:lower-2}, $\gamma$ is the
restriction of an automorphism of $(D;R^l_3)$,
which is what we had to show. 
\end{proof}

Proposition~\ref{prop:apl} has the following dual version, which has also been announced in Section~\ref{sect:sporadic}, and which can be shown analogously. 

\begin{prop}\label{prop:apu}
We have $\cl{\apu} = \Aut(D;R^u_3)$. The structure
$(D;R^u_3)$ is homogeneous. 
\end{prop}

\subsubsection{A homogeneous structure for $\cl{\apl,\apu}$}
We now present a description of the join of the previous two medium-sized groups. 
First set
\begin{align*}
S_4 := & \{(a_1,\dots,a_4) \; | \; a_1 < a_2 < a_3 < a_4 \text{ and }  
E(a_1,a_3) + E(a_1,a_4) + E(a_2,a_3) + E(a_2,a_4) = 0 \}  \; .
\end{align*}

\begin{prop}\label{prop:apl-apu}
We have $\cl{\apl,\apu} = \Aut(D;S_4)$. The structure
$(D;S_4)$ is homogeneous.
\end{prop}

\begin{definition}\label{def:layered}
Let $S = \{v_1,\dots,v_n\} \subseteq D$ be of cardinality $n$.
A function $\delta \colon S \to D$ is called 
\emph{layered}
if there exist vectors $s,t \in \{0,1\}^{n-1}$ 
such that for $i,j \in \{1,\dots,n\}$, $i<j$, 
we have $$E(\delta(v_i),\delta(v_j)) = E(v_i,v_j) + 
\sum_{1 \leq k \leq i} s_k
+ \sum_{j-1 \leq k < n} t_k  \; .$$
\end{definition}

%It is easy to show that 
%the set of all
%layered operations are 
%also closed under composition. 

\begin{lem}\label{lem:layered}
Let $S \subseteq D$ be finite, and $\gamma \colon S \to D$ be an isomorphism between the substructures of $(D;S_4)$ induced by $S$ and by $\gamma(S)$.  
Then $\gamma$ is layered. 
\end{lem}
\begin{proof}
As in the proof of Proposition~\ref{prop:lower}, our proof is by induction on the cardinality $n$ of $S$. The statement is clearly true for $n \in \{1,2,3\}$. 
Now suppose that we have shown the claim
for $n \geq 3$, and that we want to show it for $n+1$. 
Write $S = \{v_1,\dots,v_{n+1}\}$ with
$v_1 < \cdots < v_{n+1}$.
By induction hypothesis, there exist
$s,t \in \{0,1\}^{n-1}$ such that for all $i,j \in \{1,\dots,n\}$ we have
$$E(\gamma(v_i),\gamma(v_j)) = E(v_i,v_j) + \sum_{1 \leq k \leq i} s_k + \sum_{j-1 \leq k < n} t_k \; .$$
We can assume without loss of generality that $s_1=0$
(otherwise, replace $t_{n-1}$ by $1+t_{n-1}$ and $s_1$ by $1+s_1$). Define $t' \in \{0,1\}^n$ by 
\begin{itemize}
\item $t'_n := E(v_1,v_{n+1})
+E(\gamma(v_1),\gamma(v_{n+1}))$, 
\item $t'_{n-1} := t_{n-1} + t_n'$,
and 
\item $t'_i = t_i$ for $i \in \{1,\dots,n-2\}$.
\end{itemize}
Furthermore, define $s' \in \{0,1\}^n$ by $s'_i := s_i$ for $i \in \{1,\dots,n-1\}$, and 
$$s'_n := E(\gamma(v_n),\gamma(v_{n+1})) + E(v_n,v_{n+1}) + \sum_{1 \leq k < n} s_k + t'_n \, .$$.

We claim that for $i,j \in \{1,\dots,n+1\}$, $i<j$,
we have
$$E(\gamma(v_i),\gamma(v_j)) = E(v_i,v_j) + 
 \sum_{1 \leq k \leq i} s'_k +
\sum_{j-1 \leq k < n+1} t'_k \; .$$
This is true for $i,j \in \{1,\dots,n\}$ since
\begin{align*}
E(\gamma(v_i),\gamma(v_j)) & = E(v_i,v_j) + \sum_{1 \leq k \leq i} s_k +
\sum_{j-1 \leq k < n} t_k \\
& = E(v_i,v_j) + 
\sum_{1 \leq k \leq i} s'_k 
+ t_n' + \sum_{j-1 \leq k < n} t'_k\\
& = E(v_i,v_j) + 
\sum_{1 \leq k \leq i} s'_k 
+ \sum_{j-1 \leq k < n+1} t'_k
\end{align*}
by induction assumption and the definition
of $s'$ and $t'$. 

Otherwise, $j=n+1$. 
Consider first the case that $i=1$. 
Then by definition of $t'_n$ we have
\begin{align}
E(\gamma(v_1),\gamma(v_{n+1})) = & E(v_1,v_{n+1})  + t_n'  \label{eq:onenplusone} \\
= & E(v_1,v_{n+1}) + \sum_{1 \leq k \leq 1} s'_k +
\sum_{n \leq k < n+1} t'_k  \nonumber
\end{align}
which is what we had to show. 

Next, consider the case that $i \in \{2,\dots,n-1\}$. 
Since $\gamma$ preserves $S_4$ and preserves $\neg S_4$ we have that
\begin{align}
& E(v_i,v_n) + E(v_i,v_{n+1}) + E(v_1,v_n) + E(v_1,v_{n+1}) \nonumber \\
= \quad &  
E(\gamma(v_i),\gamma(v_n)) + E(\gamma(v_i),\gamma(v_{n+1})) + E(\gamma(v_1),\gamma(v_n)) + E(\gamma(v_1),\gamma(v_{n+1})) \; . \label{eq:main-case}
\end{align} 
Since 
\begin{align*}
E(\gamma(v_i),\gamma(v_n)) = & E(v_i,v_n) + \sum_{1 \leq k \leq i} s'_k + \sum_{n-1 \leq k < n+1} t'_k \\
= & E(v_i,v_n) + \sum_{1 \leq k \leq i} s'_k + (t_{n-1} + t'_n) + t'_n \\
= & E(v_i,v_n) + \sum_{1 \leq k \leq i} s'_k + t_{n-1}
\end{align*}
and
\begin{align*}
E(\gamma(v_1),\gamma(v_n)) = & E(v_1,v_n) 
+ \sum_{1 \leq k \leq 1} s'_k +
\sum_{n-1 \leq k < n+1} t'_k \\
= & E(v_1,v_n) + t_{n-1}
\end{align*}
Equation~\ref{eq:main-case} simplifies to
\begin{align*}
& E(v_i,v_{n+1}) + E(v_1,v_{n+1}) \\
= \; &  
 E(\gamma(v_i),\gamma(v_{n+1})) + E(\gamma(v_1),\gamma(v_{n+1})) + \sum_{1 \leq k \leq i} s'_k \; .
\end{align*} 
Therefore,
\begin{align*}
E(\gamma(v_i),\gamma(v_{n+1})) = & E(v_i,v_{n+1}) + E(v_1,v_{n+1}) + \sum_{1 \leq k \leq i} s'_k + E(\gamma(v_1),\gamma(v_{n+1})) \\
= & E(v_i,v_{n+1}) + E(v_1,v_{n+1}) + \sum_{1 \leq k \leq i} s'_k + E(v_1,v_{n+1}) + t'_n  \\
%= & E(v_i,v_{n+1}) + \sum_{1 \leq k \leq 1} s'_k + 
%\sum_{n \leq k < n+1} t'_k + \sum_{1 \leq k \leq i} s'_k \\
%= & \dots \\
%= & E(v_i,v_{n+1}) + \sum_{1 \leq k \leq i} s'_k + t'_n \\
= & E(v_i,v_{n+1}) + \sum_{1 \leq k \leq i} s'_k +
\sum_{n \leq k < n+1} t'_k
\end{align*}

Finally, consider the case $i=n$.
By definition of $s'_n$ we have 
\begin{align*}
E(\gamma(v_n),\gamma(v_{n+1})) = & E(v_n,v_{n+1}) + \sum_{1 \leq k < n} s_k + s_n' + t'_n  \\
= & E(v_n,v_{n+1}) + \sum_{1 \leq k < n+1} s'_k + \sum_{n \leq k < n+1} t'_k
\end{align*}
and this concludes the induction.
% TODO: try the definition of s', t' above on examples
% of size 5
\end{proof}

\begin{lem}\label{lem:layered-2}
Let $S \subseteq D$ be finite,
and let $\gamma \colon S \to D$ be layered. 
Then $\gamma$ can be extended to an element of
$\cl{\apl,\apu}$. 
%Then there exists a sequence $\beta_1,\dots,\beta_{2n-1} \in \Aut(D;E,<)$ and
%a sequence of $\alpha_1,\dots,\alpha_{2n-2} \in \{\apl,\apu,\id\}$
%such that $\gamma(x) = \beta_{2n-1} \alpha_{2n-2} \beta_{2n-2}  \cdots \alpha_1 \beta_1(x)$ for all $x \in S$.
\end{lem}
\begin{proof}
This can be shown analogously to Lemma~\ref{lem:lower-2}. 
\end{proof}

\begin{prop}\label{prop:layered}
Let $\gamma$ be a permutation of $D$. 
Then the following are equivalent. 
\begin{enumerate}
\item $\gamma$ is generated by $\{\apl,\apu\}$. 
\item $\gamma \in \Aut(D;S_4)$.
\item $\gamma$ is a layered permutation. 
\end{enumerate}
\end{prop}
\begin{proof}
For the implication $(1) \Rightarrow (2)$ 
it suffices to show that both $\apl$ and $\apu$ 
preserve $S_4$, which is straightforward.
To prove $(2) \Rightarrow (3)$, 
let $\gamma \in \Aut(D;S_4)$ be arbitrary.
Lemma~\ref{lem:layered} shows that the restriction 
of $\gamma$ to finite subsets of $D$ is layered,
and hence $\gamma$ is layered. The implication
$(3) \Rightarrow (1)$ follows from Lemma~\ref{lem:layered-2} and local closure.
\end{proof}

\begin{proof}[Proof of Proposition~\ref{prop:apl-apu}]
The equality $\Aut(D;S_4) = \cl{\apl}$ follows immediately from the equivalence of (1) and (2) in 
Proposition~\ref{prop:layered}. To show homogeneity of $(D;S_4)$,
let $\gamma$ be a partial isomorphism of $(D;S_4)$.
By Lemma~\ref{lem:layered},
$\gamma$ is layered. 
By Lemma~\ref{lem:layered-2}, the function $\gamma$ is the restriction of an automorphism of $(D;S_4)$,
which is what we had to show. 
\end{proof}

\subsubsection{All groups in $\OP$}
\label{sect:op}
\begin{proof}[Proof of Proposition~\ref{prop:homogeneous}]
It is clear that $\cl{\mix -\id}$ is contained in $\cl{\apl}$,
and analogously that it is contained in $\cl{\apu}$.
Moreover, it is clear that $\cl{\mix\sw\id,\apl}$ 
contains $\cl{\apu}$, since $\mix-\id \circ \mix\sw\id \circ \apl$ (here, we assume that the irrational number $\pi$ used to define $\apl$ equals the irrational number used to define $\mix\sw\id$) behaves as $\apu$. Dually, $\cl{\mix\sw\id,\apu}$ 
contains $\cl{\apl}$. Finally, $\cl{\apl,\apu}$ contains
$\mix\sw\id$ since $\mix-\id \circ \apl \circ \apu$ behaves
as $\mix\sw\id$. 
All these groups are contained in $\Aut(D;<)$, and
contain $\Aut(\Delta)$.
This shows that 
$\OP$ consists of precisely the groups $\Aut(\Delta)$, 
$\cl{\mix\sw\id}$, $\cl{\mix-\id}$, $\cl{\mix\sw\id,\mix-\id}$, 
$\cl{\apl}$, $\cl{\apu}$, $\cl{\apl,\apu}$, and
$\Aut(D;<)$. We verify for each of those groups $G$ that there exists a homogeneous structure $\Gamma$ with an at most 4-ary language such that $\Aut(\Gamma) = G$.
\begin{itemize}
\item Since $(D;<)$ and $(D;E,<)$ are homogeneous, we are done for those groups.
\item For the groups $\cl{\apl}$, $\cl{\apu}$, 
$\cl{\apl,\apu}$, the statement has been shown in 
Proposition~\ref{prop:apl}, Proposition~\ref{prop:apu}, and Proposition~\ref{prop:apl-apu}.
\item The groups $\cl{\mix-\id}$, $\cl{\mix\sw\id}$, and $\cl{\mix\sw\id,\mix-\id}$ have been treated in Section~\ref{sect:superposed-homogen}. 
\end{itemize}
This concludes the proof that all groups in $\OP$ are
automorphisms groups of homogeneous structures in an at most 4-ary language.  
\end{proof}

Note that when a structure is homogeneous in an at most 4-ary language, then it is also homogeneous in a 4-ary language (using dummy variables).
%(we can simply repeat the last argument of a relation of smaller arity several times to increase the arity of a relation without changing the automorphism group, and without affecting homogeneity). 
From now on, we will therefore use the shorter phrase `homogeneous in a 4-ary language' instead of `homogeneous in an at most 4-ary language'.

\subsection{Moving above $\Aut(D;<)$ by order preserving behaviors}

\begin{lem}\label{lem:idorminusbetween}
Let $G\supseteq \adel$ be a closed group generating an order preserving canonical function $f\colon\delcn\To\Delta$, where $c_1,\ldots,c_n\in D$. Suppose that there exist infinite orbits $X,Y$ of $\delcn$ satisfying $\neg(Y< X)$ such that $f$ sends all pairs $(x,y)$, where $x\in X$, $y\in Y$, and $x<y$, to edges, or, dually, to non-edges. Then $G\supseteq \Aut(D;<)$.
\end{lem}
\begin{proof}
If $f$ eradicates edges or non-edges on $Y$ then we are done by Lemma~\ref{lem:eradicatesonorbit}, so we may assume that it keeps or flips the graph relation on $Y$. If $f$ flips the graph relation on $Y$, then pick any self-embedding $e$ of $\Delta$ whose range is contained in $Y$; replacing $f$ by $f\circ e\circ f$ we then have a function which still satisfies the assumptions of the lemma, and which keeps the graph relation on $Y$.

Let $a_1,\ldots ,a_k$ be elements of $D$ with $a_1<\cdots < a_k$, let $1\leq i\leq k$, and let $\alpha\in\adel$ be such that $\alpha[\{a_1, \ldots, a_i\}]\subseteq X$ and $\alpha[\{a_{i+1}, \ldots, a_k\}]\subseteq Y$. Then $f\circ \alpha$ is order preserving, keeps the graph relation on $\{a_{i+1}, \ldots, a_k\}$, and eradicates non-edges between $\{a_1, \ldots, a_i\}$ and $\{a_{i+1}, \ldots, a_k\}$. By applying this step iteratively starting with $i=k$ and finishing with $i=1$, we have that the set $\{a_1, \ldots, a_k\}$ is mapped to a complete graph by a function generated by $f$. Thus, $G$ contains $\Aut(D; <)$ by Lemma~\ref{lem:eradicatesonfinite}.
\end{proof}

\begin{lem}\label{lem:differentbehavior}
Let $G\supseteq \adel$ be a closed group generating an order preserving canonical function $f\colon\delcn\To\Delta$, where $c_1,\ldots,c_n\in D$. 
  Assume that $f$ keeps the graph relation on $D\setminus\{c_1,\ldots,c_n\}$. Assume moreover that there exist infinite orbits $X,Y$ of $\delcn$ on the same level and $1\leq i\leq n$ such that $f$ keeps the graph relation between $\{c_i\}$ and $X$, and flips the graph relation between $\{c_i\}$ and $Y$. Then $G\supseteq \Aut(D;<)$.
\end{lem} 
\begin{proof}
Without loss of generality assume that $\{c_i\}< X$ (and hence $\{c_i\}< Y$). Assume first that there are edges between $c_i$ and all elements of $X\cup Y$. Then let $Z$ be an infinite orbit at the level of $X$ and $Y$ whose elements are not adjacent to $c_i$. 
Because $f$ is canonical, it keeps or flips the graph relation between $\{c_i\}$ and $Z$. 
Replacing $X$ by $Z$ in the first case and $Y$ by $Z$ in the latter case, we then have that $c_i$ is adjacent to the elements of precisely one of the sets $X$ and $Y$. Arguing dually in the case where $c_i$ is adjacent to none of the elements in $X\cup Y$, we obtain the same situation. Without loss of generality, we henceforth assume that $c_i$ is adjacent to the elements of $X$, and not adjacent to those of $Y$.
 
 Let $a_1,\ldots ,a_k$ be elements of $D$ such that $a_1<\cdots < a_k$, and let $1\leq j\leq k$.
 There is a permutation $\delta\in\Aut(\Delta)$ that maps $a_j$ to $c_i$ and $\{a_{j+1}, \ldots, a_{k}\}$ into $X\cup Y$ (by our assumption above, those elements adjacent to $a_j$ must go into $X$, and the others into $Y$). Then $f\circ \delta$ does not modify the graph relation on $\{a_{j+1}, \ldots, a_{k}\}$, and eradicates non-edges between $\{a_j\}$ and $\{a_{j+1}, \ldots, a_{k}\}$. Hence, applying such functions from $j=k$ until $j=1$, $\{a_{1}, \ldots, a_{k}\}$ can be mapped to an independent set, and the lemma follows from Lemma~\ref{lem:eradicatesonfinite}.
 \end{proof}

\begin{lem}\label{lem:2levelsaboveconstant}
Let $G\supseteq \adel$ be a closed group generating an order preserving canonical function $f\colon\delcn\To\Delta$, where $c_1,\ldots,c_n\in D$. 
  Assume that $f$ keeps the graph relation on $D\setminus\{c_1,\ldots,c_n\}$. Assume moreover that there exist distinct levels $L_1,L_2,L_3$ of $\delcn$ and $1\leq i\leq n$ with $L_1<\{c_i\}<L_2, L_3$ such that $f$ keeps the graph relation between $\{c_i\}$ and $L_1\cup L_3$, and such that $f$ flips the graph relation between $\{c_i\}$ and $L_2$. Then $G\supseteq \Aut(D;<)$.
\end{lem}
\begin{proof}
 Let $\gamma\in \adel$ be such that $\gamma(c_i)=c_i$, everything above $c_i$ is mapped into $L_2$, and everything below $c_i$ is mapped into $L_1$. Then $f\circ \gamma\circ f$ is a canonical function $\delcn\rightarrow \Delta$ generated by $G$, keeps the graph relation between $c_i$ and $L_1$, and switches the role of $L_2$ and $L_3$. Hence, we may assume that $L_2<L_3$. Let $A$ be an arbitrary finite subset of $D$ that consists of the elements $a_1<\cdots <a_k$. Let $1\leq j<m\leq k$. There is a permutation $\delta\in\Aut(\Delta)$ that maps $a_j$ to $c_i$, $\{a_{1}, \ldots, a_{j-1}\}$ into $L_1$, $\{a_{j+1}, \ldots, a_{m-1}\}$ into $L_2$, and $\{a_{m}, \ldots, a_{k}\}$ into $L_3$. Then $f\circ \delta$ flips the graph relation between a pair of elements in $A$ if and only if one of them is $a_j$ and the index $r$ of the other is such that $j<r<m$. Applying the same modification with the pair $(j, m+1)$ instead of $(j,m)$, these two modifications combined flip the graph relation between $a_j$ and $a_m$, and nowhere else in $A$. Hence, $A$ can be mapped to an independent set in a finite number of such steps, and the lemma follows from Lemma~\ref{lem:eradicatesonfinite}.
\end{proof}

\subsection{Moving above $\apu$ and $\mix\sw\id$}
\begin{lem}\label{lem:apu}
Let $G\supseteq \adel$ be a closed group, and let $f\colon\delcn\To\Delta$ be an order preserving canonical function generated by $G$, where $c_1,\ldots,c_n\in D$. 
 Assume that there exists $1\leq i\leq n$ and levels $X,Y$ of $\delcn$ with $X<\{c_i\}<Y$ such that $f$ keeps the graph relation on $X\cup Y$ as well as between $X$ and $\{c_i\}$, and such that $f$ flips the graph relation between $\{c_i\}$ and $Y$. Then $G$ contains $\apu$.
\end{lem}
\begin{proof}
Let $A$ be an arbitrary finite subset of $D$ that consists of the elements $a_1<\cdots <a_k$. Let $1\leq j\leq k$ be arbitrary. There exists $\delta \in\adel$ which maps $\{a_1,\ldots,a_{j-1}\}$ into $X$, $a_j$ to $c_i$, and $\{a_{j+1},\ldots,a_{k}\}$ into $Y$. Hence, application of $f\circ \delta$ flips the graph relation between $\{a_j\}$ and $\{a_{j+1},\ldots,a_{k}\}$, and keeps the graph relation otherwise on $A$. Now if we fix any $1\leq m\leq k$ and apply this process iteratively to all $m\leq j\leq k$, we get that the graph relation is flipped on $\{a_{m},\ldots,a_{k}\}$, kept on $\{a_1,\ldots,a_{m-1}\}$, and kept between $\{a_1,\ldots,a_{m-1}\}$ and $\{a_{m},\ldots,a_{k}\}$. Hence, $\apu$ is generated by $f$.
\end{proof}

\begin{lem}\label{lem:said}
Let $G\supseteq \adel$ be a closed group, and let $f\colon\delcn\To\Delta$ be an order preserving canonical function generated by $G$, where $c_1,\ldots,c_n\in D$. 
 Assume that there exists $1\leq i\leq n$ such that $f$ keeps the graph relation on $D\setminus\{c_1,\ldots,c_n\}$, and such that $f$ flips the graph relation between $\{c_i\}$ and $D\setminus\{c_1,\ldots,c_n\}$. Then $G$ contains $\mix\sw\id$.
\end{lem}
\begin{proof}
There exists a self-embedding $e$ of $\Delta$ which fixes $c_i$ and whose range is contained in $D\setminus\{c_1,\ldots,c_{i-1},c_{i+1},\ldots,c_n\}$. Then $f\circ e$ flips edges and non-edges between $\{c_i\}$ and its complement, and keeps the edge relation on this complement.
\end{proof}
  
\subsection{Fitting behaviors with groups}

\begin{definition}\label{def:sec6maindef}
Let $c_1,\ldots,c_n\in D$, and let $B$ be a behavior between $\delcn$ and $\Delta$. Let moreover $C\supseteq \adel$ be a closed group. We say that
\begin{itemize}
\item $B$ \emph{forces} $C$ iff every canonical function $f\colon\delcn\To\Delta$ satisfying $B$ generates all functions in $C$;
\item $B$ is \emph{compatible} with $C$ iff there exists a canonical function $f\colon\delcn\To\Delta$ satisfying $B$ which preserves all relations invariant under $C$;
\item $B$ \emph{fits} $C$ iff $B$ forces $C$ and $B$ is compatible with $C$.
\end{itemize}

Now let $\S$ be a set of closed groups above $\adel$. Then we say that 
\begin{itemize}
\item $B$ is $\S$-\emph{fittable} iff there exists a closed group $C$ in $\S$ such that $B$ fits $C$;
\item $B$ is \emph{fittable} iff there exists a closed group $C\supseteq \adel$ such that $B$ fits $C$.
\end{itemize}
\end{definition}

We will use the above terminology mainly for specific sets of type conditions satisfied by canonical functions from $\delcn$ to $\Delta$.

\begin{definition}\label{def:sec6maindef2}
Let $c_1,\ldots,c_n\in D$, let $f\colon \delcn\To \Delta$ be canonical, and let $k\geq 1$. 
Then we call any restriction of the behavior of $f$ (i.e., the set of type conditions satisfied by $f$) to types involving $k$ fixed infinite $1$-types a \emph{$k$-constellation} of $f$.
\end{definition}

For example, if $X$ and $Y$ are infinite orbits of $\delcn$, then the set of all type conditions satisfied by $f$ which say something about its behavior on $X\cup Y$ is a 2-constellation of $f$. 

\begin{lem}\label{lem:compatible}
Let $c_1,\ldots,c_n\in D$, let $f\colon \delcn\To \Delta$ be canonical, and let $B$ be a $k$-constellation of $f$, where $k\geq 1$. Then $B$ is compatible with a given closed group $C\supseteq \adel$ if and only if it preserves all relations invariant under $C$ on the union of the orbits concerned by $B$.
\end{lem}
\begin{proof}
Write $S$ for the union of the orbits concerned by $B$, and $\Gamma$ for the structure on $D$ whose relations are precisely those invariant under $C$. Clearly, if $B$ is compatible with $C$, then $f$ preserves all relations of $\Gamma$ on $S$. Assume now the latter; we have to show that there exists a canonical function $g\colon \delcn\To \Delta$ which preserves all relations of $\Gamma$ and satisfies $B$. Let $T$ be the set of finite partial functions on $D$ which preserve all relations of $\Gamma$ and satisfy $B$. We claim that every finite $A\subseteq D$ is the domain of a function in $T$. To see this, consider the restriction of $f$ to $A\cap S$, which is a partial isomorphism on $\Gamma$; by the homogeneity of $\Gamma$, it extends to an automorphism $\alpha$ of $\Gamma$, and the restriction of $\alpha$ to $A$ is an element of $T$. Note that if $t\in T$ and $\beta\in\adel$, then $\beta\circ t\in T$. Hence, by the same standard compactness argument needed to prove Lemma~\ref{lem:compactness}, there exists a function $h\colon D\To D$ whose restriction to any finite subset of $D$ is an element of $T$. We then have that $h$ preserves the relations of $\Gamma$ and satisfies $B$. By Proposition~\ref{prop:ramseystrong},
$$
\overline{\{\alpha\circ h\circ \beta\;|\; \alpha\in \adel,\; \beta\in \adelcn\}}
$$
contains a canonical function $g$ which obviously still preserves the relations of $\Gamma$ and satisfies $B$, proving the lemma.
\end{proof}

The remainder of this section is devoted to the proof of the following lemma, which, as we will show immediately, implies Proposition~\ref{prop:ordered}. Recall that together with Proposition~\ref{prop:unorderedCanonical}, which we proved in Section~\ref{sect:unordered}, we then obtain a complete proof of Proposition~\ref{prop:ji}.

\begin{lem}\label{lem:complete}
Let $1\leq k\leq 4$, and let $B$ be a $k$-constellation of an order preserving canonical function $f\colon \delcn\To \Delta$, where $c_1,\ldots,c_n\in D$. Then $B$ is $\OP$-fittable.
\end{lem}

In the following lemma, we are given a canonical function $f\colon \delcn\To \Delta$ and a closed group $C$ which can locally invert the effect of $f$ on a fixed collection of infinite orbits; the lemma roughly states that we can then produce from $f$ and $C$ a canonical function which behaves like the identity function on the collection.

\begin{lem}\label{lem:invertonconstellation}
Let 
\begin{itemize}
\item $C\supseteq \adel$ be a closed group;
\item $f\colon \delcn\To \Delta$ be a canonical function, where $c_1,\ldots,c_n\in D$; 
\item $B$ be a $k$-constellation of $f$ which is compatible with $C$.
\end{itemize}
Then $$\overline{\{\gamma\circ f\circ \beta\;|\; \gamma\in  C,\; \beta\in \adelcn\}}$$ contains a canonical function $g\colon \delcn\To \Delta$ whose restriction to the union of the orbits concerned by $B$ behaves as the identity function does (i.e., that restriction is a partial self-embedding of $\Delta$).
\end{lem}
\begin{proof}
Let $S$ be the union of the infinite orbits of $\delcn$ concerned by $B$. 
Let $\Gamma$ be the structure whose relations are those preserved by the permutations in $C$. Clearly, $\Gamma$ is homogeneous, and the restriction of $f$ to any finite subset of $S$ is a partial isomorphism of $\Gamma$ whose inverse extends to an element of $C$. 
Consider the set $T$ of all partial functions $p$ such that 
\begin{itemize}
\item $p$ is the restriction of a function of the form $\gamma\circ f$, where $\gamma\in C$, to some finite subset of $D$;
\item $p$ behaves as the identity function does on the intersection of its domain with $S$.
\end{itemize}
By the above observation, every finite subset of $D$ is the domain of a function in $T$. Since $\Gamma$ is $\omega$-categorical, by Lemma~\ref{lem:compactness} there exists a function $h$ in $\overline{\{\gamma\circ f\;|\; \gamma\in  C\}}$ which behaves like the identity function on $S$. Now Proposition~\ref{prop:ramseystrong} tells us that $$
\overline{\{\alpha\circ h\circ \beta\;|\; \alpha\in\adel, \beta\in \adelcn\}}
$$ contains a canonical function $g\colon \delcn\To \Delta$; clearly, the restriction of this function to $S$ still behaves as the identity function does.
%Write $D=\{d_0,d_1,\ldots\}$, and consider the set $T$ of all partial functions $p$ such that 
%\begin{itemize}
%\item $p$ is the restriction of a function of the form $\gamma\circ f$, where $\gamma\in C$, to $\{d_0,\ldots,d_i\}$ for some $i\in\omega$;
%\item $p(s)=s$ for all $s\in S\cap\{d_0,\ldots,d_i\}$.
%\end{itemize}
%By the above observation, $T$ contains functions of arbitrarily large domain. Now consider $p,p'\in T$ equivalent iff $p$ and $p'$ have equal domain $\{d_0,\ldots,d_i\}$, and the type of $(p(d_0),\ldots,p(d_i))$ equals the type of $(p'(d_0),\ldots,p'(d_i))$ in $\Delta$. For $p\in T$, write $[p]$ for the equivalence class of $p$ with respect to this equivalence relation. Setting $[p]\leq[p']$ iff there exists $q\in [p']$ extending $p$ then defines a partial order on the equivalence classes, which is in fact a finitely branching infinite tree. Picking an infinite branch $B$, we can then find a sequence of $(p_i)_{i\in\omega}$ of elements of $T$ such that $[p_i]\in B$ and the domain of $p_i$ equals $\{d_0,\ldots,d_i\}$ for all $i\in\omega$, and such that $p_j$ extends $p_i$ whenever $i<j$. The union of the $p_i$ is a function $h$ in $\overline{\{\gamma\circ f\;|\; \gamma\in  C\}}$ which behaves like the identity function on $S$. Now Proposition~\ref{prop:ramseystrong} tells us that $\overline{\{\alpha\circ h\circ \beta\;|\; \alpha\in\adel, \beta\in \adelcn\}}$ contains a canonical function $g\colon \delcn\To \Delta$; clearly, the restriction of this function to $S$ still behaves like the identity function.
\end{proof}

\begin{definition}\label{def:subconstellation}
Let $\Lambda, \Omega$ be structures, and let $B$ be a behavior between $\Lambda$ and $\Omega$. Then we call any subset of $B$ a \emph{subbehavior} of $B$. If $B$ is a $k$-constellation, and $B'\subseteq B$ is an $m$-constellation, then we will call $B'$ an \emph{$m$-subconstellation} of $B$.
\end{definition}

\begin{lem}\label{lem:usehomogeneous}
Let $C\supseteq \adel$ be the automorphism group of a homogeneous structure in an $m$-ary relational language, where $m\geq 1$. Let $f\colon \delcn\To \Delta$ be a canonical function, where $c_1,\ldots,c_n\in D$. If $k\geq m$, then a $k$-constellation of $f$ is compatible with $C$ if and only if its $m$-subconstellations are compatible with $C$.
\end{lem}
\begin{proof}
Let $B$ be any $k$-constellation of $f$, where $k\geq m$. If $B$ is compatible with $C$, then trivially so are its $m$-subconstellations. For the converse, let $S$ be the union of the orbits concerned by $B$; we claim that $f$ preserves all relations invariant under $C$ on $S$. Otherwise, $f$ would violate an $m$-ary relation $R$ invariant under $C$, since $C$ is homogeneous in an $m$-ary relational language. So there would be a tuple $(a_1,\ldots,a_m)$ of elements in $S$ such that $(a_1,\ldots,a_m)\in R$ and $(f(a_1),\ldots,f(a_m))\notin R$. Writing $O_i$ for the orbit of $a_i$, for all $1\leq i\leq m$, we would then have that any $m$-constellation of $f$ concerning all of the $O_i$ would be incompatible with $C$, a contradiction. Hence, $f$ indeed preserves all relations invariant under $C$ on $S$, and so we are done by Lemma~\ref{lem:compatible}.
\end{proof}

We are now ready to prove Proposition~\ref{prop:ordered} (assuming the truth of Lemma~\ref{lem:complete}).

\begin{proof}[Proof of Proposition~\ref{prop:ordered}]

 Let $M$ be the largest group in $\OP$ that is contained in $G$; in other words, $M$ is the join of $\adel$ and all groups in $\JIo$ contained in $G$. Let $\Gamma$ be the homogeneous structure in a $4$-ary language such that $M=\Aut(\Gamma)$, guaranteed by Proposition~\ref{prop:homogeneous}. 
 
 Suppose that $f$ has a $4$-constellation $B$ which is not compatible with $M$. By Lemma~\ref{lem:complete}, $B$ is $\OP$-fittable; let $C\in \OP$ be so that $B$ fits $C$. Since $B$ forces $C$, we have that $H$ contains $C$. On the other hand, since $B$ is compatible with $C$, we have that $C$ cannot be contained in $M$:  otherwise, since $B$ is not compatible with $M$, it would not be compatible with  $C$ either. By the definition of $M$, we conclude that $C$ is not contained in $G$ either, and so we are done.
 
We may thus assume that all $4$-constellations of $f$ are compatible with $M$; by Lemmas~\ref{lem:usehomogeneous} and~\ref{lem:invertonconstellation}, there is a canonical function $g\colon \delcn\To \Delta$ which is not generated by $G$ and whose restriction to $D\setminus\{c_1,\ldots,c_n\}$ behaves like the identity function. By canonicity, we have that $g$ keeps or flips the graph relation between any $\{c_i\}$ and any infinite orbit; using Lemma~\ref{lem:differentbehavior} we may even assume that $g$ keeps or flips the graph relation between any $\{c_i\}$ and any level of $\delcn$.
 
 Assume that there exist $1\leq i\leq n$ and two levels such that $g$ keeps the graph relation between $\{c_i\}$ and one of them, and $g$ flips the graph relation between $\{c_i\}$ and the other. There exist two levels $L_1, L_2$ such that $L_1<\{c_i\}<L_2$, such that $g$ keeps the graph relation between $\{c_i\}$ and one of them, and such that $g$ flips the graph relation between $\{c_i\}$ and the other. Assume without loss of generality that $g$ keeps the graph relation between $\{c_i\}$ and $L_1$. Then $H$ contains $\apu$ by Lemma~\ref{lem:apu}. If $\apu\not\in G$ then we are done, so we may assume that $\apu\in G$. Then it is easy to see that $G$ generates a canonical function $h \colon (\Delta,c_i)\To\Delta$ which preserves the order, which flips the graph relation between $g(c_i)$ and all $d\in D$ with $g(c_i)< d$, and which keeps the graph relation otherwise. By Lemma~\ref{lem:2levelsaboveconstant}, we may assume that $g$ keeps the graph relation between a level $L$ and $\{c_i\}$ if and only if $L<\{c_i\}$.
 Replacing $g$ by $h\circ g$ we then may assume that $g$ keeps the graph relation between $\{c_i\}$ and the union of all infinite orbits. Repeating this process for all $1\leq i\leq n$, we then have that for each $i$ the function $g$ keeps or flips the graph relation between $\{c_i\}$ and the union of all infinite orbits.  

 Assume now that there exists $1\leq i\leq n$ such that $g$ flips the graph relation between $\{c_i\}$ and the union of all infinite orbits. Then $H$ contains $\mix \sw \id$ by Lemma~\ref{lem:said}, and we may assume that $G$ contains this function as well, for otherwise we are done. Then it is easy to see that $G$ generates a canonical function $h \colon (\Delta,c_i)\To\Delta$ which preserves the order, which flips the graph relation between $g(c_i)$ and its complement, and which keeps the graph relation otherwise. Replacing $g$ by $h\circ g$ we then may assume that $g$ keeps the graph relation between $\{c_i\}$ and the union of all infinite orbits. Repeating this process for all $1\leq i\leq n$, we then have that for each $i$ the function $g$ keeps or flips the graph relation between $\{c_i\}$ and the union of all infinite orbits.
 
We now continue with this assumption, i.e., if $x,y\in D$ are so that the edge relation is altered between $x$ and $y$, then $x,y\in\{c_1,\ldots,c_n\}$.
As $g$ violates a relation definable in $\Delta$, there exist $1\leq i<j\leq n$ such that the graph relation is flipped between $\{c_i\}$ and $\{c_j\}$. Using Lemma~\ref{lem:eradicatesonfinite}, it is then straightforward to see that $H$ contains $\Aut(D;<)$.
\end{proof}

\subsection{Fitting the constellations}\label{sect:fitting}
We will now prove Lemma~\ref{lem:complete} by considering all possible $k$-constellations $B$, for $1\leq k\leq 4$. For every such $B$ we find a group $C\in\OP$ such that $B$ fits $C$. From $k= 2$ on, we will draw a picture describing $B$, and the corresponding group $C$ will be indicated in the lower right corner. The proof that a constellation $B$ fits a group $C$ consists of two parts: proving that $B$ forces $C$ and proving that $B$ is compatible with $C$. Verifying compatibility can be automated using Lemma~\ref{lem:compatible}, and we shall omit these verifications: for example, to verify that $R^{(3)}$ is preserved, one has to consider all 3-element ordered graphs, distribute their vertices in all possible ways among the orbits, and check that the modification of the graph relations according to the behavior really agrees with $R^{(3)}$ in all these cases.

Of course, the number of $k$-constellations grows with $k$, and we will be obliged to systematically use knowledge on smaller constellations (i.e., $(k-1)$-subconstellations). An example of this is the following easy but useful observation concerning constellations which fit the specific group $\Aut(D;<)$; since they appear often in our analysis, it is convenient to give them their own name.

\begin{defn}
We call a $k$-constellation of an order preserving canonical function $f\colon \delcn\To\Delta$ \emph{full} iff it fits $\Aut(D;<)$.
\end{defn}

\begin{lem}\label{full}
Let $f\colon\delcn\To\Delta$ be canonical and order preserving. Then any $k$-constellation of $f$ with a subconstellation forcing $\Aut(D;<)$ also fits $\Aut(D;<)$ (and in particular, is full and $\OP$-fittable).
\end{lem}
\begin{proof}
Clearly such a constellation forces $\Aut(D;<)$; because $f$ is order preserving, it also is compatible with $\Aut(D;<)$.
\end{proof}

As a consequence, once we know that a certain $(k-1)$-constellation forces $\Aut(D;<)$, we only have to consider $k$-constellations which do not have this constellation as a subconstellation, reducing the number of cases.

Another lemma that reduces the number of cases is the following.

\begin{lem}\label{lem:reduction}
Let $f\colon\delcn\To\Delta$ be an order preserving canonical function, and let $B$ be a $k$-constellation of $f$. Let $S$ be the union over all orbits concerned by $B$. Suppose that $B$ has a subconstellation fitting a group $C$, suppose there is a canonical function in 
$$\overline{\{\gamma\circ f\circ \beta \;|\; \gamma\in C,\; \beta\in\Aut(\Delta,c_1,\ldots,c_n)\}}$$ whose restriction to $S$ has an $\OP$-fittable behavior $B'$. Then $B$ is $\OP$-fittable.
\end{lem}
\begin{proof}
Let $C'$ fit $B'$. Clearly, $B$ forces $C\vee C'$. On the other hand, suppose that $R$ is a relation which is invariant under $C\vee C'$ and which is violated by $f$ on $S$. Then note that all functions of the form $f\circ\beta$, where $\beta\in\Aut(\Delta,c_1,\ldots,c_n)$, violate $R$ on precisely the same sets as $f$, since they all have the same behavior. Because $R$ is invariant under $C$, even all functions of the form $\gamma\circ f\circ\beta$, where $\gamma\in C$ and $\beta\in \Aut(\Delta,c_1,\ldots,c_n)$,  violate $R$ on precisely the same sets as $f$. Hence, any function in the closure of the set of such functions violates $R$ on $S$. But this contradicts the assumption that $R$ is invariant under $C'$, proving that $B$ is compatible with $C\vee C'\in \OP$.
\end{proof}

\begin{lem}\label{lem:identityonsub}
Assume that all $(k-1)$-constellations of order preserving canonical functions are $\OP$-fittable. Assume moreover that all $k$-constellations of order preserving canonical functions which have a $(k-1)$-subconstellation with the behavior of the identity function are $\OP$-fittable. Then all $k$-constellations of order preserving canonical functions are $\OP$-fittable.
\end{lem}
\begin{proof}
Let $B$ be a $k$-constellation of an order preserving canonical function $f\colon\delcn\To\Delta$, and let $B'$ be an arbitrary $(k-1)$-subconstellation. Then $B'$ is $\OP$-fittable by assumption; let $C\in\OP$ be so that $B'$ fits $C$. Let $S$ be the union of the orbits concerned by $B$, and $S'$ the union of those concerned by $B'$. By Lemma~\ref{lem:invertonconstellation}, there exists a canonical $g\in\overline{\{\gamma\circ f\circ \beta \;|\; \gamma\in C,\; \beta\in\Aut(\Delta,c_1,\ldots,c_n)\}}$ which behaves like the identity function on $S'$. The constellation given by $S$ for $g$ is $\OP$-fittable by assumption, and so $B$ is $\OP$-fittable by Lemma~\ref{lem:reduction}.
\end{proof}

%\begin{lem}\label{lem:identityonsub}
%Assume that all $(k-1)$-constellations of order preserving canonical functions are $\OP$-fittable. 
%Let $S$ be the union of $k$ infinite $1$-types, and let $S'$ be a subset of $S$ that is the union of $k-1$ infinite $1$-types. 
%Assume that all $k$-constellations of order preserving canonical functions whose underlying set is $S$ and whose restriction to $S'$ is the behavior of the identity function are $\OP$-fittable. 
%Then all $k$-constellations of order preserving canonical functions whose underlying set is $S$ are $\OP$-fittable.
%\end{lem}
%\begin{proof}
%Let $B$ be a $k$-constellation with underlying set $S$ of an order preserving canonical function $f\colon\delcn\To\Delta$, and let $B'$ be the restriction of $B$ to $S'$. 
%Then $B'$ is $\OP$-fittable by assumption; let $C\in\OP$ be so that $B'$ fits $C$. 
%By Lemma~\ref{lem:invertonconstellation}, there exists a canonical $g\in\overline{\{\gamma\circ f\circ \beta \;|\; \gamma\in C,\; \beta\in\Aut(\Delta,c_1,\ldots,c_n)\}}$ which behaves like the identity function on $S'$. 
%The constellation given by $S$ for $g$ is $\OP$-fittable by assumption, and so $B$ is $\OP$-fittable by Lemma~\ref{lem:reduction}.
%\end{proof}

\subsubsection{Fitting the 1-constellations}
Clearly, all 1-constellations fit either $\Aut(D;<)$ (when eradicating edges or non-edges; cf.~Lemma~\ref{lem:eradicatesonorbit}), or $\cl{\mix -\id}$ (when flipping the graph relation), or $\adel$ (when keeping the graph relation).

\subsubsection{Fitting the 2-constellations}\label{subsec2C}
\begin{figure}
\includegraphics[width= 10cm]{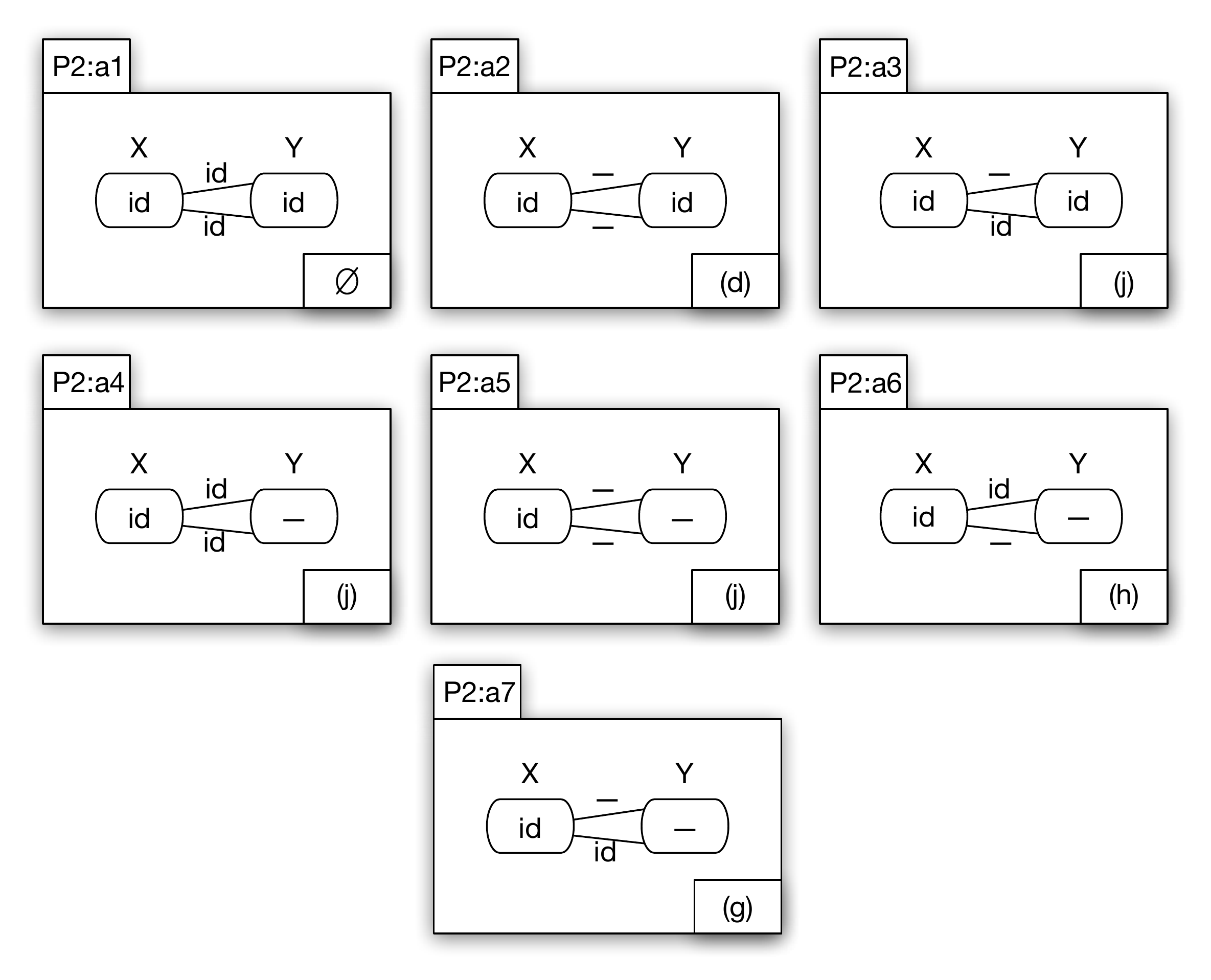}
%\caption{}{}
\caption{Cases C2:a}\label{C2a}
\end{figure}

We first treat the 16 different 2-constellations with the two orbits on the same level.
Figure~\ref{C2a} lists all those constellations up to symmetry with the behavior $id$ on $X$ (see Lemma~\ref{lem:identityonsub} with $k=2$ and $S'=X$). As we have mentioned, we only verify that each of these 2-constellations $B$ forces the group $C$ indicated in the lower right of the respective picture. 
In each case, $f$ denotes a canonical function from $(\Delta,c_1,\dots,c_n)$ to $\Delta$ satisfying $B$, and we have to show that $f$ generates all functions in $C$. 
On these pictures, $X$ and $Y$ are orbits, and the two lines connecting them indicate two different kinds of pairs $(x,y)\in X\times Y$: those with $x<y$ indicated by the ascending line and those with $x<y$ by the descending one. 
The signs $\id$ and $-$ written in the orbits or on the lines specify the behavior of $f$, where $-$ corresponds to the behavior of $\mix - \id$. 

\vspace{.2cm}

C2:a1. There is nothing to prove, every function ``generates'' $\Aut(\Delta)$ by definition. 

C2:a2. Let $A$ be an arbitrary finite subset of $D$ and let $S\subseteq A$. There exists a $\gamma\in \Aut(\Delta)$ such that $\gamma[S]\subseteq X$ and $\gamma[A\setminus S]\subseteq Y$. Then $f\circ \gamma$ flips the graph relation between $S$ and $A\setminus S$, and keeps it otherwise on $A$. Thus, $f$ generates $\mix \sw \id$.
%
%C2:a3. As $f$ behaves like $\mix - \id$ on $X$ it is clear that $f$ generates $\cl{\mix - \id}$. By composing $f$ with $\mix - \id$ we obtain a constellation as in C2:a2. Hence, $f$ generates $\mix \sw \id$.
%
%C2:a4. As $f$ behaves like $\mix - \id$ on $X$, it generates $\cl{\mix - \id}$.

C2:a3, C2:a4, C2:a5. Let $A$ be an arbitrary finite subset of $D$ that consists of the elements $a_1<\cdots <a_k$. Let $1\leq i<j\leq k$. There exist $\gamma_1,\gamma_2,\gamma_3,\gamma_4 \in \Aut(\Delta)$ such that 
\begin{itemize}
\item $\gamma_1[A] \subseteq X$, 
\item $\gamma_2[A \setminus \{a_i\}] \subseteq X$, $\gamma_2(a_i) \in Y$,
\item $\gamma_3[A \setminus \{a_j\}] \subseteq X$, $\gamma_3(a_j) \in Y$,
\item $\gamma_4[A \setminus \{a_i,a_j\}] \subseteq X$, $\gamma_4[\{a_i,a_j\}] \subseteq Y$.
\end{itemize}
We would like to `combine' the four behaviors $f \circ \gamma_i$, for $i \in \{1,2,3,4\}$,
in order to flip the graph relation on the pair $(a_i,a_j)$ and to keep it otherwise on $A$.
To formalize this, we choose 
$\delta_1,\delta_2,\delta_3,\delta_4 \in \Aut(\Delta)$ such that 
\begin{itemize}
\item $\delta_1[A] \subseteq X$, 
\item $\delta_2\circ f \circ \delta_1[A \setminus \{a_i\}] \subseteq X$, $\delta_2 \circ f \circ \delta_1(a_i) \in Y$,
\item $\delta_3\circ f \circ \delta_2 \circ f \circ \delta_1[A \setminus \{a_j\}] \subseteq X$, $\delta_3\circ f \circ \delta_2 \circ f \circ \delta_1(a_j) \in Y$,
\item $\delta_4\circ f \circ \delta_3 \circ f \circ \delta_2 \circ f \circ \delta_1[A \setminus \{a_i,a_j\}] \subseteq X$, $\delta_4\circ f \circ \delta_3 \circ f \circ \delta_2 \circ f \circ \delta_1[\{a_i,a_j\}] \subseteq Y$.
\end{itemize}
Consider now $f \circ \delta_4 \circ f \circ \delta_3 \circ \cdots \circ \delta_1$; this function indeed flips the graph relation on the pair $(a_i,a_j)$
and otherwise behaves as the identity on $A$. One can now use Lemma~\ref{lem:eradicatesonfinite} to show that $f$ generates $\Aut(D;<)$. 

C2:a6. Let $A$ be an arbitrary finite subset of $D$ and let $S\subseteq A$ be an upward closed subset. There exists a $\gamma\in \Aut(\Delta)$ such that $\gamma[S]\subseteq Y$ and $\gamma[A\setminus S]\subseteq X$. Then $f\circ \gamma$ flips the graph relation on $S$ and keeps it otherwise on $A$. Thus, $f$ generates $\apu$.

C2:a7. This can be shown analogously to the previous case.
%Let $A$ be an arbitrary finite subset of $D$ that consists of the elements $a_1<\cdots <a_k$. Let $S\subseteq A$ be an upward closed subset. There exists a $\gamma\in \Aut(\Delta)$ such that $\gamma[S]\subseteq Y$ and $\gamma[A\setminus S]\subseteq X$. Then $f\circ \gamma$ flips the graph relation on $S$ and keeps it otherwise on $A$. Thus $f$ generates $\apu$.

%The combination of the four modifications that the application of $f$ yields flips the graph relation on the pair $(a_i,a_j)$ and keeps it otherwise on $A$. Hence, $A$ can be mapped to an empty graph, and thus $f$ generates $\Aut(D;<)$.

\vspace{.2cm}

\begin{figure}
\includegraphics[width= 10cm]{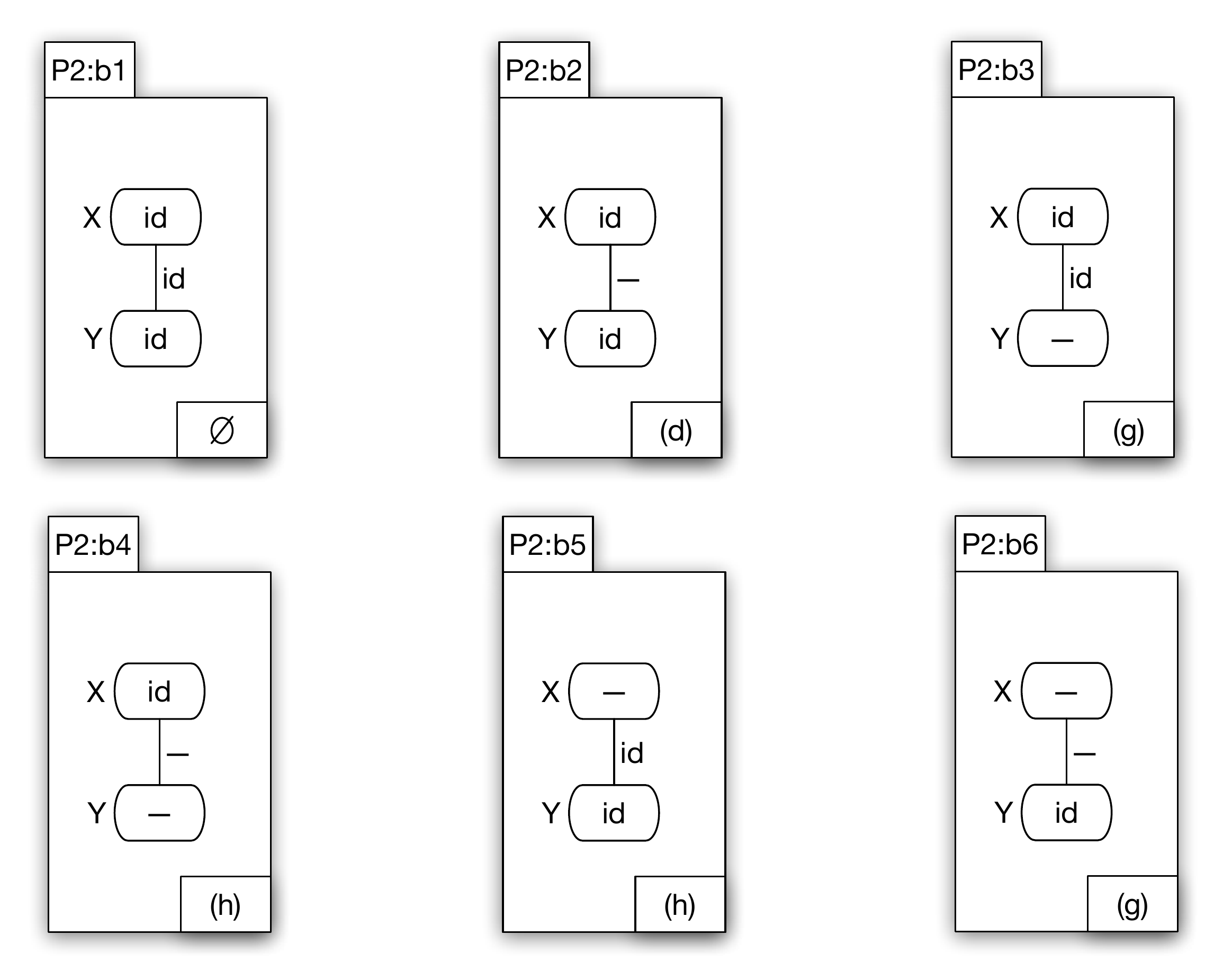}
%\caption{}{}
\caption{Cases C2:b}\label{2Cb}
\end{figure}

We now discuss the 2-constellations with the two orbits on different levels.
There are eight such constellations. 
We do not consider those which have the behavior of $\mix - \id$ in both orbits (see Lemma~\ref{lem:identityonsub}), the rest of the cases are illustrated in Figure~\ref{2Cb}.

\vspace{.2cm}

C2:b1. There is nothing to prove, every function generates $\Aut(\Delta)$.

C2:b2. Let $A$ be an arbitrary finite subset of $D$ that consists of the elements $a_1<\cdots <a_k$. Let $1\leq i\leq k$. Let $A_1=\{a_1, \ldots, a_{i-1}\}$, 
$A_2=\{a_{i+1}, \ldots, a_{k}\}$. There exists a $\gamma \in \Aut(\Delta)$ such that 
$\gamma[A_1] \subseteq X$, $\gamma(a_i)\in X$, and $\gamma[A_2] \subseteq Y$.
We can also choose $\delta \in \Aut(\Delta)$ 
such that $\delta(a_i) \in Y$, $\delta[A_1] \subseteq X$, and $\delta[A_2] \subseteq Y$. 
Then combining the effect of $f \circ \delta$ and of $f \circ \gamma$ (formalized as
in case C2:a3), we obtain a function that
flips the graph relation between $\{a_i\}$ and $A\setminus\{a_i\}$ and keeps it otherwise on $A$. Hence, $f$ generates $\mix \sw \id$.

C2:b3. Let $A$ be an arbitrary finite subset of $D$ that consists of the elements $a_1<\cdots <a_k$. Let $S\subseteq A$ be a downward closed subset. There exists a $\gamma\in \Aut(\Delta)$ such that $\gamma[S]\subseteq Y$ and $\gamma[A\setminus S]\subseteq X$. Then $f\circ \gamma$ flips the graph relation on $S$ and keeps it otherwise on $A$. Thus, $f$ generates $\apl$.

C2:b5. Analogously to the previous case it can be shown that $f$ generates
$\apu$. 
%Let $A$ be an arbitrary finite subset of $D$ that consists of the elements $a_1<\cdots <a_k$. Let $S\subseteq A$ be an upward closed subset. There exists a $\gamma\in \Aut(\Delta)$ such that $\gamma[S]\subseteq Y$ and $\gamma[A\setminus S]\subseteq X$. Then $f\circ \gamma$ flips the graph relation on $S$ and keeps it otherwise on $A$. Thus $f$ generates $\apu$.

C2:b4. As $f$ behaves like $\mix - \id$ on $Y$ we have that $f$ generates $\mix - \id$. By composing $\mix - \id$ with $f$ we obtain a constellation as in C2:b5, hence $f$ generates $\apu$.

C2:b6. As $f$ behaves like $\mix - \id$ on $X$ we have that $f$ generates $\mix - \id$. By composing $\mix - \id$ with $f$ we obtain a constellation as in C2:b3, hence $f$ generates $\apl$.

%C2:b7. As $f$ behaves like $\mix - \id$ on $Y$ we have that $f$ generates $\mix - \id$. By composing $\mix - \id$ with $f$ we obtain a constellation as in C2:b2, hence $f$ generates $\mix \sw \id$.
%
%C2:b8. As $f$ behaves like $\mix - \id$ on $Y$ we have that $f$ generates $\mix - \id$.

\subsubsection{Fitting the 3-constellations}\label{subsec3C}
\begin{figure}
\includegraphics[width= 10cm]{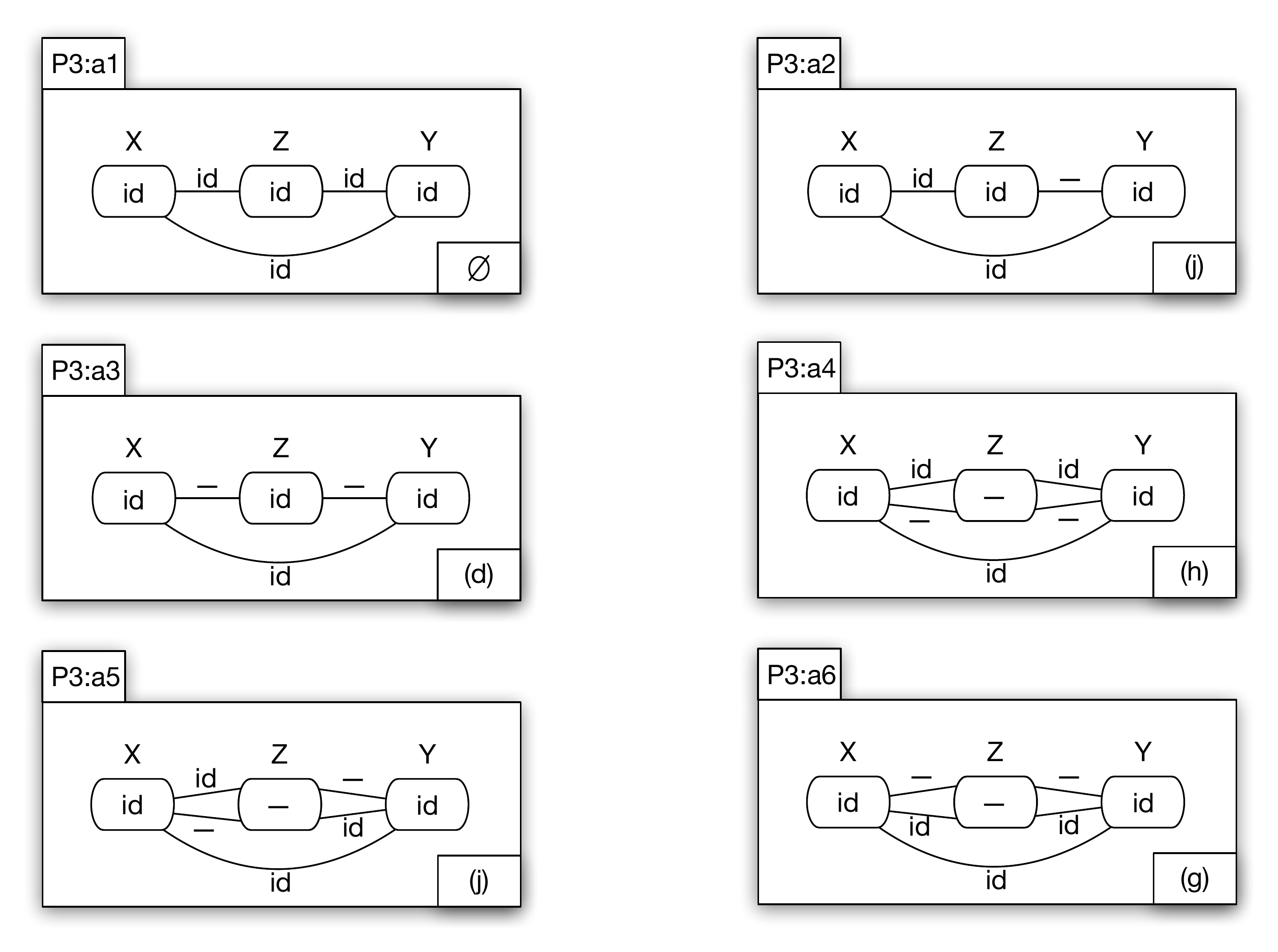}
%\caption{}{}
\caption{Cases C3:a}\label{C3a}
\end{figure}

We check the 3-constellations according to the following case distinction: either all three orbits are on the same level, or two of them are on the same level, or all three orbits are on different levels. In each of these cases, we may refer to Lemma~\ref{lem:identityonsub} and assume that for two orbits $X, Y$ of our choice we have the identity behavior on and between $X$ and $Y$.
There are eight 3-constellations with all three orbits $X,Y,Z$ on the same level with the identity behavior on and between $X$ and $Y$, and such that the constellation does not contain a full 2-subconstellation. We only need to consider these 
3-constellations up to symmetry of $X$ and $Y$, which leads to six cases, drawn in Figure~\ref{C3a}.

\vspace{.2cm}

C3:a1. There is nothing to prove, every function generates $\Aut(\Delta)$.

C3:a2. Let $A$ be an arbitrary finite subset of $D$ that consists of the elements $a_1<\cdots <a_k$. Let $1\leq i<j\leq k$. There exists a $\gamma\in \Aut(\Delta)$ such that $\gamma(a_i)\in Z$, $\gamma(a_j)\in Y$ and $\gamma[A\setminus \{a_i,a_j\}]\subseteq X$. Then $f\circ \gamma$ flips the graph relation between $a_i$ and $a_j$, and keeps it otherwise on $A$. Thus, $f$ generates $\Aut(D;<)$.

C3:a3. The constellation induced by $X\cup Z$, treated in C2:a2, shows that $f$ generates $\mix \sw \id$.

C3:a4. The constellation induced by $X\cup Z$, treated in C2:a6, shows that $f$ generates $\apu$.  

C3:a5. Let $A$ be an arbitrary finite subset of $D$ that consists of the elements $a_1<\cdots <a_k$. Let $1\leq i<j\leq k$. 
Let $A_1=\{a_1, \ldots, a_{i-1}\}$, $A_2=\{a_{i+1}, \ldots, a_{j-1}\}$, and $A_3=\{a_{j+1}, \ldots, a_{k}\}$. 
There exist $\gamma_1, \gamma_2, \gamma_3\in \Aut(\Delta)$ such that 
\begin{itemize}
\item $\gamma_1[A_1]\subseteq X$, $\gamma_1[A_2]\subseteq Z$, $\gamma_1[A_3]\subseteq Y$, 
\item $\gamma_2[A_1]\subseteq Z$, $\gamma_2[A_2]\subseteq Y$, 
$\gamma_2[A_3]\subseteq X$, 
\item $\gamma_3[A_1]\subseteq Y$, $\gamma_3[A_2]\subseteq X$, 
$\gamma_3[A_3]\subseteq Z$, and 
\item $\gamma_m(a_i), \gamma_m(a_j) \in Z$ for $m \in \{1,2,3\}$. 
\end{itemize}
Then the combined effec to $f \circ \gamma_i$ for $m \in \{1,2,3\}$ (for a formalization of this, see case C2:a3) flips the graph relation between $a_i$ and $a_j$, and keeps it otherwise on $A$. Thus, $f$ generates $\Aut(D;<)$.

C3:a6. The constellation induced by $X\cup Z$, treated in C2:a7, shows that $f$ generates $\apl$.

\vspace{.2cm}

\begin{figure}
\includegraphics[width= 10cm]{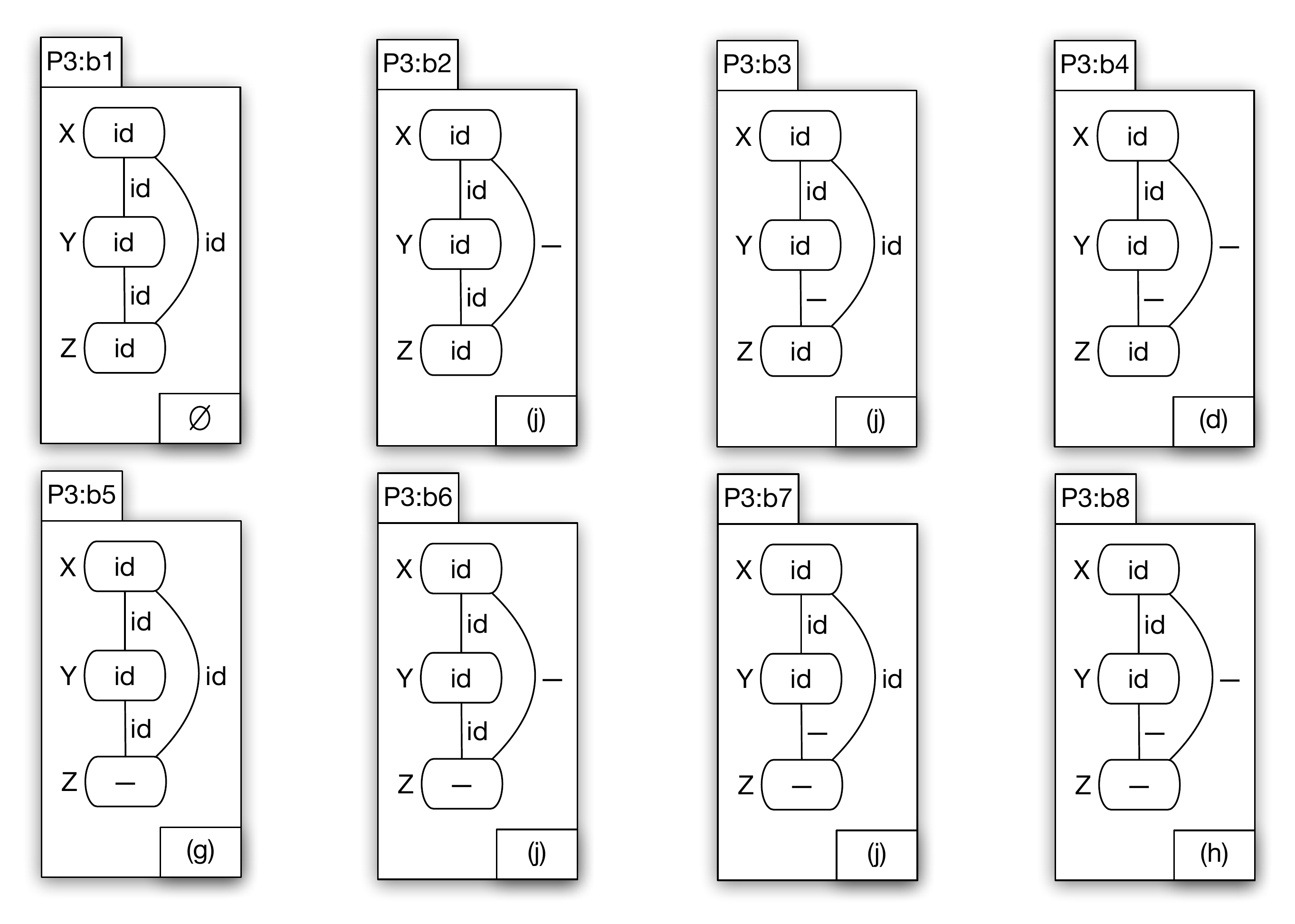}
%\caption{}{}
\caption{Cases C3:b}\label{C3b}
\end{figure}

Figure~\ref{C3b} contains the 3-constellations with orbits $X>Y>Z$, such that the behavior is identical on and between $X$ and $Y$, and the constellation does not have a full 2-subconstellation. 

%For the remaining arguments, we make a convention that allows a more concise presentation. Please revisit the cases C2:a7, C2:a8, C2:a9, C2:a10, 
%or the case C3:a5 above. In those cases we have described how to use the constellation of $f$ to obtain several functions $f_1,f_2,\dots$ with a certain behavior on an arbitrary finite subset of $A$. 
%We have then chained those functions to obtain yet another behavior on $A$.
%In order to do so, we had to take certain care to choose automorphisms to shift the image of $A$ under $f$

\vspace{.2cm}

C3:b1. There is nothing to prove, every function generates $\Aut(\Delta)$.

C3:b2. Let $A$ be an arbitrary finite subset of $D$ that consists of the elements $a_1<\cdots <a_k$. Let $1\leq i<j\leq k$. Let $A_1=\{a_1, \ldots, a_{i-1}\}$, $A_2=\{a_{i+1}, \ldots, a_{j-1}\}$, $A_3=\{a_{j+1}, \ldots, a_{k}\}$. There exist $\gamma_1, \gamma_2, \gamma_3, \gamma_4\in \Aut(\Delta)$ such that
\begin{itemize}
\item $\gamma_m(A_1)\subseteq Z$, $\gamma_m(A_2)\subseteq Y$, $\gamma_m(A_3)\subseteq X$ for $m=1,2,3,4$ 
\item $\gamma_1(a_i) \in Z$, $\gamma_1(a_j) \in Y$,
\item $\gamma_1(a_i) \in Y$, $\gamma_1(a_j) \in Y$,
\item $\gamma_1(a_i) \in Z$, $\gamma_1(a_j) \in X$,
\item $\gamma_1(a_i) \in Y$, $\gamma_1(a_j) \in X$.
%\item $\gamma_m(a_i)$ is in $X$ or in $Y$, and 
%\item $\gamma_m(a_j)$ is in $Y$ or $Z$ such that all four possibilities occur once. 
\end{itemize}
By combining the effect of the functions $f \circ \gamma_m$, for $m \in \{1,2,3,4\}$ (as formalized in C2:a3) we obtain a function that flips the graph relation between $a_i$ and $a_j$, and keeps it otherwise on $A$. Thus, $f$ generates $\Aut(D;<)$.

C3:b3. The constellation induced by $Y\cup Z$ shows that $f$ generates $\mix \sw \id$. 
By applying an order preserving permutation that flips the graph relation between $Z$ and its complement, we arrive at C3:b2. Hence, $f$ generates $\Aut(D;<)$.

C3:b4. The constellation induced by $Y\cup Z$, treated in C2:b2, shows that $f$ generates $\mix \sw \id$.

C3:b5. The constellation induced by $Y\cup Z$, treated in C2:b3, shows that $f$ generates $\apl$.

C3:b6.  Let $A$ be an arbitrary finite subset of $D$ that consists of the elements $a_1<\cdots <a_k$. Let $1\leq i<j\leq k$. Let $A_1=\{a_1, \ldots, a_{i-1}\}$, $A_2=\{a_{i+1}, \ldots, a_{j-1}\}$, $A_3=\{a_{j+1}, \ldots, a_{k}\}$. There exist $\gamma_1, \gamma_2, \gamma_3, \gamma_4\in \Aut(\Delta)$ such that $\gamma_m[A_1]\subseteq Z$, $\gamma_m[A_2]\subseteq Y$, $\gamma_m[A_3]\subseteq X$ for $m=1,2,3,4$, and $\gamma_m(a_i)$ is in $X$ or in $Y$ and $\gamma_m(a_j)$ is in $Y$ or $Z$ such that all four possibilities occur once. Then the combined effect of $f\circ \gamma_m$ for $m \in \{1,2,3,4\}$ (as formalized in C2:a3) flips the graph relation between $a_i$ and $a_j$, and keeps it otherwise on $A$. Thus, $f$ generates $\Aut(D;<)$.

C3:b7. There exists some $\gamma\in \Aut(\Delta)$ such that $\gamma[f[X]\cup f[Y]]\subseteq Y$ and $\gamma[f[Z]]\subseteq Z$. Then $f\circ \gamma\circ f$ has the same constellation as in C3:b2. Hence, $f$ generates $\Aut(D;<)$.

C3:b8. The constellation induced by $Y\cup Z$, treated in C2:b4, shows that $f$ generates $\apu$.

\vspace{.2cm}

\begin{figure}
\includegraphics[width= 10cm]{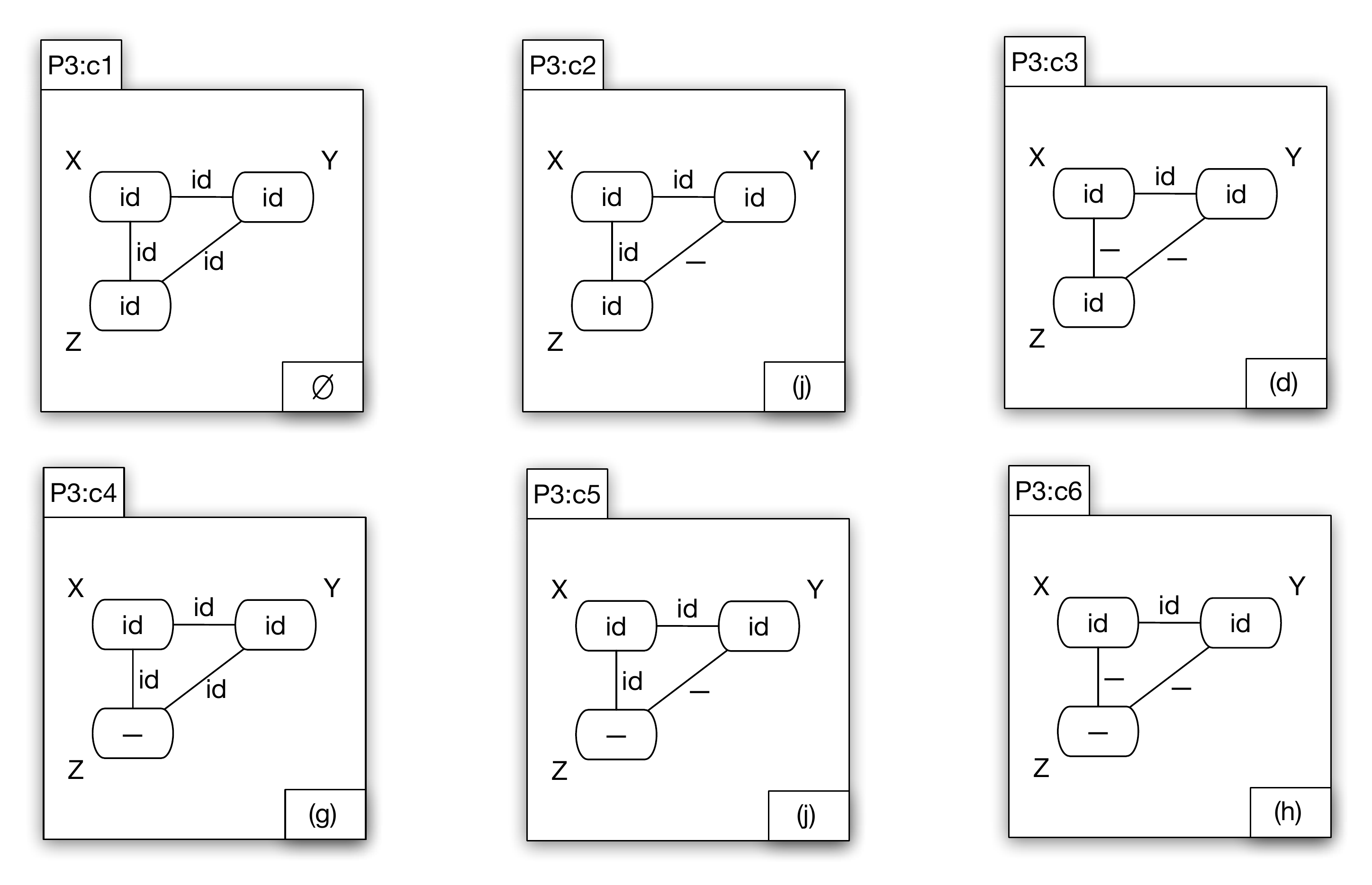}
%\caption{}{}
\caption{Cases C3:c}\label{C3c}
\end{figure}

Figure~\ref{C3c} shows the 3-constellations with orbits $X,Y,Z$ such that $X$ and $Y$ are on the same level, $Z<X$, the behavior is identical on and between $X$ and $Y$, and the constellation does not contain a full 2-subconstellation.
Note that the cases C3:c have the same group in the right lower corner
than the corresponding cases in C3:b. And indeed, it is easy to see that each case in C3:c forces at least the groups that are forced by the corresponding constellation in C3:b.

\begin{figure}
\includegraphics[width= 10cm]{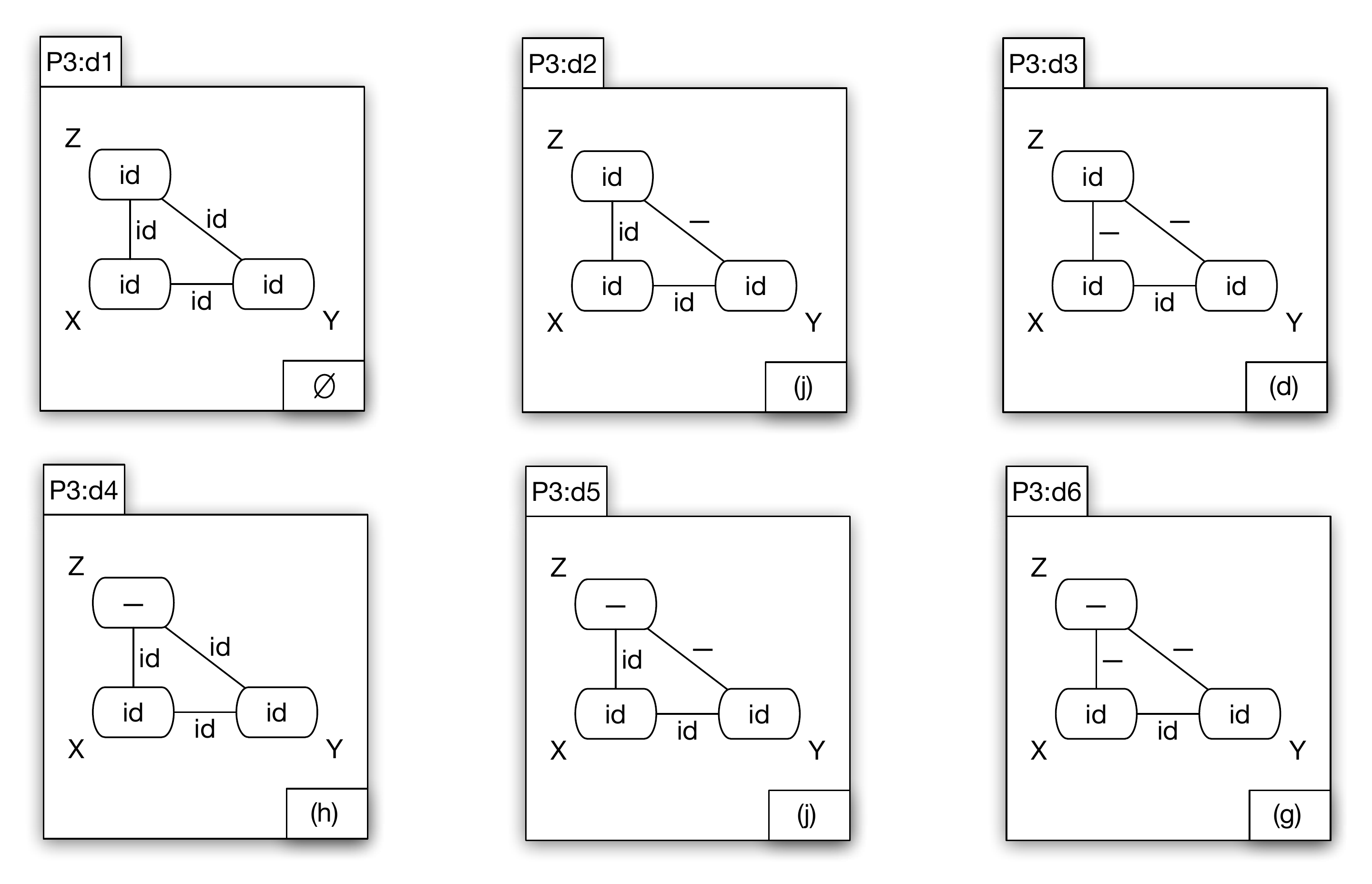}
%\caption{}{}
\caption{Cases C3:d}\label{C3d}
\end{figure}

Figure~\ref{C3d} contains the 3-constellations with orbits $X,Y,Z$ such that $X$ and $Y$ are on the same level, $Z>X$, the behavior is identical on and between $X$ and $Y$, and the constellation does not contain a full 2-subconstellation.
%The proof is analogous to the proof for Figure~\ref{C3c}.
The situation here is analogous to the Cases~C3:c. 

\subsubsection{Fitting the 4-constellations}\label{subsec4C}

Observe that so far, every $k$-constellation which is identical on all orbits has been either full or compatible with $R^{(3)}$. This leads to the following simplification which reduces the number of $4$-constellations to consider.

\begin{lem}\label{constidonorbits}
 Let $k\geq 1$ and let $B$ be a $k$-constellation of an order preserving canonical function $f\colon \delcn\To\Delta$. If the restriction of $B$ to any orbit is the identical behavior, then $B$ is $\OP$-fittable.
\end{lem}
\begin{proof}
 According to Lemma~\ref{full} we may assume that $B$ has no subconstellation that forces $\Aut(D;<)$. By Lemma~\ref{lem:eradicatebetweenorbits} we may moreover assume that $B$ either keeps or flips the graph relation between any two orbits, or otherwise it is full. 
 If $B$ is identical between any pair of orbits, then $B$ fits $\Aut(\Delta)$.  Hence, we may assume that there exist two orbits such that the behavior of $B$ flips the graph relation between them. According to the Cases C2:a2 and C2:b2 we have that $B$ forces $\cl{\mix \sw \id}$. Assume that $B$ violates $R^{(3)}$. As $R^{(3)}$ is a ternary relation, there exist three orbits $X,Y,Z$ such that the subconstellation $B'$ induced by  $X\cup Y\cup Z$ violates $R^{(3)}$. So $B'$ has the identity behavior on every orbit, and violates $R^{(3)}$; according to the case-by-case analysis of the 3-constellations, we have that all such $B'$ are full. Thus, $B$ is compatible with $R^{(3)}$, and so $B$ fits $\Aut(D;R^{(3)})=\cl{\mix \sw \id}$.
\end{proof}

	We check the 4-constellations according to a similar case distinction as before. In all the cases we may refer to Lemma~\ref{lem:identityonsub} and 
choose three orbits for which we assume that the behavior is identical on and between them. According to Lemma~\ref{constidonorbits} we may assume that the behavior flips the graph relation on the remaining fourth orbit. We only analyze those constellations which do not contain a subconstellation forcing $\Aut(D;<)$ -- see Lemma~\ref{full}.

\begin{figure}
\includegraphics[width= 10cm]{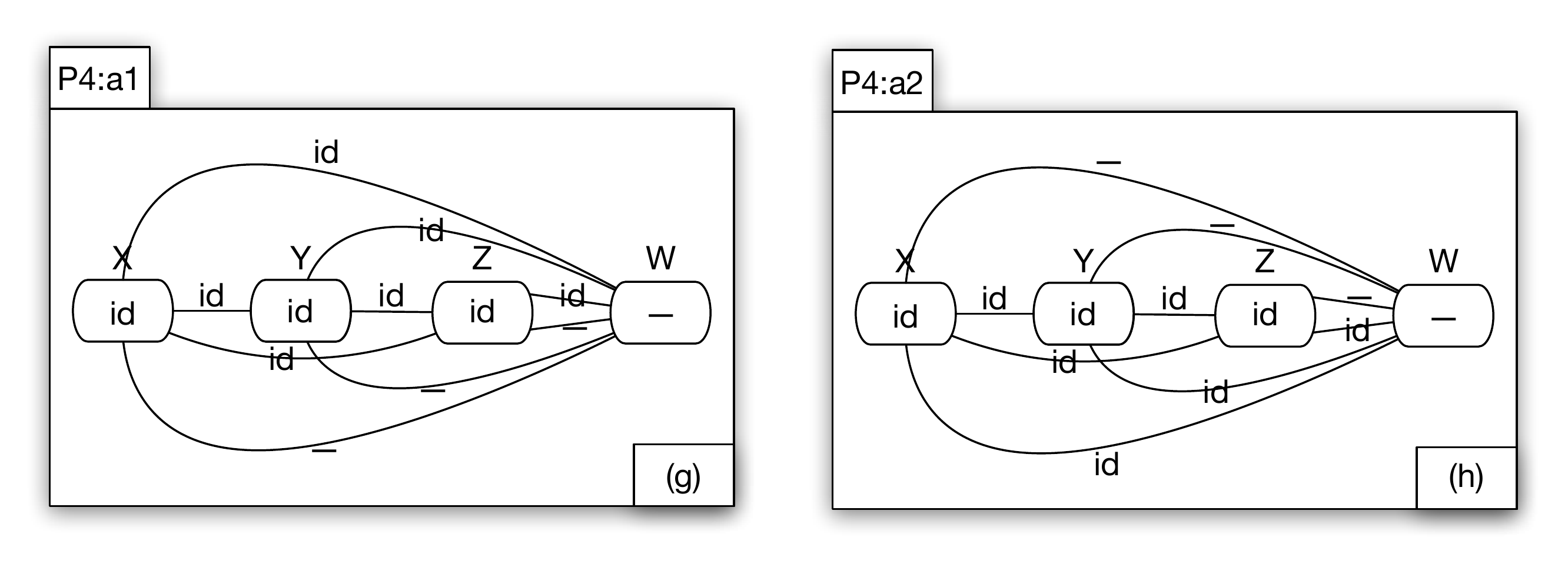}
%\caption{}{}
\caption{Cases C4:a}\label{C4a}
\end{figure}

Figure~\ref{C4a} contains the remaining 4-constellations with all orbits on the same level. There are only two constellations to be checked as we assume that the behavior is identical on and between $X, Y, Z$, flips the graph relation on $W$, and that there is no full subconstellation. In particular, any 3-subconstellation belongs to one of the Cases C3:a1, C3:a4, C3:a6.

\vspace{.2cm}

C4:a1.  The constellation induced by $Z \cup W$, treated in C2:a7, shows that $f$ generates $\apl$.

C4:a2.  The constellation induced by $Z \cup W$, treated in C2:a6, shows that $f$ generates $\apu$.

\vspace{.2cm}

\begin{figure}
\includegraphics[width= 10cm]{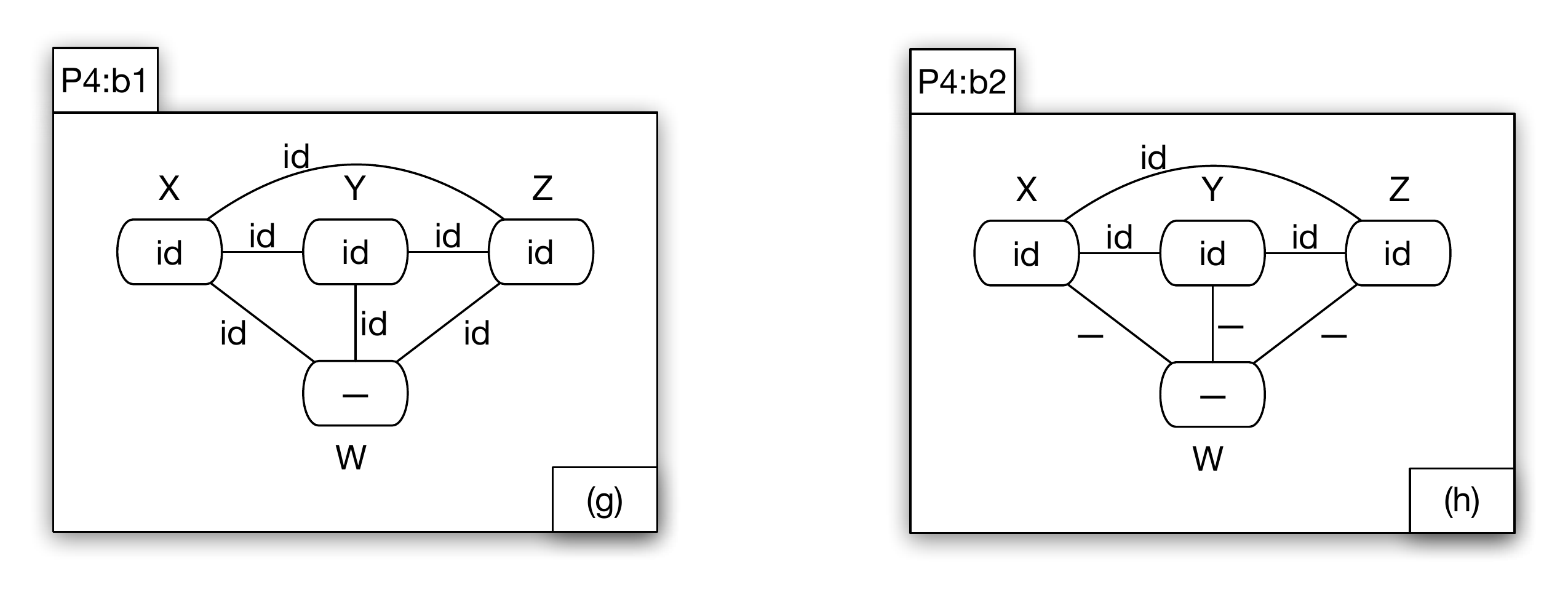}
%\caption{}{}
\caption{Cases C4:b}\label{C4b}
\end{figure}

Figure~\ref{C4b} shows the relevant 4-constellations with three orbits $X,Y,Z$ on the same level and the fourth orbit $W<X$. There are only two constellations to be checked as we assume that the behavior is identical on and between $X, Y, Z$, flips the graph relation on $W$, and that there is no full subconstellation. 
In particular, either the graph relation is kept everywhere else, or the graph relation is flipped everywhere else. 

\vspace{.2cm}

C4:b1. The constellation induced by $X\cup W$, treated in C2:b3, shows that $f$ generates $\apl$.

C4:b2. The constellation induced by $X\cup W$, treated in C2:b4, shows that $f$ generates $\apu$.

\vspace{.2cm}

\begin{figure}
\includegraphics[width= 8cm]{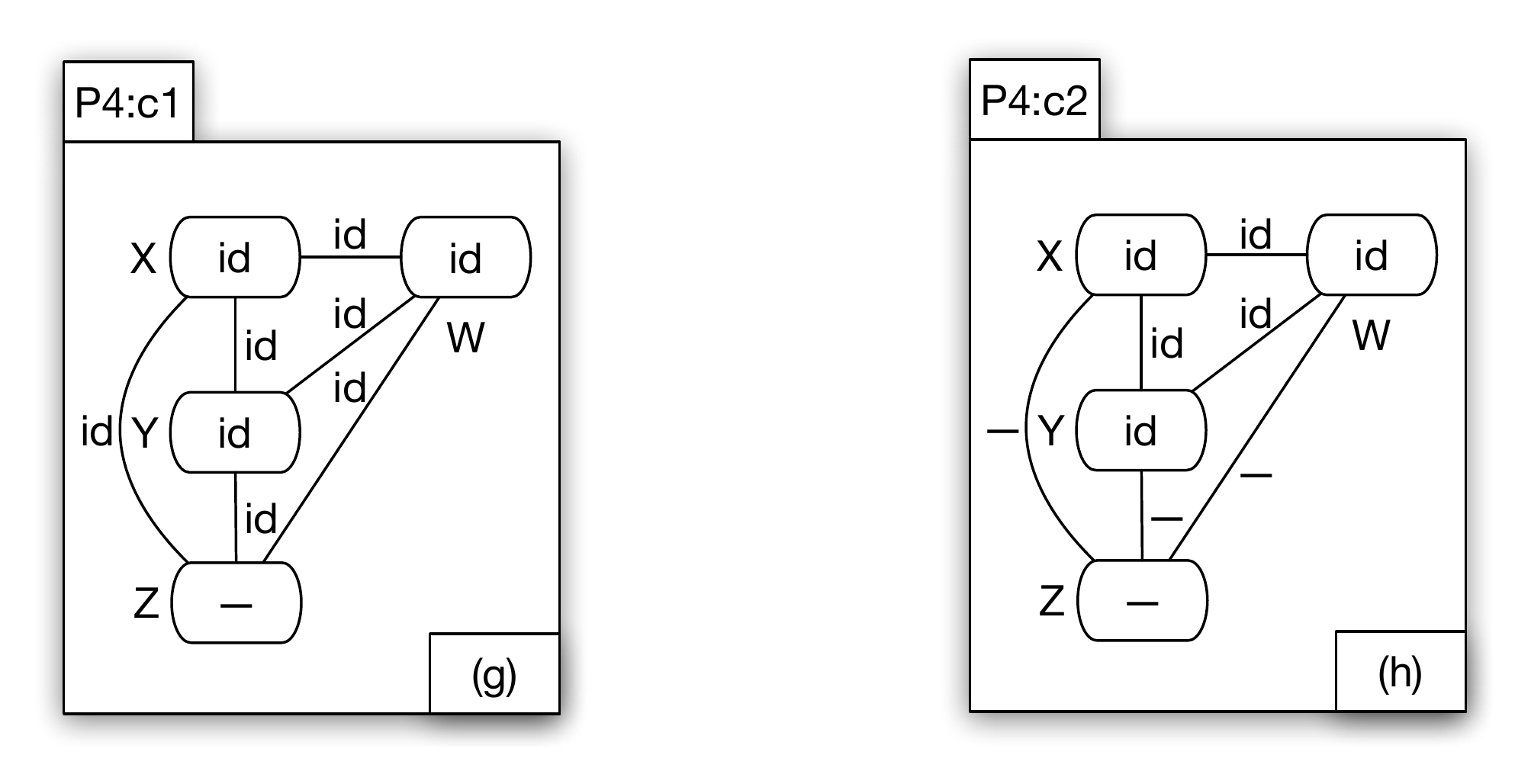}
%\caption{}{}
\caption{Cases C4:c}\label{C4c}
\end{figure}

The case with $X,Y,Z$ on the same level and the fourth orbit $W>X$ is analogous to the previous case, Case~C4b, and we omit it. 
Figure~\ref{C4c} contains the relevant 4-constellations with three orbits $X>Y>Z$ and the fourth orbit $W$ on the same level as $X$. There are only two constellations to be checked as we assume that the behavior is identical on and between $X, Y, W$, flips the graph relation on $Z$, and that there is no full subconstellation. 
In particular, either the graph relation is kept everywhere else, or the graph relation is flipped everywhere else. 

\vspace{.2cm}

C4:c1.  The constellation induced by $Z \cup W$, treated in C2:b3, shows that $f$ generates $\apl$.

C4:c2.  The constellation induced by $Z\cup W$, treated in C2:b4, shows that $f$ generates $\apu$.

\vspace{.2cm}

\begin{figure}
\includegraphics[width= 10cm]{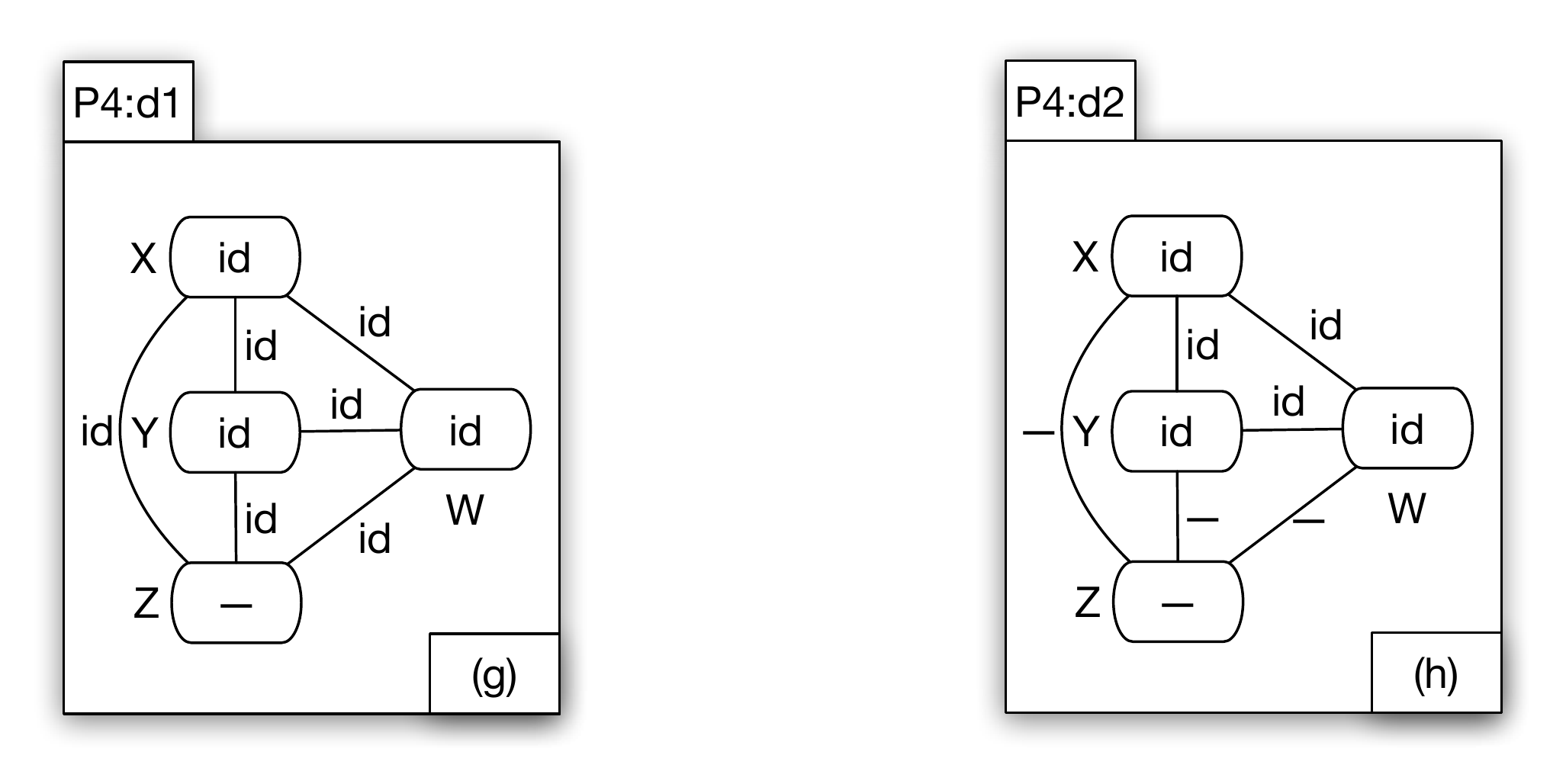}
%\caption{}{}
\caption{Cases C4:d}\label{C4d}
\end{figure}

The case with $X>Y>Z$ and the fourth orbit $W$ on the same level as $Z$ is analogous to the previous case, Case~C4:c, so we omit it. 
Figure~\ref{C4d} contains the relevant 4-constellations with three orbits $X>Y>Z$ and the fourth orbit $W$ on the same level as $Y$. There are only two constellations to be checked as we assume that the behavior is identical on and between $X, Y, W$, flips the graph relation on $Z$, and that there is no full subconstellation. In particular, either the graph relation is kept everywhere else, or the graph relation is flipped everywhere else.

\vspace{.2cm}

C4:d1. %The group $\cl{\apl}$ is forced because of the subconstellation induced by $X\cup W$.
 The constellation induced by $Z \cup W$, treated in C2:b3, shows that $f$ generates $\apl$.

C4:d2. %The group $\cl{\apu}$ is forced because of the subconstellation induced by $X\cup W$.
 The constellation induced by $Z \cup W$, treated in C2:b4, shows that $f$ generates $\apu$.

\vspace{.2cm}

\begin{figure}
\includegraphics[width=7cm]{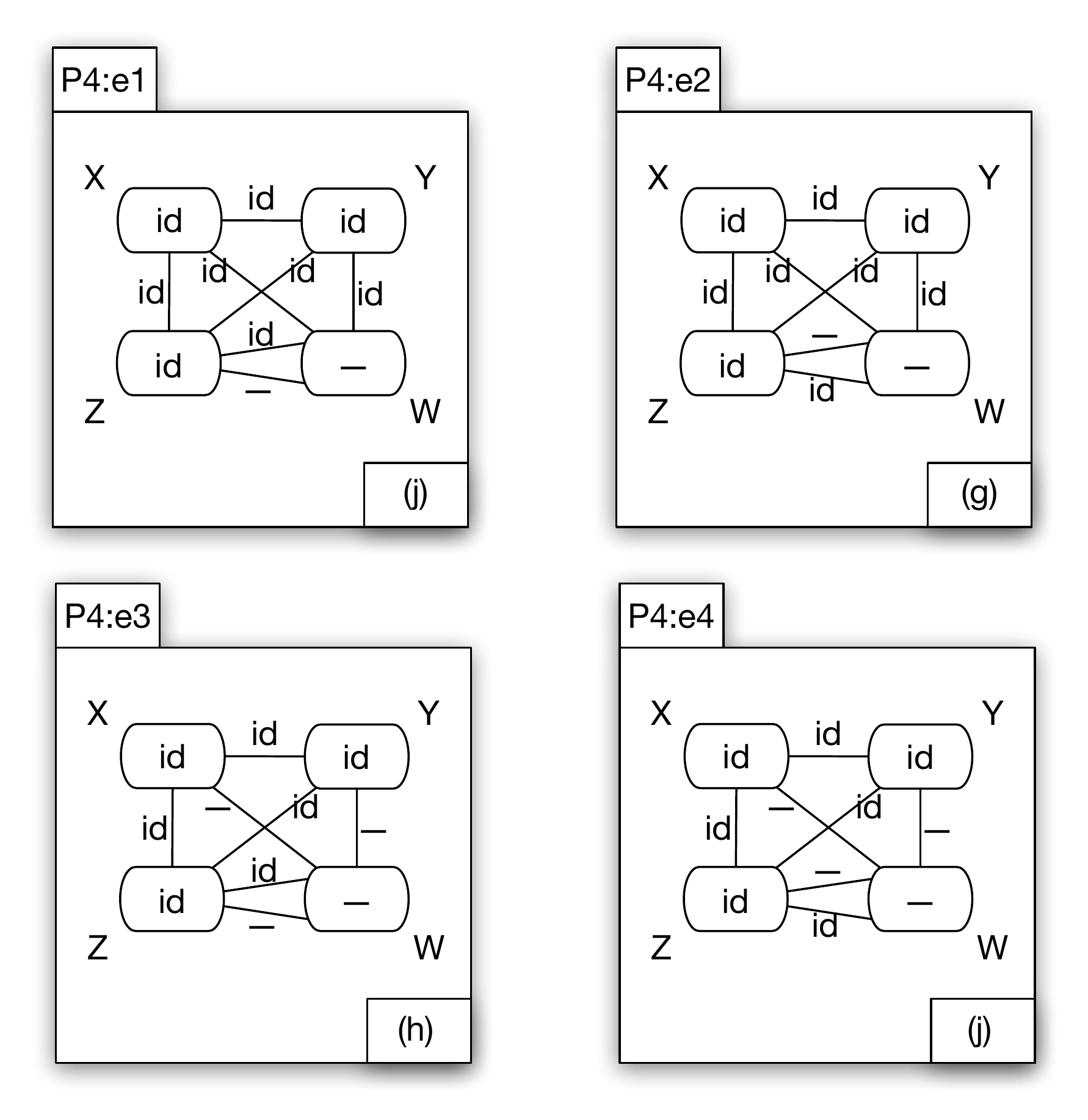}
%\caption{}{}
\caption{Cases C4:e}\label{C4e}
\end{figure}

Figure~\ref{C4e} contains the relevant 4-constellations with orbits $X$ and $Y$ on the same level, $Z$ and $W$ on the same level, and $X>Z$. There are only four constellations to be checked, because we assume that the behavior is identical on and between $X, Y, Z$, flips the graph relation on $W$, and that there is no full subconstellation. In particular, either the graph relation is kept between $X$ and $W$ and between $Y$ and $W$, or the graph relation is flipped between $X$ and $W$ and between $Y$ and $W$. There are also two possibilities for the behavior between $Z$ and $W$ according to Subsection~\ref{subsec2C}, Case C2.

\vspace{.2cm}

C4:e1. Let $B'$ be the constellation induced by $Y\cup Z\cup W$. Then $B'$ is $\OP$-fittable by Subsection~\ref{subsec3C}. 
Consider $a \in W$, $b,c \in Z$, $d \in Y$, such that $a<b<c<d$ and no edges between vertices from $\{a,b,c,d\}$. Then $(a,b,c,d) \in S_4$, 
and $B'$ violates $S_4$ on $(a,b,c,d)$.
Thus, $B'$ is full.

C4:e2.  The constellation induced by $Z \cup W$, treated in C2:a7, shows that $f$ generates $\apl$.

C4:e3.  The constellation induced by $Z \cup W$, treated in C2:a6, shows that $f$ generates $\apu$.

C4:e4. Let $B'$ be the subconstellation induced by $Y\cup Z\cup W$. Then $B'$ is $\OP$-fittable by Subsection~\ref{subsec3C}. 
Consider $a \in W$, $b,c \in Z$, $d \in Y$, such that $a<b<c<d$ and no edges between vertices from $\{a,b,c,d\}$. Then $(a,b,c,d) \in S_4$, and $B'$ violates $S_4$ on $(a,b,c,d)$. Thus, $B'$ is full.

\vspace{.2cm}

\begin{figure}
\includegraphics[width= 6cm]{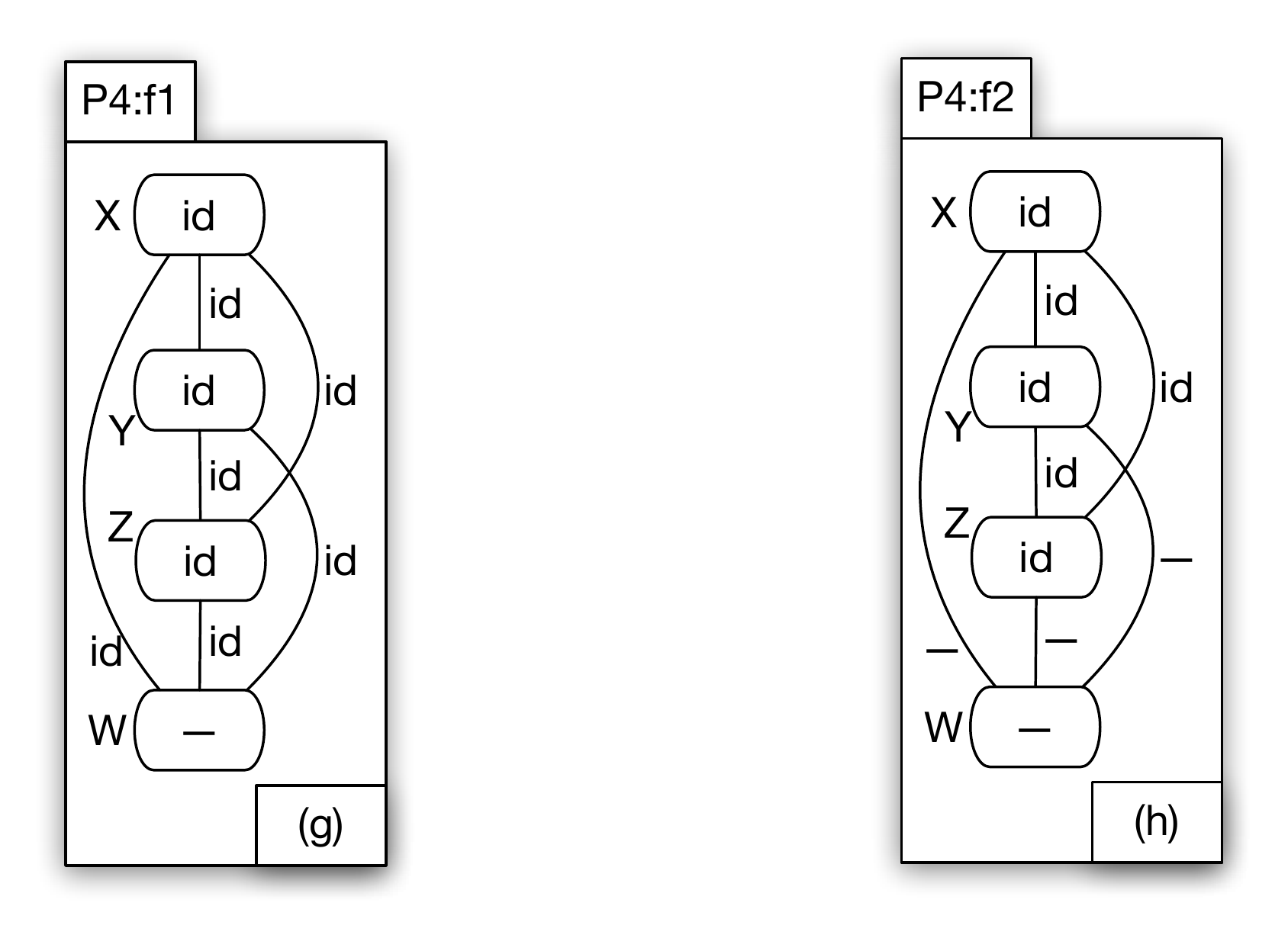}
%\caption{}{}
\caption{Cases C4:f}\label{C4f}
\end{figure}

Figure~\ref{C4f} contains the relevant 4-constellations with three orbits $X>Y>Z>W$. There are only two constellations to be checked as we assume that the behavior is identical on and between $X, Y, Z$, flips the graph relation on $W$, and that there is no full subconstellation. In particular, either the graph relation is kept everywhere else, or the graph relation is flipped everywhere else.

\vspace{.2cm}

C4:f1. %The group $\cl{\apl}$ is forced because of the subconstellation induced by $X\cup W$.
 The constellation induced by $X\cup W$, treated in C2:b3, shows that $f$ generates $\apl$.

C4:f2. %The group $\cl{\apu}$ is forced because of the subconstellation induced by $X\cup W$.
 The constellation induced by $X\cup W$, treated in C2:b4, shows that $f$ generates $\apu$.

\section{The 42 Reducts}
\label{sect:verification}

We have so far established that every group in $\L$ is the join of groups in $\JI$ (Corollary~\ref{cor:join}); in particular, we already know that $\L$ is finite. The goal of this section is to obtain a precise picture of $\L$, and a description of all 42~non-trivial proper reducts.

Since all elements of $\L$ are joins of groups in $\JI$,  we have in particular that the join irreducibles of $\L$ are contained in $\JI$. Assuming that all elements of $\JI$  are also join irreducible in $\L$, which we are going to see in this section, we then get that the elements of $\L$ correspond precisely to the \emph{ideals} of $\JI$.

\begin{definition}
A set $\I\subseteq \JI$ is called an \emph{ideal} iff it contains all elements of $\JI$ contained in ${\bigvee \{H:H\in \I\}}$.
\end{definition}

It is well-known and easy to see that every finite lattice is isomorphic to the lattice of ideals of its join irreducibles. We therefore have the following.

\begin{prop}
The lattice $\L$ of closed supergroups of $\adel$ is isomorphic to the lattice of ideals of $\JI$ with the order of inclusion.
\end{prop}

It will be convenient to represent the ideals of $\JI$ by their maximal elements, which makes the following definition useful.

\begin{definition}
For a set $\S\subseteq\JI$ we write $\dcl \S$ for the \emph{downward closure} of $\S$ in $\JI$, i.e., the subset of those groups in $\JI$ which are contained in some group in $\S$.
\end{definition}

We will now systematically list all ideals of $\JI$ by their maximal elements. We also identify them by their name in Figure~\ref{fig:reducts}.

\subsection{Ideals containing $\Aut(D;E)$, $\Aut(D;<)$, or $\Aut(D;T)$.}

These ideals are easily determined using the existing classifications for the random graph, the rationals, and the random tournament. Recall the names of the elements of $\JI$ defined in Figures~\ref{fig:JI1} and~\ref{fig:JI2}.

\begin{lem}\label{lem:ideals-E}
The ideals of $\JI$ containing $\Aut(D;E)$ are the following:
\begin{enumerate}
\item the trivial ideal $\JI$ (ijk);
\item $\dcl \{\Aut(D;E)\}$ (i);
\item $\dcl \{\Aut(D;E), \cl{\mix -\id},\cl{\mix -\rev}\}$ (dfi);
\item $\dcl \{\Aut(D;E),\cl{ \mix \sw\id}, \cl{\mix\id\turn},\cl{\mix \sw\turn}\}$ (cei);
\item $\dcl \{\Aut(D;E), \cl{\mix -\id},\cl{\mix -\rev},\cl{\mix \sw\id}, \cl{\mix \sw\turn}\}$ (cdefi).
\end{enumerate}
\end{lem}
\begin{proof}
This follows from Thomas' classification~\cite{RandomReducts} of the reducts of $\Aut(D;E)$ and the obvious inclusions which hold between the elements of $\JI$ (e.g., $ \cl{ \mix \id\rev}\subseteq \Aut(D;E)$). The groups in (3) are precisely those elements of $\JI$ which preserve $R^{(4)}$, (4) those which preserve $R^{(3)}$, and (5) those which preserve $R^{(5)}$.
\end{proof}

\begin{lem}\label{lem:ideals-<}
The ideals containing $\Aut(D;<)$ are the following:
\begin{enumerate}
\item the trivial ideal $\JI$ (ijk);
\item $\dcl \{\Aut(D;<)\}$ (j);
\item $\dcl \{\Aut(D;<), \cl{ \mix \id\rev},\cl{\mix -\rev}\}$ (bfj);
\item $\dcl \{\Aut(D;<), \cl{\mix\id\turn},\cl{\mix \sw\turn}\}$ (aej);
\item $\dcl \{\Aut(D;<),\cl{\mix \id\rev},\cl{\mix -\rev},\cl{\mix\id\turn},\cl{ \mix \sw\turn}\}$ (abefj).
\end{enumerate}
\end{lem}
\begin{proof}
This follows from Cameron's classification~\cite{Cameron5} of the reducts of $\Aut(D;<)$  and the obvious inclusions  which hold between the elements of $\JI$ (e.g., $ \cl{ \mix-\id}\subseteq \Aut(D;<)$). The groups in (3) are precisely those elements of $\JI$ which preserve $\Betw$, (4) those which preserve $\Cycl$, and (5) those which preserve $\Sep$.
\end{proof}

\begin{lem}\label{lem:ideals-T}
The ideals containing $\Aut(D;T)$ are the following:
\begin{enumerate}
\item the trivial ideal $\JI$;
\item $\dcl \{\Aut(D;T)\}$ (k);
\item $\dcl \{\Aut(D;T), \mix -\id, \mix \id \rev\}$ (ack);
\item $\dcl \{\Aut(D;T), \mix \sw\id, \mix\id\turn\}$ (bdk);
\item $\dcl \{\Aut(D;T), \mix -\id,\mix \id\rev,\mix \sw\id,\mix\id\turn\}$ (abcdk).
\end{enumerate}
\end{lem}
\begin{proof}
This follows from Bennett's classification~\cite{Bennett-thesis} of the reducts of $\Aut(D;T)$  and the obvious inclusions  which hold between the elements of $\JI$ (e.g., $ \cl{ \mix-\rev}\subseteq \Aut(D;T)$).  The groups in (3) are precisely those elements of $\JI$ which preserve $\Betw_T$, (4) those which preserve $\Cycl_T$, and (5) those which preserve $\Sep_T$.
\end{proof}

\subsection{Ideals containing $\cl{\apl}$ or $\cl{\apu}$.}

Recall the names of the elements of $\JI$ defined in Figures~\ref{fig:JI1} and~\ref{fig:JI2}.

\begin{lem}\label{lem:ideals-andras}
Let $\I$ be an ideal of $\JI$ which contains $\cl{\apl}$ or $\cl{\apu}$. Then one of the following holds:
\begin{enumerate}
\item  $\I=\dcl \{\cl{\apl}\}$ (g);
\item  $\I=\dcl \{\cl{\apu}\}$ (h);
\item  $\I=\dcl\{\cl{\mix \sw\id},\cl{\apl}, \cl{\apu}\}$ (dgh);
\item  $\I=\dcl\{\cl{\mix \id \rev},\cl{\mix \sw\id},\cl{\mix -\rev} ,\cl{\apl}, \cl{\apu}\}$ (adegh);
\item  $\I=\dcl\{\cl{\mix \id \turn},\cl{\mix \sw\id},\cl{\apl}, \cl{\apu}\}$ (bdfgh);
\item  $\I=\dcl\{\cl{\mix \id \rev},\cl{\mix \id \turn},\cl{\mix \sw\id},\cl{\mix -\rev},\cl{\apl}, \cl{\apu}\}$ (abdefgh).
\item $\I$ contains $\Aut(D;<)$.
\end{enumerate}
\end{lem}
\begin{proof}
We first show that $\cl{\apl,\apu}$ contains $\mix \sw\id$. It is clear that the group contains $\mix -\id$. Let $\pi$ be the irrational which defines both $\apl$ and $\apu$, and assume that $\apl$ and $\apu$ send the interval $(-\infty,\pi)$ onto itself. Then $\mix -\id\circ \apl\circ\apu$ switches the graph relation between the intervals $(-\infty,\pi)$ and $(\pi,\infty)$, and keeps it otherwise. From this we see that $\mix \sw\id$ is contained in $\cl{\apl,\apu}$.

If $\I\subseteq \OP$, then this implies that (1), (2), (3), or (7) holds. Otherwise, Propositions~\ref{prop:unorderedCanonical} and~\ref{prop:ordered} imply that $\I$ contains some $G\in \JI\setminus \JIo$. From the obvious inclusions, it follows that $\I$ contains either $\cl{\mix \id\rev}$, $\cl{\mix -\rev}$, $\cl{\mix \id\turn}$ or $\cl{\mix \sw\turn}$. 

Note that $\I$ contains $\cl{\mix \id\rev}$ if and only if it contains $\cl{\mix -\rev}$, since $\cl{\apl}$ and $\cl{\apu}$ both contain $\mix -\id$. Moreover, if $\I$ contains those functions, then it contains both $\cl{\apl}$ and $\cl{\apu}$. To see this, suppose it contains $\cl{\apl}$, and let $\pi$ be the irrational defining $\apl$. Assuming that $\apl$ and $\mix \id\rev$ send the interval $(-\infty,\pi)$ onto itself, we then have $\mix\id\rev\circ \apl\circ \mix \id\rev$ behaves precisely like $\apu$.

If $\I$ contains $\cl{\mix \id\turn}$ or $\cl{\mix \sw\turn}$, then it contains both $\cl{\apl}$ and $\cl{\apu}$: assuming that $\mix \id\turn$ sends the interval $(-\infty,\pi)$ onto itself, we have that $\mix \id\turn\circ \apl\circ \mix\id\turn$ behaves like $\apu$; we obtain the same if we use $\mix \sw\turn$ instead of $\mix \id\turn$. Hence, in this situation $\I$ also contains $\cl{\mix \sw\id}$, by our observation of the first paragraph.
It follows that $\I$ contains $\cl{\mix \id\turn}$ if and only if it contains $\cl{\mix \sw\turn}$, and the lemma is proven.

\ignore{

We first show that $\bigvee \I$ contains an order-preserving function $f$ which for some irrational $\pi$ flips the graph relation between the intervals $(-\infty,\pi)$ and $(\pi,\infty)$, and which keeps the graph relation otherwise. To this end, we may without loss of generality assume that $\I$ contains $\cl{\apl}$. It is clear that $\I$ then also contains $\cl{\mix -\id}$.
In the following, let $\pi$ be the irrational defining all the functions $\apl$, $\mix \id\turn$ and $\mix \sw\turn$. Assuming moreover that $\apl$ sends $(-\infty,\pi)$ onto itself, and that $\mix \id\rev$, $\mix -\rev$, $\mix \id\turn$, and $\mix \sw\turn$ send $(-\infty,\pi)$ onto $(\pi,\infty)$, we then obtain the desired $f$ as $\mix-\id\circ h \circ\apl \circ h\circ \apl$, where $h$ is any of the functions $\mix \id\rev$, $\mix -\rev$, $\mix \id\turn$ or $\mix \sw\turn$.

Now fix $c_1,c_2,c_3\in D$ such that $c_1,c_2\in (-\infty,\pi)$, and $c_3\in (\pi,\infty)$. Then the canonical function $g$ guaranteed by Proposition~\ref{prop:ramseystrong} satisfies the assumptions of Lemma~\ref{lem:2levelsaboveconstant}: it is order-preserving
}
\end{proof}

\subsection{The remaining ideals.}

Finally we consider those ideals of $\JI$ which do not contain $\cl{\apl}$, $\cl{\apu}$, $\Aut(D;<)$, $\Aut(D;E)$, or $\Aut(D;T)$; in other words, those ideals which are subsets of $\{\cl{\mix - \id},\cl{\mix \id \rev},\cl{\mix -\rev},\cl{\mix \sw\id},
\cl{\mix \id \turn},\cl{\mix \sw\turn}\}$. Since they are numerous, we do not list them, but state the only restriction that holds for such subsets.

\begin{lem}\label{lem:ideals-atoms}
Let $\I$ be an ideal of $\JI$. Then:
\begin{itemize}
\item $\I$ contains none, exactly one, or all of $\{\mix -\id, \mix \id\rev, \mix -\rev\}$;
\item $\I$ contains none, exactly one, or all of $\{\mix \sw\id, \mix \id\turn, \mix \sw\turn\}$.
\end{itemize}
\end{lem}
\begin{proof}
This is an easy exercise of composing functions, and we do not give the proof.
\end{proof}

\subsection{The final reduct count.}
Lemmas~\ref{lem:ideals-E}, \ref{lem:ideals-<}, \ref{lem:ideals-T}, \ref{lem:ideals-andras}, and \ref{lem:ideals-atoms} give us 4+4+4+6+(5$\mult 5-1$)=42 possible non-empty proper ideals. 
It remains to show that each of those sets is in fact an ideal. To this end, we list for every possible ideal $\I$ a set of relations which is invariant under the group defined by $\I$ (i.e., the group $\bigvee\I$). The relations show that all those groups are distinct.
We need to define one relation which we have not encountered so far.

For a permutation $\sigma$ of $\{1,2,3,4\}$, 
write $S_4^\sigma$ for the $4$-ary relation 
$$\{(a_1,a_2,a_3,a_4) \; | \; S_4(a_{\sigma(1)},a_{\sigma(2)},a_{\sigma(3)},a_{\sigma(4)})\}.$$ 
Write $D_4$ for the dihedral group on $1,2,3,4$. Set
\begin{align*}
S^D&:=\bigcup_{\sigma \in D_4} S_4^{\sigma}.
\end{align*}
Figure~\ref{fig:pres} shows which ideals preserve which relations (presence of a cross in the table stands for ``preserves'', and the absence for ``violates''; the names of the elements of $\JI$ are defined in Figures~\ref{fig:JI1} and~\ref{fig:JI2}). Checking the table can be automated and is left to the reader. We remark that the groups in $\OP$ correspond precisely to the rows in which the order relation $<$ is preserved.

\subsection{The lattice inclusions.}
We finally describe how we have verified the edges in the Hasse diagram of $\mathfrak L$ in Figure~\ref{fig:reducts}. It is straightforward to verify that if there is an ascending edge from vertex $u$ to vertex $v$ in the diagram, then the group corresponding to $u$ is contained in the group corresponding to $v$: we have labelled the vertices by the maximal elements of ideals of JI, and from those maximal elements one can calculate the ideals by adding $\{a,b\}$ when there is $i$ among the maximal elements,
$\{g,h,d\}$ when there is $j$, $\{c\}$ when there is $g$ or when there is $h$, and $\{f,e\}$ when there is $k$. The containments between elements of JI we use here have been verified previously. 

Next, we have to verify that whenever there is no edge between vertices $u$ and $v$ in the graph then either there is an ascending path in the diagram from one vertex to the other (and hence the edge is not displayed since it is in the transitive closure of the drawn edges of the Hasse diagram), or indeed the group corresponding to $v$ does not contain the group corresponding to $u$. For this task we have used Figure~\ref{fig:pres}; note that for a given $u$ we only have to do the check for the $v$ that are \emph{maximal} with the property above, since then also all groups contained in the group corresponding to $v$ will not contain the group corresponding to $u$ as well (this reduces the work considerably). 

We finally would like to comment on the visible symmetry of the lattice: besides the lattice automorphism that switches $g$ and $h$ and fixes all other elements of the lattice, there is also
a lattice automorphism that acts as the permutation $(ik)(ae)(bf)(j)(g)(h)(c)(d)$ on JI. 
We have not found a concise argument for the existence of this automorphism without explicitly describing $\mathfrak L$.  
\begin{figure}
\begin{center}
{\footnotesize
\begin{tabular}{|l|c|c|c|c|c|c|c|c|c|c|c|c|c|c|c|}\hline
  &                		$E$&$R^{(3)}$&$R^{(4)}$&$R^{(5)}$&$<$&$\betw$&$\cyc$&$\sep$&$T$&$\betw_T$&$\cyc_T$&$\sep_T$&$R_3^l$&
  $R_3^u$&$S^D$\\\hline
%JIs
a				&\bf x&x&x&x		&&\bf x&&x		&&x&&x		&&		&x\\\hline
b				&\bf x&x&x&x		&&&\bf x&x		&&&x&x		&&		&x\\\hline
c				&&&\bf x&x		&\bf x&x&x&x		&&x&&x		&x&x		&x\\\hline
d				&&\bf x&&x		&\bf x&x&x&x		&&&x&x		&&		&x\\\hline
e				&&&\bf x&x		&&\bf x&&x		&\bf x&x&x&x		&&		&x\\\hline
f				&&\bf x&&x		&&&\bf x&x		&\bf x&x&x&x		&&		&x\\\hline
g				&&&&		&x&x&x&x		&&&&		&\bf x&		&x\\\hline
h				&&&&		&x&x&x&x		&&&&		&&\bf x		&x\\\hline
i				&\bf x&x&x&x		&&&&		&&&&		&&		&\\\hline
j				&&&&		&\bf x&x&x&x		&&&&		&&		&\\\hline
k				&&&&		&&&&		&\bf x&x&x&x		&&		&\\\hline\hline
%joins of two JIs
ab				&\bf x&x&x&x		&&&&\bf x		&&&&x		&&		&x\\\hline
ad				&&\bf x&&x		&&\bf x&&x		&&&&x		&&		&x\\\hline
af				&&\bf x&&x		&&&&x		&&\bf x&&x		&&		&x\\\hline
bc				&&&\bf x&x		&&&\bf x&x		&&&&x		&&		&x\\\hline
be				&&&\bf x&x		&&&&x		&&&\bf x&x		&&		&x\\\hline
cd				&&&&\bf x		&\bf x&x&x&x		&&&&x		&&		&x\\\hline
cf				&&&&\bf x		&&&\bf x&x		&&\bf x&&x		&&		&x\\\hline
de				&&&&\bf x		&&\bf x&&x		&&&\bf x&x		&&		&x\\\hline
ef				&&&&\bf x		&&&&x		&\bf x&x&x&x		&&		&x\\\hline
%joins of three JIs
ace				&&&\bf x&x		&&\bf x&&x		&&x&&x		&&		&x\\\hline
bdf				&&\bf x&&x		&&&\bf x&x		&&&x&x		&&		&x\\\hline
abce				&&&\bf x&x		&&&&\bf x		&&&&x		&&		&x\\\hline
acde				&&&&\bf x		&&\bf x&&x		&&&&x		&&		&x\\\hline
acef				&&&&\bf x		&&&&x		&&\bf x&&x		&&		&x\\\hline
abdf				&&\bf x&&x		&&&&\bf x		&&&&x		&&		&x\\\hline
bcdf				&&&&\bf x		&&&\bf x&x		&&&&x		&&		&x\\\hline
bdef				&&&&\bf x		&&&&x		&&&\bf x&x		&&		&x\\\hline
abcdef			&&&&\bf x		&&&&\bf x		&&&&x		&&		&x\\\hline
dgh				&&&&		&\bf x&x&x&x		&&&&		&&		&\bf x\\\hline
% the new ones
adegh			&&&&		&&\bf x&&x		&&&&		&&		&\bf x\\\hline
bdfgh			&&&&		&&&\bf x&x		&&&&		&&		&\bf x\\\hline
abdefgh			&&&&		&&&&\bf x		&&&&		&&		&\bf x\\\hline\hline
% above i
cei				&&&\bf x&x		&&&&		&&&&		&&		&\\\hline
dfi				&&\bf x&&x		&&&&		&&&&		&&		&\\\hline
cdefi				&&&&\bf x		&&&&		&&&&		&&		&\\\hline\hline
% above j neu
aej				&&&&		&&\bf x&&x		&&&&		&&		&\\\hline
bfj				&&&&		&&&\bf x&x		&&&&		&&		&\\\hline
abefj				&&&&		&&&&\bf x		&&&&		&&		&\\\hline\hline
% above k neu
ack				&&&&		&&&&		&&\bf x&&x		&&		&\\\hline
bdk				&&&&		&&&&		&&&\bf x&x		&&		&\\\hline
abcdk			&&&&		&&&&		&&&&\bf x		&&		&\\\hline

\end{tabular}
}
\caption{Preservation table. An $x$ in an entry
indicates that the relation of the column is preserved by the group $G$ indicated by the row. Moreover, %we marked for each group a subset
%of the entries using bold font such that 
each group $G$ equals the set of all permutations that preserve the relations that are marked by $\bf x$ in bold font.}
\label{fig:pres}
\end{center}

\end{figure}

\bibliographystyle{alpha}
\bibliography{orderedrandom.bib}
\end{document}